\documentclass[a4paper,11pt]{amsart}
\usepackage[foot]{amsaddr}
\usepackage[margin=2.7cm]{geometry}
\usepackage{graphicx,color, amsfonts,amsmath}
\usepackage[titletoc]{appendix}
\usepackage{xcolor}
\usepackage{mathtools} 
\usepackage{hyperref}

\usepackage{subfig}
\usepackage{graphicx}

\theoremstyle{definition}

\theoremstyle{remark}

\numberwithin{equation}{section} 

\usepackage{color}
\newcommand{\green}{\color{black}}
\newcommand{\blue}{\color{black}} 

\newcommand{\magenta}{\color{black}}

\newcommand{\ldb}{\mathopen{\lbrack\!\lbrack}}
\newcommand{\rdb}{\mathclose{\rbrack\!\rbrack}}
\begin{document}
\title[DEIM and unfitted mesh FEMs]{Discrete Empirical Interpolation and unfitted mesh FEMs: application in PDE--constrained optimization}

\author{Georgios Katsouleas\textsuperscript{1}}
\address{\textsuperscript{1}Department of Mathematics, National Technical University of Athens, Zografou Campus, 15780, Greece.}
\thanks{}
\email{gekats@mail.ntua.gr}

 \author{Efthymios N. Karatzas\textsuperscript{1,2,3}}
\address{\textsuperscript{2}FORTH, Institute of Applied and Computational Mathematics, Heraclion, Crete, Greece.}
\address{\textsuperscript{3}SISSA (affiliation),  International School for Advanced Studies, Mathematics Area, mathLab, Via Bonomea 265, Trieste, 34136, Italy.}
\email{karmakis@math.ntua.gr \& efthymios.karatzas@sissa.it}

\author{Fotios Travlopanos\textsuperscript{1}}
\thanks{}
\email{ftravlo@gmail.com}

\thanks{}%
\subjclass[2000]{Primary}%
\date{\today}

 \begin{abstract}
In this work, we investigate the performance CutFEM as  a high fidelity solver as well as we construct a competent and economical reduced order solver for PDE--constrained optimization problems in parametrized domains that live in a fixed background geometry and mesh.  Its effectiveness and reliability will be assessed through its application for the numerical solution of quadratic optimization problems with elliptic equations as constraints, examining 
an archetypal case. The reduction strategy will be via Proper Orthogonal Decomposition of suitable FE snapshots, using an aggregated state and adjoint test space, while the efficiency of the offline-online decoupling will be ensured by means of Discrete Empirical Interpolation of the optimality system matrix and right--hand side, enabling thus a rapid resolution of the reduced order model for each new spatial configuration. 
\end{abstract}

\maketitle

\section{Introduction}


Parametrized partial differential equations (PDEs) arise naturally in diverse scientific fields, ranging from engineering to finance. The inclusion of {\itshape physical  parameters} in the formulation of PDEs is necessitated by the need to model possible variations in the systems'  physical properties, boundary conditions or source terms. On the other hand, {\itshape geometric parametrization} is even more challenging and reflects the desire to study phenomena in complex geometrical configurations involving large deformations and even topological changes. Repeated queries to the underlying solver  for different parameter values are typical in the context of uncertainty quantification, optimal control and shape optimization. In such cases, standard discretization techniques become impractical or infeasible.

It should be noted that, even in the non--parametric case, the numerical solution of PDE--constrained optimization is a computationally demanding task, since the optimality conditions require the solution of the state problem, the adjoint problem, as well as a further set of equations ensuring the optimality of the solution. 
{{We highlight that the solution of such systems is not an easy task while a large number of variables originating from the spatial discretization presents an additional challenge that makes the direct solution of the above optimality system difficult. 
{\green{With the technology introduced in this paper, one may efficiently compute approximate solutions of parametrized optimal control problems with PDEs as constraints. Examples where one needs to repeatedly solve parametrized optimal control problems, 
are optimal control for fluid flows or shape optimization, where a number of physical
properties and/or geometrical variables or material parameters are often not precisely known, which challenges
input data under uncertainties,  often approximated with optimal control techniques. 
For instance, it is important in certain problems to estimate
the shape and the location of a system for the production policies to be optimized, see e.g. the
boundary of an oil reservoir, the internal combustion engine and the vibrated geometry of an exhaust system
which need to be optimized so as to maximize the power output of the engine. In several cases, parametrized
geometries may characterize materials in micro-structure and macroscopic level and mechanics of deformable
bodies. Other sources of geometric parametrization may come from manufacturing tolerances which
leads engineering analysis to a more robust design, for instance, the thermostat housing pipe optimal
control problem, also, the distribution of gas within a pipeline may be parametrized and/or random, since a pipeline operator,
typically does not know in advance whether a power plant will come online and for how long. In all such
situations, parametrization in modeling parameters, in geometry, in initial conditions, or in spatially varying
material properties induces uncertainty where parameterization and the finding/recognition of the optimal control may be involved in the outputs of the model and in any quantities of interest
derived from these outputs. For all aforementioned in this paragraph, the interested reader may see for example \cite{AK20} 
and references therein.}}
 {\green{So, we will}} verify the numerical efficiency of the reduced basis method, for the geometrical parametrization bottleneck in FEM when topology deformations are taking place see e.g.  uncertainty quantification,  shape optimization etc, since we avoid remeshing for each different domain configuration as well as transformations to reference domains, underlying that in this framework we avoid the bottleneck of the computationally and time expensive assembling  of the final finite element matrices system.
}}


{\green{Generally, for linear quadratic elliptic optimal control problems we refer  to {\cite{D10}} for reduced basis and a posteriori error estimation on parametrized  cases, to \cite{L1971} for a classical development of a theory of optimal control, as well as, Sections 2.1--2.15  in the recent book {\cite{T2010} for optimal control theory and discretization.}}}

{\blue{Furthermore,}} strategies for quick and efficient resolution of optimal control problems governed by parametrized PDEs in a low--dimensional framework have been developed{\blue{.
For reduced basis methods for the solution of parametrized optimal control problems to the case of (noncoercive/elliptic) Stokes equations, we refer to \cite{NMR15}. For  certified reduced basis method (for reliable solution of parametrized optimal control problems governed by partial differential equation with elliptic equations as constraints and infinite-dimensional control variable) using a Galerkin projection onto a low-dimensional space of basis functions and an efficient and rigorous a posteriori error estimate on the state-control-adjoint variables, we cite \cite{NRMQ13,KG14}.
 In the work of \cite{DH13}, a parameter optimization problem subject to constraints parametrized partial differential equations where the derivative information can be calculated efficiently in the reduced basis framework (in the case of a general linear output functional and parametrized evolution problems with linear parameter separable operators), is considered. In the same work,  the sensitivity information is calculated directly instead of applying the more widely used adjoint approach with rigorous a-posteriori error estimators for the solution, the gradient and the optimal parameters that are computed online.
 Furthermore in \cite{KV08}, the problem of unmodelled dynamics in the proper orthogonal decomposition approach to optimal control is avoided, employing the optimality system proper orthogonal decomposition (OS-POD) approach. A--posteriori error estimates for reduced order models (ROMs) in the parametrized linear--quadratic case related to parabolic partial differential and stationary Helmholtz equation as constraints respectively, are provided in {\blue{\cite{D10, TUV11}}}}}. The usual ingredients of the reduced basis methodology include a Galerkin projection performed onto a suitable low--dimensional subspace of basis functions and then an efficient offline--online decoupling of the computational procedure, due to an affine parametric dependence of the system components. The offline stage is a time--consuming, parameter--independent pre--processing phase, while the online state should involve an inexpensive calculation for each new parameter value.  Although typical for physical parametrizations \cite{Strazzullo18}, affine parametric dependence does not follow readily for geometric parametrizations.

The  overall objective of this manuscript  is to explore {\blue{for the first time to our knowledge}} {\textit{efficient reduction strategies}} for optimal control problems governed by geometrically parametrized PDEs {{{\blue{discretized by {\textit{embedded finite element methods}} in a fixed background geometry and mesh employing the advantages of these methods in cases of geometrical deformations}}}}. 
The combination of  the Proper Orthogonal Decomposition Galerkin (POD) and the Discrete Empirical method (DEIM) procedure with an {\textit{unfitted method}} to alleviate the problem of {\textit{offline/online decoupling}} in the case of  parametrized domains appears beneficial. {\blue{In particular,}} we avoid remeshing for each different domain configuration as well as transformations to reference domains, noting that in this framework we avoid the bottleneck of the computationally and time expensive assembling  of the final finite element matrices system. To the authors' best knowledge, this is the first time this has been attempted in the literature.

More precisely, we focus on the case of quadratic cost functionals constrained to linear PDEs in parametrized domains.  Representing the boundary of the geometry through level--set functions, we follow the strategy in \cite{Hansbo_etc_1} and consider the original physical domain as being embedded in a fixed Cartesian background geometry which encompasses all possible configurations of the domain for different parameter values. Then, using the associated fixed background mesh, it is possible to combine unfitted mesh finite element methods (FEMs)  with reduced order modeling techniques to decrease the overall computational burden of the numerical simulation. 
 
Embedded and immersed methods have a long history, dating back to the pioneering  work of Peskin \cite{P1972}. Several improved variants can be found in the recent literature, including such methods as the ghost-cell finite difference method \cite{WFC13}, cut--cell volume method \cite{PHO16}, immersed interface \cite{KB18}, ghost fluid \cite{BG14},  shifted boundary methods \cite{MS18}, $\phi$--FEM \cite{DL19}. {\blue{Related to the cut finite element methods (CutFEM):  in \cite{Hansbo_etc_1} a review on advances on robust unfitted finite element methods on cut meshes on complex geometries is developed, while in
\cite{BH2010} cut elements based on a stabilized Lagrange multiplier method are examined.
In \cite{BH2012}, the classical Nitsche type weak boundary conditions are  extended to a fictitious domain setting, 
while in \cite{L19}, a new way to impose Dirichlet and Neumann boundary conditions on unfitted meshes without cutting the mesh cell is introduced. 
 For elliptic interface problem, using a Nitsche approach allowing discontinuities internal to the elements, and in the approximation across the interface and a priori error estimates we refer to \cite{HH2002}.
In \cite{AK20}, elliptic optimal control problem defined on uncertain domains are discretized by a fictitious domain finite element method and cut elements while for a properly stabilized CutFEM method and 
 an effective reduced order model for a  embedded boundary parametrized Cahn-Hilliard phase-field system based on cut finite elements, one can see \cite{KR20}}}, among others. Finally for a comprehensive overview of this research area, the interested reader is referred to the review paper \cite{MI05}. 
 
 The strategy of building ROMs for geometrically parametrized  PDEs upon embedded methods has been previously exploited in  the works {\blue{\cite{KBR2019,KSASR2019,KSNSRa2019,KSNSRb2019} (all based on the proper orthogonal decomposition method) considering  
\begin{itemize}
    \item projection-based reduced order models for a cut finite element method in parametrized domains for Poisson and Stokes systems based on a fixed background domain and comparing the efficiency of several solution extensions in a background geometry, 
    \item a reduced order approach for the embedded shifted boundary finite element method and a heat exchange system on parametrized geometries,
    \item a reduced basis approach for partial differential equations on parametrized geometries based on the shifted boundary finite element method and application to Stokes flows, and
    \item a reduced order shifted boundary method for parametrized incompressible Navier-Stokes equations
\end{itemize}
respectively.}} 
This approach avoids remeshing and the need to develop a reference domain formulation, as typically done in fitted mesh finite element formulations \cite{RHP08}.
  
However, the resulting reduced system is non-affine, hence the fast evaluation of a reduced solution for each new parameter value is somewhat restricted, requiring large computational effort. The computational efficiency of this approach requires the further approximation of the non--affine elements of the reduced system (i.e., system matrix and right hand side) by an affine decomposition. The present work aims to address this challenge through the combination of a ROM with the {\textit{discrete empirical interpolation method (DEIM)}} \cite{CS10}. This allows us to recover an offline/online computational procedure. The reduced order model is generated through a Galerkin projection to a reduced basis obtained though a proper orthogonal decomposition (POD) approach (POD--Galerkin), using the CutFEM method as high fidelity solver. 
{\blue{We underline that in the  classical finite element methods literature there are two kinds of DEIM that may employed, one with assembled formulation in the online phase used for the reduced solution and one with unassembled, namely UDEIM, which may
lead to additional gains in the online cost of the reduced
order models although with additional costs in the offline. The latter types of reduced order models depend on significantly  fewer components of the arguments than that resulting from the DEIM applied to the assembled form similar to \cite{TiRi13, AHS14}.
%
%
This modified DEIM / UDEIM algorithm minimizes the number of element function calls, operating with unassembled vectors in order
for each DEIM selected point to be linked with just one element.
%
%
%
In the present work even if the latter appears to be more efficient in the online phase (although it appears more time expensive in the offline) we are investigating the assembled elements case. This is an  important and state of the art research next step and extension of  the works \cite{KBR2019,KSASR2019,KSNSRa2019,KSNSRb2019} that  allows the uncoupling of the offline/online procedures efficiently (and up to our knowledge for the first time in the literature)  considering the more demanding parametrically varying domains and CutFEM (based on fixed background mesh), solving optimal control systems that could not resolved in the past.  A full investigation for the unassembled formulation case will be treated in a future work.  
%
%
%
}}

The paper is structured as follows. In Section \ref{preliminaries}, we introduce our model problem: a linear/quadratic optimal control problem constrained by an  elliptic PDE defined in a general parametrized domain. We recall the necessary background for its discretization via cut elements  and recast  the problem  in the framework of saddle--point problems, whereby the well-posedeness of the FE truth  approximation readily follows. The reduction strategy and the main features of the method are detailed in Section \ref{ROMsection}. In particular, we first describe how  a Garlerkin projection onto a low--dimensional subspace of basis functions built via POD on optimal solutions for all state, adjoint state and control  variables is performed. Implementation details for the DEIM approximation of the system matrix and its right--hand side vector follow, concluding with a summary of all specific aspects  for an efficient offline/online decoupling.   The decreasing computational cost in the proposed procedure is demonstrated in section \ref{numer} via numerical evidence. Moreover, the convergence of the reduced POD-DEIM solution to the CutFEM solution with respect to the number of POD basis functions is verified.

\section{The model  problem and preliminaries}\label{preliminaries}
\subsection{Problem formulation}
 For any parameter  $\boldsymbol\mu \in \mathcal{P}$,  let $\Omega(\boldsymbol\mu)\subset {\mathbb R}^d$ ($d=2 ,3$) be a bounded, $\boldsymbol\mu$--dependent domain of interest and $\Gamma(\boldsymbol\mu)=  \partial \Omega(\boldsymbol\mu)$ its  boundary. 
 This geometric description of the domain $\Omega(\boldsymbol\mu)$  in terms of parameters in the space $\mathcal{P}$ reflects the fact that its exact shape is governed by uncertainty. The state and (distributed) control variables are denoted by $\left(y(\boldsymbol\mu), u(\boldsymbol\mu)\right) \in H^1(\Omega(\boldsymbol\mu)) \times L^2(\Omega(\boldsymbol\mu))$.
We consider a quadratic-linear optimization problem on $\Omega(\boldsymbol\mu)$; namely, for a given target function $y_d(\boldsymbol\mu) \in L^2(\Omega(\boldsymbol\mu))$ and a fixed control stabilization constant $\alpha>0$, we seek to minimize the quadratic cost functional
\begin{equation}\label{1}
J(y(\boldsymbol\mu),u(\boldsymbol\mu))=\frac{1}{2}\int_{\Omega(\boldsymbol\mu)}\left(y(\boldsymbol\mu)-y_d(\boldsymbol\mu)\right)^2dx+\frac{\alpha}{2}\int_{\Omega(\boldsymbol\mu)}u(\boldsymbol\mu)^2dx,
\end{equation}
 subject to the linear state equation:
\begin{eqnarray} \label{The Poisson Problem} 
 -\Delta  y(\boldsymbol\mu)&=&f(\boldsymbol\mu)+u(\boldsymbol\mu)\quad\,\,\,\,\, \text{in $\Omega(\boldsymbol\mu)$},
 \nonumber\\
  y(\boldsymbol\mu) &=&g_D(\boldsymbol\mu) \qquad\quad\, \, \quad \text{on $\Gamma_D(\boldsymbol\mu)$},
  \\
  \nonumber {\bf{n_\Gamma}}\cdot \nabla y(\boldsymbol\mu) &=&g_N(\boldsymbol\mu) \quad\quad\quad\quad\, \, \text{on $\Gamma_N(\boldsymbol\mu)$}.
\end{eqnarray}
In (\ref{The Poisson Problem}), $\Gamma_D(\boldsymbol\mu)$, $\Gamma_N(\boldsymbol\mu)\subset \Gamma(\boldsymbol\mu)$ denote those boundary parts  where Dirichlet and Neumann boundary conditions are applied,  while the forcing term $f(\boldsymbol\mu)$ and the boundary data $g_D(\boldsymbol\mu)$,   $g_N(\boldsymbol\mu)$ are given functions in $\Omega(\boldsymbol\mu)$ and on the boundary parts $\Gamma_D(\boldsymbol\mu)$, $\Gamma_N(\boldsymbol\mu)$,  respectively. {\magenta{We also use the notation $(\cdot, \cdot)_{\Omega(\boldsymbol\mu)}$, $(\cdot, \cdot)_{\Gamma_D(\boldsymbol\mu)}$ and $(\cdot, \cdot)_{\Gamma_N(\boldsymbol\mu)}$ for the $L^2(\Omega(\boldsymbol\mu))$, $L^2({\Gamma_D(\boldsymbol\mu)})$ and $L^2({\Gamma_N(\boldsymbol\mu)})$ inner products respectively.}} The parametric weak formulation of   (\ref{The Poisson Problem})  reads as follows: find a  state 
$$ y(\boldsymbol\mu) \in V_{g_D}(\boldsymbol\mu):=\left\{w(\boldsymbol\mu) \in H^{1}(\Omega(\boldsymbol\mu)): w|_{\Gamma_D(\boldsymbol\mu)}=g_D(\boldsymbol\mu)\right\},
$$ 
such that
\begin{equation}\label{weak}
\left(\nabla y(\boldsymbol\mu), \nabla \upsilon(\boldsymbol\mu)\right)_{\Omega(\boldsymbol\mu)}=\left(f(\boldsymbol\mu) + u(\boldsymbol\mu), \upsilon(\boldsymbol\mu)\right)_{\Omega(\boldsymbol\mu)}+\left(g_N(\boldsymbol\mu), 
\magenta{{\bf{n_\Gamma}}\cdot \nabla \upsilon (\boldsymbol\mu)}\right) 
_{\Gamma_N(\boldsymbol\mu)},
\end{equation}
for all $ \upsilon(\boldsymbol\mu) \in V_0(\boldsymbol\mu)${\magenta{, where  $V_0(\boldsymbol\mu)$ is equal to $V_{g_D(\boldsymbol\mu)}$ for $g_D \equiv 0$.}}

In the following, we will use an \itshape optimize--then--discretize \normalfont approach to approximate the solution of (\ref{1})--(\ref{weak}) in a discrete solution space through an unfitted method. 
{\magenta{Regarding the existence and the uniqueness of an optimal solution of the system (\ref{1})--(\ref{The Poisson Problem}) we refer to \cite{T2010} and references therein.
}}

\subsection{Discretization with Cut Finite Elements} 

Implementation of an embedded method for the discretization of (\ref{weak}) requires the definition of a suitable fixed  background domain $\mathcal B$ which contains all possible  configurations of $\Omega(\boldsymbol\mu )$ for different values of $\boldsymbol\mu \in \mathcal{P}$. Letting $\mathcal B_h$ its corresponding mesh,  the continuous boundary value problem (\ref{weak}) is then discretized on an \itshape extended domain \normalfont  ${\Omega} _{\mathcal{T}}(\boldsymbol\mu)\supset{\Omega}(\boldsymbol\mu)$ corresponding to a minimal submesh $\mathcal{T}_h(\boldsymbol\mu)\subset \mathcal B_h$ which covers the original domain ${\Omega}(\boldsymbol\mu)$, but is not fitted to its boundary. Therefore, discrete solutions will be sought in the piecewise linear finite element space
\begin{equation}\label{space}
V_h(\boldsymbol\mu)=\left \{  \upsilon \in C^0({\bar {\Omega}}_{\mathcal T}(\boldsymbol\mu))\,:\,\upsilon |_{T}  \in P^1(K), \, \forall T\in \mathcal{T}_h(\boldsymbol\mu)    \right \}.
\end{equation}
As usual, the subscript $h=\max_{T \in \mathcal{B}_h} diam(T)$ indicates global mesh size.

According to the CutFEM paradigm, the boundary conditions at $\Gamma(\mu)$ are satisfied weakly through a variant of Nitsche's method. On the other hand, coercivity over the whole computational mesh is ensured by means of additional ghost penalty terms in the discrete counterpart of (\ref{weak}) which act on the gradient jumps in the boundary zone. Therefore, a more delicate analysis of boundary intersecting elements is required; we use
$$
{G}_h(\boldsymbol\mu):=\{T\in \mathcal{T}_h(\boldsymbol\mu):T\cap\Gamma(\boldsymbol\mu)\neq\emptyset\}
$$
to denote  the relevant set of simplices and consider the set of associated element faces by 
\begin{equation}\label{FG}
\mathcal{F}_G(\boldsymbol\mu):=\left\{F: F \ \text{is a face {\magenta{of}}} \ T, F\notin \partial\Omega_{\mathcal{T}}(\boldsymbol\mu)\right\};
 \end{equation}
 The jump of the gradient of $v_h(\boldsymbol\mu)\in V_h(\boldsymbol\mu)$ on a face $F=T\cap T'$ is defined by $$\ldb{{\bf n_F}}\cdot\nabla  v_h(\boldsymbol\mu)\rdb={\bf{n_F}}\cdot\nabla v_h(\boldsymbol\mu) |_T -{\bf{n_F}}\cdot \nabla v_h(\boldsymbol\mu)|_{T'},$$
  where $ { \bf n_T}$ denotes the outward--pointing unit normal vector to $T$. 

Using the previous notation, we actually seek a pair of discrete  control and state ${\magenta{u_h(\boldsymbol\mu), y_h(\boldsymbol\mu)}} \in V_h(\boldsymbol\mu)$, such that  
 \begin{equation}\label{Discr1}
 J(y_h(\boldsymbol\mu),u_h(\boldsymbol\mu))=\frac{1}{2}\int_{\Omega(\boldsymbol\mu)}\left(y_h(\boldsymbol\mu)-y_d(\boldsymbol\mu)\right)^2dx+\frac{\alpha}{2}\int_{\Omega(\boldsymbol\mu)}u_h(\boldsymbol\mu)^2dx
\end{equation}
is minimized under the constraint
 \begin{equation}\label{Discr}
\alpha\left(y_h(\boldsymbol\mu), \upsilon_h(\boldsymbol\mu) \right) =\left(u_h(\boldsymbol\mu), \upsilon_h(\boldsymbol\mu)\right )_{\Omega(\boldsymbol\mu)} + c(\upsilon_h(\boldsymbol\mu)), \ \ \ \forall  \upsilon_h(\boldsymbol\mu) \in V_h(\boldsymbol\mu),
\end{equation}
where the action of the control is separated from the remaining linear forcing term in the right hand side of (\ref{Discr}), {\magenta{while $c(\upsilon_h(\boldsymbol\mu))$ in defined in (\ref{def:(uh)})}}. Here, the bilinear term is defined by
\begin{align}
\alpha\left(y_h(\boldsymbol\mu), \upsilon_h(\boldsymbol\mu) \right) & := \left( \nabla y_h(\boldsymbol\mu) ,\nabla \upsilon_h(\boldsymbol\mu)\right)_{\Omega(\boldsymbol\mu)}-\left({\bf n_\Gamma}\cdot \nabla y_h(\boldsymbol\mu),\upsilon_h(\boldsymbol\mu)\right)_{\Gamma_D(\boldsymbol\mu)}-  \label{bilinear} \\
& \quad -\left(y_h(\boldsymbol\mu),{\bf n_\Gamma} \cdot \nabla\upsilon_h(\boldsymbol\mu)\right)_{\Gamma_D(\boldsymbol\mu)}+\left(\gamma_Dh^{-1}y_h(\boldsymbol\mu),\upsilon_h(\boldsymbol\mu)\right )_{\Gamma_D(\boldsymbol\mu)} + \nonumber \\
&\quad +\left(\gamma_N h {\bf n_\Gamma}\cdot \nabla y_h(\boldsymbol\mu),{\bf n_\Gamma}\cdot \nabla \upsilon_h(\boldsymbol\mu) \right)_{\Gamma _N(\boldsymbol\mu)} + { j(y_h(\boldsymbol\mu),\upsilon_h(\boldsymbol\mu))  } ,\nonumber
\end{align}
 where
 \begin{eqnarray}
  j(y_h(\boldsymbol\mu),\upsilon_h(\boldsymbol\mu)):=\sum_{F\in \mathcal{F}_G(\boldsymbol\mu)}\left(\gamma_1 h\ldb{{\bf n_F}}\cdot\nabla  y_h(\boldsymbol\mu)\rdb,\ldb{{{\bf n_F}}\cdot\nabla} \upsilon_h(\boldsymbol\mu)\rdb\right)_F \label{ghost},
  \end{eqnarray}
{\magenta{the subscript F indicates an $L^2$ inner product on the facet $F$,}}
and $\gamma_D$, $\gamma_N$, and $\gamma_1$  are positive penalty parameters, see for instance \cite{BH2012}. 
The role of the stabilization term $  j(y_h(\boldsymbol\mu),v_h(\boldsymbol\mu))$ is to actually extend  the coercivity of the bilinear form from the physical domain $\Omega(\boldsymbol\mu)$ to the whole computational domain ${\Omega} _{\mathcal{T}}(\boldsymbol\mu)$. Moreover,
\begin{align}
c\left(\upsilon_h(\boldsymbol\mu) \right) & := \left(f_h(\boldsymbol\mu), \upsilon_h(\boldsymbol\mu)\right )_{\Omega(\boldsymbol\mu)} +\left(g_D(\boldsymbol\mu),\gamma_D h^{-1}\upsilon_h(\boldsymbol\mu)+ {\bf n_\Gamma}\cdot \nabla \upsilon_h(\boldsymbol\mu)\right)_{\Gamma _D(\boldsymbol\mu)}+\nonumber \\ & \quad +\left(g_N(\boldsymbol\mu), \upsilon_h(\boldsymbol\mu)+\gamma_N h {\bf n_\Gamma}\cdot \nabla \upsilon_h(\boldsymbol\mu)\right)_{\Gamma _N(\boldsymbol\mu)} .\label{def:(uh)}
\end{align}

Note that the solution of (\ref{Discr1}), (\ref{Discr}) is sought in the $\boldsymbol\mu$--dependent FE space (\ref{space}). Since, we are interested in obtaining and manipulating solutions (snapshots) for different values of $\boldsymbol\mu$, it is instructive to introduce a $\boldsymbol\mu$--independent extension of $V_h(\boldsymbol\mu)$ on the common background domain ${\mathcal B}$ as
$$
\hat{V}_{h}:=\left \{  \upsilon \in C^0({\bar {\mathcal B_h}})\,:\,\upsilon |_K  \in P^1(K), \, \forall K\in \mathcal{B}_h    \right \}.
$$
Hence, the resulting solutions are extended from ${\Omega}_{\mathcal T}(\boldsymbol\mu)$ to the whole (common) background mesh $\mathcal{B}_h$ and have identical degrees of freedom. 
{\blue{
Related to the aforementioned extension, as we will see in Section \ref{ROMsection}, the construction of the reduced order basis is based on the whole background domain.
For this reason, the manipulation of the out of interest -outside- the true geometry area, namely
``ghost area", needs particular care.  In the present work, following \cite{KBR2019}, we
use the solution values as they are computed applying the cut finite element method and the smooth
solution from the true to the extended domain. This allows a smooth extension of the solution to
the neighboring ghost elements with values which are decreasing smoothly to zero, 
This approach guarantees a regular ``solution" in the background
domain based on $\Omega_T(\boldsymbol\mu)$  and permits the construction of a reduced basis with better approximation properties. For
more details and a full investigation of the possible choices and of the handling of the ghost area we also refer to
\cite{KBR2019}.
}}

{\magenta{For the interested reader, theoretical convergence results and a priori error estimates for a  CutFEM fine model are available in \cite{AK20} and \cite{T2010}, see also references therein.}}
\subsection{Saddle--point formulation.} The above optimal control problem may equivalently be formulated as a saddle--point problem \cite{QMN16}. 
{\magenta{Let $N^j$, $j = 1,...,\mathcal{N}$, be a numbering of the mesh nodes. There exists a unique family $(\phi_
j)_{i=1,...,\mathcal{N}}$ such that $ \phi_j \in \hat{V}_h$ and $\phi_j(N^j) = \delta_{i j}$. This family
is a basis of $\in \hat{V}_h$, which is of dimension $\mathcal{N}$, and for all $v_h \in \hat{V}_h$, we have $v_h =\sum^{\mathcal{N}}_{j=1}v_h(N^j)\phi_j = \sum^{\mathcal{N}}_{j=1}a_j\phi_j$, where $a_j=v_h(N^j)$ are appropriate coefficients.}} Indeed, denoting $\hat{V}_{h}=\text{span}\left\{\phi_j\right\}_{j=1}^{\mathcal{N}}$ the finite element basis functions with $\mathcal{N} \equiv \text{dim} \hat{V}_{h}$, then (\ref{Discr1}), (\ref{Discr}) are rewritten at an algebraic level as follows:
$$
\min_{\mathbf{y}_{\boldsymbol\mu}, \mathbf{u}_{\boldsymbol\mu} \in \mathbb{R}^{\mathcal{N}}}J(\mathbf{y}_{\boldsymbol\mu}, \mathbf{u}_{\boldsymbol\mu})=\min_{\mathbf{y}_{\boldsymbol\mu}, \mathbf{u}_{\boldsymbol\mu} \in \mathbb{R}^{\mathcal{N}}}\left[\frac{1}{2}\mathbf{y}_{\boldsymbol\mu}^T\mathbf{M}_{\boldsymbol\mu}\mathbf{y}_{\boldsymbol\mu}+\frac{1}{2}\left\|y_d\right\|_{L^2(\Omega(\boldsymbol\mu))}^2-\mathbf{y}^T_{\boldsymbol\mu}\mathbf{b}_{\boldsymbol\mu} +\frac{\alpha}{2} \mathbf{u}^T_{\boldsymbol\mu}\mathbf{M}_{\boldsymbol\mu}\mathbf{u}_{\boldsymbol\mu}\right]
$$
under the constraint
$$
\mathbf{A}_{\boldsymbol\mu}\mathbf{y}_{\boldsymbol\mu}=\mathbf{M}_{\boldsymbol\mu}\mathbf{u}_{\boldsymbol\mu}+\mathbf{c}_{\boldsymbol\mu},
$$
where $\mathbf{y}_{\boldsymbol\mu}, \mathbf{u}_{\boldsymbol\mu} \in \mathbb{R}^{\mathcal{N}}$ are the vectors of coefficients in the expansions of $\hat{y}_h^{\boldsymbol\mu}, \hat{u}_h^{\boldsymbol\mu} \in \hat{V}_h$ with respect to the finite element basis respectively, 
$\mathbf{b}_{\boldsymbol\mu}\equiv \left(\left(\phi_j, y_d\right)_{\Omega(\boldsymbol\mu)} \right)_{j=1}^{\mathcal{N}}\in \mathbb{R}^{ \mathcal{N}}$,  $\mathbf{c}_{\boldsymbol\mu}\equiv \left(c\left(\phi_j\right)\right)_{j=1}^{\mathcal{N}}\in \mathbb{R}^{ \mathcal{N}}$, and  
\begin{align*}
\mathbf{A}_{\boldsymbol\mu}&:=\left(\mathbf{A(\boldsymbol\mu)}_{ij}\right)_{i,j=1}^{\mathcal{N}}=\left(\alpha\left(\phi_j, \phi_i \right) \right)_{i,j}^{\mathcal{N}}\in \mathbb{R}^{\mathcal{N}\times \mathcal{N}}, \nonumber \\
\mathbf{M}_{\boldsymbol\mu}&:=\left(\mathbf{M(\boldsymbol\mu)}_{ij}\right)_{i,j=1}^{\mathcal{N}}=\left(\left(\phi_j, \phi_i\right)_{\Omega(\boldsymbol\mu)} \right)_{i,j=1}^{\mathcal{N}}\in \mathbb{R}^{\mathcal{N}\times \mathcal{N}}
\end{align*}
are the \itshape stiffness \normalfont and \itshape mass \normalfont matrices. Here the subscripts are meant to imply that  the respective quantities are $\boldsymbol\mu$--dependent.

Considering the stationary conditions of the Lagrangian
\begin{eqnarray*}
\mathcal{L}(\mathbf{y}_{\boldsymbol\mu}, \mathbf{u}_{\boldsymbol\mu}, \mathbf{p}_{\boldsymbol\mu})=\frac{1}{2}\mathbf{y}^T_{\boldsymbol\mu}\mathbf{M}_{\boldsymbol\mu}\mathbf{y}_{\boldsymbol\mu}+\frac{1}{2}\left\|y_d\right\|_{L^2(\Omega(\boldsymbol\mu))}^2-\mathbf{y}^T_{\boldsymbol\mu}\mathbf{b}_{\boldsymbol\mu} +\frac{\alpha}{2} \mathbf{u}^T_{\boldsymbol\mu}\mathbf{M}_{\boldsymbol\mu}\mathbf{u}_{\mu}
\\
{\magenta{ + \mathbf{ \mathbf{p}}^T_{\boldsymbol\mu}\left[\mathbf{A}_{\boldsymbol\mu}\mathbf{y} -
\mathbf{M}_{\boldsymbol\mu}\mathbf{u}_{\boldsymbol\mu}-\mathbf{c}_{\boldsymbol\mu}\right],}}
\end{eqnarray*}
where $ \mathbf{p}_{\mu} \in \mathbb{R}^{\mathcal{N}}$ the vector of Lagrange multipliers, we obtain the $3\mathcal{N} \times 3\mathcal{N}$ linear system
\begin{equation}\label{system}
\underbrace{\begin{bmatrix} \mathbf{M}_{\boldsymbol\mu} & \mathbb{O} & \mathbf{A}^T_{\boldsymbol\mu} \\   \mathbb{O}  &\alpha \mathbf{M}_{\boldsymbol\mu} & -\mathbf{M}^T_{\boldsymbol\mu}  \\ \mathbf{A}_{\boldsymbol\mu} & -\mathbf{M}_{\boldsymbol\mu} &  \mathbb{O} 
\end{bmatrix}}_{\mathbf{\mathcal{A}}(\boldsymbol\mu)}\begin{bmatrix} \mathbf{y}_{\boldsymbol\mu} \\  \mathbf{u}_{\boldsymbol\mu} \\ \mathbf{ \mathbf{p}}_{\boldsymbol\mu}
\end{bmatrix}=\underbrace{\begin{bmatrix} \mathbf{b}_{\boldsymbol\mu} \\ 0  \\ \mathbf{c}_{\boldsymbol\mu}
\end{bmatrix}}_{\beta(\boldsymbol\mu)}
\end{equation}
which clearly has saddle--point structure. Hence, 
letting {\magenta{$\Phi=\begin{bmatrix}\phi_1 & \cdots & \phi_{\mathcal{N}}\end{bmatrix}$
,
the triple $\left(y_h(\boldsymbol\mu), u_h(\boldsymbol\mu), p_h(\boldsymbol\mu)\right)=
\left(\sum_{j=1}^{\mathcal{N}}{y}_{\boldsymbol\mu, j}\phi_j, \sum_{j=1}^{\mathcal{N}}{u}_{\boldsymbol\mu, j}\phi_j , \sum_{j=1}^{\mathcal{N}}{p}_{\boldsymbol\mu, j}\phi_j\right)
$}} is a high--fidelity FE approximation of the optimal state, control and adjoint state respectively, {\magenta{or with simpler notation}} 
$\left(y_h(\boldsymbol\mu), u_h(\boldsymbol\mu), p_h(\boldsymbol\mu)\right)=
{{\left(\Phi\mathbf{y}_{\boldsymbol\mu}, \Phi\mathbf{u}_{\boldsymbol\mu}, \Phi\mathbf{ \mathbf{p}}_{\boldsymbol\mu}\right)}}
.$

\section{The reduced basis method for the parametrized model problem}\label{ROMsection}

System (\ref{system})  typically features a very large size and a process for its solution entails huge computational costs. This is especially true in a many--query context, when optimization is desirable for many different configurations of the domain $\Omega(\boldsymbol\mu)$ of interest; i.e., for many different parameter values of $\boldsymbol\mu \in \mathcal{P}$. In such cases, reduction of the computational complexity is mandatory  for a fast resolution of the problem and for achieving real--time control. 

In this direction, developing a reduced order model (ROM) of much smaller dimension would be highly desirable. In its essence, the procedure to obtain such a ROM involves performing a Galerkin projection on an $N$--dimensional approximation space, where $N$ is typically much smaller than the original dimension $\mathcal{N}$ of $V_h(\boldsymbol\mu)$. The smaller the magnitude of $N$ is, the cheaper the ROM is to solve. 

\subsection{Reduction strategy.}\label{red} RB spaces are built upon high--fidelity FE approximations for the state, adjoint and control variables from (\ref{system}) and may be generated either by performing \itshape Proper Orthogonal Decomposition \normalfont  (POD) on suitably selected snapshots or by a greedy procedure with respect to an a--posteriori error bound. For our purposes, we pursue the former approach and proceed in this section to highlight its most important aspects.

For RB generation according to the POD strategy, one needs to  explore the solution manifold via (\ref{system}) for different domain configurations and compress the resulting data, retaining the most important and essential information. More precisely, one needs to:
 
 \begin{enumerate}
 \item  Sample the parameter space and select
 \begin{equation}\label{train}
 \mathcal S := \left\{\boldsymbol\mu_1, \dots, \boldsymbol\mu_M\right\}\subset \mathcal{P}
 \end{equation}
 of cardinality $M:=\left|\mathcal S\right|$.  
 \item  Compute a set of corresponding \itshape snapshots\normalfont, i.e., FE solutions 
 $$\left\{y_h(\boldsymbol\mu_j), u_h(\boldsymbol\mu_j), p_h(\boldsymbol\mu_j)\right\}_{j=1}^{M}=\left\{\Phi\mathbf{y}_{\boldsymbol\mu}, \Phi\mathbf{u}_{\boldsymbol\mu}, \Phi\mathbf{ \mathbf{p}}_{\boldsymbol\mu}\right\}_{j=1}^{M}$$
  of (\ref{system}), realized  for the parameter values sampled in $\mathcal S$. These are then stored in the $\mathcal N \times M$ \itshape snapshot matrices  \normalfont $\mathbf S_{y} :=\begin{bmatrix}\mathbf{y}_{\boldsymbol\mu_1} & \cdots & \mathbf{y}_{\boldsymbol\mu_M}\end{bmatrix}$,  $\mathbf S_{u} :=\begin{bmatrix}\mathbf{u}_{\boldsymbol\mu_1} & \cdots & \mathbf{u}_{\boldsymbol\mu_M}\end{bmatrix}$ and  $\mathbf S_{p} :=\begin{bmatrix}\mathbf{ \mathbf{p}}_{\boldsymbol\mu_1} & \cdots & \mathbf{ \mathbf{p}}_{\boldsymbol\mu_M}\end{bmatrix}$ and define the related correlation matrices $\mathbf C_{j} := \frac{1}{M}\mathbf S_{j}^T \mathbf M\mathbf S_{j}$ ($j=y, u, p$).
 \item Solve the $M \times M$ eigenvalue problems
 $$
 \mathbf C_{j} x_{ji}=\lambda_{ji}x_{ji},  \ \ i=1, \dots ,M \ \text{and} \ \ j=y, u, p,
 $$
 where each set of eigenvalues $\left\{\lambda_{ji}\right\}_{i=1}^{M}$  is indexed in non--increasing order. Truncating the full orthogonal bases $\left\{x_{ji}\right\}_{i=1}^{M}$  and retaining the first $N_j$ elements in each case ($j=y, u, p$), we define 
 $$
 V_j:=\frac{1}{\sqrt{M}}\Phi\mathbf S_{j}\begin{bmatrix}x_{j1} & \cdots & x_{jN_j}\end{bmatrix} \in \mathbb{R}^{\mathcal{N} \times N_j}, \ \ \text{for} \ \  j=y, u, p.
 $$ 
 \item Following the methodology in {\blue{\cite{D10, KG14, NRMQ13}}}, 
 the state, control and adjoint variables are approximated in lower-dimensional spaces substituting
 $$
\mathbf{y}_{\boldsymbol\mu} \simeq V_{yp}\mathbf{y}_{\boldsymbol\mu}^N, \ \ \ \ \mathbf{u}_{\boldsymbol\mu} \simeq V_{u}\mathbf{u}_{\boldsymbol\mu}^N, \ \ \ \  \mathbf{p}_{\boldsymbol\mu} \simeq V_{yp}\mathbf{p}_{\boldsymbol\mu}^N,
 $$
 where the range of $V_{yp} : =\begin{bmatrix}V_y & V_p\end{bmatrix}\in \mathbb{R}^{\mathcal{N} \times (Ny+N_p)}$ is an \itshape aggregated space \normalfont generated  \itshape both \normalfont by state and adjoint snapshots, ensuring the stability and the well--posedeness of the Reduced Basis approximation \cite[Lemma 3.1]{NRMQ13}.  
{\blue{
We note that the states and adjoints are approximated by the same space spanned by the columns of $V_{yp}$, whereas the control is approximated by the space $V_u$. Using the same $V_{yp}$ for states and adjoints implies that the system (\ref{ROM}) inherits the same properties as the system (\ref{system}) (in particular $A_\mu$ and $V_{yp}^T A_\mu V_{yp}$ are invertible). Subsequently, for the problem (\ref{Discr1})-(\ref{Discr}), the optimal control is a multiple of the optimal adjoint, {\green{thus, the optimal control $u$ will be approximated numerically by the same space as the adjoint variable $p$ 
(see \cite{convectiveOptimal_Hinze17,T2010} and references therein) namely $V_{yp}= V_u$}}. 
}}
Enforcing the orthogonality of the residual of (\ref{system}) to  the range of
 \begin{equation}\label{RB}
 V:=V_{yp}\oplus V_{u}\oplus V_{yp} \in \mathbb{R}^{\mathcal{N} \times [2(N_y+N_p)+N_u]},
 \end{equation}
 we project the original $\mathcal{N} \times\mathcal{N}$ model to the dense but reduced $[2(N_y+N_p)+N_u] \times [2(N_y+N_p)+N_u]$ system
 \begin{equation}\label{ROM}
\underbrace{\begin{bmatrix}V_{yp}^T \mathbf{M}_{\boldsymbol\mu}V_{yp} & \mathbb{O} & V_{yp}^T\mathbf{A}^T_{\boldsymbol\mu}V_{yp} \\   \mathbb{O}  &\alpha V_{u}^T\mathbf{M}_{\boldsymbol\mu}V_{u} & -V_{u}^T\mathbf{M}^T_{\boldsymbol\mu}V_{yp}  \\ V_{yp}^T\mathbf{A}_{\mu}V_{yp} & -V_{yp}^T\mathbf{M}_{\boldsymbol\mu}V_{u} &  \mathbb{O} 
\end{bmatrix}}_{\mathbf{\mathcal{A}}_N(\boldsymbol\mu):=V^T\mathbf{\mathcal{A}}(\boldsymbol\mu)V}\begin{bmatrix} \mathbf{y}_{\boldsymbol\mu}^N \\  \mathbf{u}_{\boldsymbol\mu}^N \\ \mathbf{ \mathbf{p}}_{\boldsymbol\mu}^N
\end{bmatrix}=\underbrace{\begin{bmatrix} V_{yp}^T\mathbf{b}_{\boldsymbol\mu} \\ 0  \\ V_{yp}^T\mathbf{c}_{\boldsymbol\mu}
\end{bmatrix}}_{\beta_N(\boldsymbol\mu):=V^T\beta(\boldsymbol\mu)}.
 \end{equation} 
 \end{enumerate}

A complete exploration of all permissible domain configurations $\Omega(\boldsymbol\mu)$ leads potentially to a cardinality $M$ much larger than the original dimension $\mathcal{N}$, in which case it is more practical to consider in step (3) the alternative eigenvalue problems for  $\tilde{ \mathbf C_{j}}:=\frac{1}{M}\mathbf M^{1/2}\mathbf S_{j}\mathbf S_{j}^T\mathbf M^{1/2} \in \mathbb{R}^{\mathcal{N} \times \mathcal{N}}$ ($j=y, u, p$). In any case, the above procedure incurs high offline computational costs, requiring as it does the solution of three large and dense eigenvalue problems. {\magenta{We remark that  
one might also think of utilizing a singular value decomposition.}}

The previous discussion already reveals the comparative advantage of embedded methods for the numerical approximation of the optimal control problem in (\ref{1})--(\ref{The Poisson Problem}). Indeed, standard FEM would necessitate remeshing, while RB generation would involve algebraic manipulations of snapshots defined on different spatial configurations, which is typically handled via change of variables and introduction of a parameter--independent reference domain \cite{RHP08}. Within the framework of cut finite elements, all this is deemed unnecessary, since all snapshots are already computed on a common background mesh.
{\magenta{Although, we refer to the recent work \cite{convectiveOptimal_Hinze17} on POD for adaptive finite element meshes where also snapshots with different sizes of their finite element coefficient vectors arise.}}

\subsection{Discrete Empirical Interpolation for $\mathbf{\mathcal{A}}_N(\boldsymbol\mu)$ and $\beta_N(\boldsymbol\mu)$.}\label{DEIMsection} An essential feature of an efficient reduced order model should be its rapid online resolution for each new shape of the domain. Even though  (\ref{ROM}) is a  system of order $[2(N_y+N_p)+N_u]$,  the online construction of $\mathbf{\mathcal{A}}_N(\boldsymbol\mu)$ and  $\beta_N(\boldsymbol\mu)$ still depends directly on the dimension $\mathcal{N}$ of the original full--order model, since all four components $\mathbf{A}_{\boldsymbol\mu}, \mathbf{M}_{\boldsymbol\mu} \in \mathbb{R}^{\mathcal{N} \times \mathcal{N}}$ and $\mathbf{b}_{\boldsymbol\mu}, \mathbf{c}_{\boldsymbol\mu}\in \mathbb{R}^{\mathcal{N}}$ would have to {\magenta{be} }assembled. To ensure the efficiency of the reduced model (\ref{ROM}), it is therefore imperative to reduce its computational complexity further, by approximating each  of the components $\mathbf{A}_{\boldsymbol\mu}$, $\mathbf{M}_{\boldsymbol\mu}$, $\mathbf{b}_{\boldsymbol\mu}$ and $\mathbf{c}_{\boldsymbol\mu}$ via projection onto a suitable low-dimensional subspace spanned by an $\boldsymbol\mu$-independent basis. This procedure can be performed by means of the \itshape Discrete Empirical Interpolation Method \normalfont (DEIM) \cite{CS10}. To cast all four cases into a unified framework, we also convert the matrices  $\mathbf{A}_{\boldsymbol\mu}$, $\mathbf{M}_{\boldsymbol\mu} \in \mathbb{R}^{\mathcal{N} \times \mathcal{N}}$ into vectors by stacking their columns on top of one another and describe the procedure to derive an approximate expression for the \itshape vectorization  \normalfont
$$\text{vec}(\mathbf{A}_{\mu})=\begin{bmatrix}{\mathbf{A}_{\boldsymbol\mu}^1}^T  & \cdots & {\mathbf{A}_{\boldsymbol\mu}^{\mathcal{N}}}^T\end{bmatrix} ^T \in \mathbb{R}^{\mathcal{N}^2}$$
 of $\mathbf{A}_{\boldsymbol\mu}=\begin{bmatrix}\mathbf{A}_{\boldsymbol\mu}^1 & \cdots & \mathbf{A}_{\boldsymbol\mu}^{\mathcal{N}}\end{bmatrix} \in \mathbb{R}^{\mathcal{N} \times \mathcal{N}}$. The three remaining cases of $\mathbf{M}_{\boldsymbol\mu}$, $\mathbf{b}_{\boldsymbol\mu}$ and $\mathbf{c}_{\boldsymbol\mu}$ are then treated in a similar fashion. 
Recall from \cite{CS10} that the DEIM approximation $\widehat{\text{vec}(\mathbf{A}_{\boldsymbol\mu})}$ of $\text{vec}(\mathbf{A}_{\boldsymbol\mu})$ is obtained by
\begin{equation}\label{DEIM}
\text{vec}(\mathbf{A}_{\boldsymbol\mu})\simeq \widehat{\text{vec}(\mathbf{A}_{\boldsymbol\mu})}\equiv U_{\mathbf{A}}(P_{\mathbf{A}}^TU_{\mathbf{A}})^{-1}P_{\mathbf{A}}^T\text{vec}(\mathbf{A}_{\boldsymbol\mu}),
\end{equation}
which is efficiently computable in two stages (offline/online).
Here, $U_{\mathbf{A}}$ is the basis matrix for DEIM approximation and $P_{\mathbf{A}}$ is the row selection matrix corresponding to the DEIM interpolation indices. More details about the construction of these quantities and  the decoupling of the  procedure in offline/online stages are in order. A first offline stage consists of the following:
\begin{itemize}
\item [(i.)] A \textit{training procedure}, in which the full-order stiffness matrix $\mathbf{A}_{\boldsymbol\mu}$ is assembled for a training set of $M$ randomly chosen parameter values of $\boldsymbol\mu$ in the parameter space $\mathcal{P}$.
\item [(ii.)] A \textit{POD procedure} performed on  the ensemble of the whole full order snapshots
$$\mathbf{S}_{\text{vec}(\mathbf{A})}:=\begin{bmatrix}\text{vec}(\mathbf{A}_{\boldsymbol\mu_1})& \dots & \text{vec}(\mathbf{A}_{\boldsymbol\mu_M})\end{bmatrix} \in \mathbb{R}^{\mathcal{N}^2\times M}$$ 
of the stiffness matrix for the training set. This results in the definition of an ($m_{\mathbf{A}}$-dimensional and $\boldsymbol\mu$--independent) POD basis $U_{\mathbf{A}}\in \mathbb{R}^{\mathcal{N}^2 \times m_{\mathbf{A}}}$ with $m_{\mathbf{A}}<\mathcal{N}^2$. The number $m_{\mathbf{A}}$ is simply the number of interpolation indices used in the  approximation of the system matrix $\mathbf{A}_{\boldsymbol\mu}$ and will be referred to as  \itshape dimension of DEIM\normalfont.
\item [(iii.)] A \textit{DEIM procedure} for the selection of DEIM interpolation indices. This is an inductive greedy procedure on the input basis $U_{\mathbf{A}}$ from the previous step. The process starts from selecting the first interpolation index $p_1 \in \left\{1,2,\dots, \mathcal{N}^2\right\}$ which corresponds to that entry of the first input basis vector $u_1$ which has the largest magnitude. Then, we initialize a partial matrix $U_{\mathbf{A},1}=u_1$ by specifying the first column of $U_{\mathbf{A}}$ and $P_{\mathbf{A},1}=e_{p_1}$ where $e_{p_{1}}=\left(0,\dots, \underbrace{1}_{\text{$p_{1}$-th entry}}, 0, \dots, 0\right)^T$ denotes the standard basis vector of $\mathbb{R}^{\mathcal{N}^2}$. The remaining interpolation indices $p_{\ell}$ $(\ell =2, \dots, m_{\mathbf{A}})$ are selected so that each of them corresponds to the largest magnitude entry of the residual $r_{\ell}=u_{\ell}-U_{\mathbf{A},(\ell-1)}c_{\ell}$, with $c_{\ell} \in \mathbb{R}^{\ell-1}$  the solution to the linear system $(P_{\mathbf{A}, (\ell-1)}^TU_{\mathbf{A}, (\ell-1)})c_{\ell}=P^T_{\mathbf{A}, (\ell-1)}u_{\ell}$, and then  we inductively augment $U_{\mathbf{A},\ell}=\begin{bmatrix}U_{\mathbf{A}, (\ell-1)} & u_{\ell}\end{bmatrix} \in \mathbb{R}^{\mathcal{N}^2 \times \ell}$; i.e., the first $\ell$ columns of $U_{\mathbf{A}}$. During this stage, the matrices $P_{\mathbf{A}}=\begin{bmatrix}e_{p_1} & \cdots & e_{p_{m_{\mathbf{A}}}}\end{bmatrix} \in \mathbb{R}^{\mathcal{N} \times m_{\mathbf{A}}}$ and $U_{\mathbf{A}}(P^T_{\mathbf{A}}U_{\mathbf{A}})^{-1} \in \mathbb{R}^{\mathcal{N}^2 \times m_{\mathbf{A}}}$ are computed and stored.
\end{itemize} 

Then, during an online stage, an approximation for the stiffness matrix $\mathbf{A}_{\boldsymbol\mu}$ is recovered for each new parameter value $\boldsymbol\mu$ by the previously mentioned relation (\ref{DEIM}). 
Indeed, having stored the $\boldsymbol\mu$--independent first factor $U_{\mathbf{A}}(P_{\mathbf{A}}^TU_{\mathbf{A}})^{-1}\in \mathbb{R}^{\mathcal{N}^2 \times m_{\mathbf{A}}}$ in the offline stage, it clearly suffices to compute only the $m_{\mathbf{A}}$ components of  $\text{vec}(\mathbf{A}_{\boldsymbol\mu})$ in the reduced vector $\tilde{\theta}_{\mathbf{A}}(\boldsymbol\mu)\equiv P^T_{\mathbf{A}}\text{vec}(\mathbf{A}_{\boldsymbol\mu}) \in \mathbb{R}^{m_{\mathbf{A}}}$, instead of the full $\mathcal{N}^2$-dimensional vector $\text{vec}(\mathbf{A}_{\boldsymbol\mu})$. From an implementation viewpoint, a final objective in the offline stage emerges:
\begin{itemize}
\item [(iv.)] Detection of a suitable {\itshape reduced mesh}, whereupon it should be sufficient to partially assemble the stiffness matrix $\mathbf{A}_{\boldsymbol\mu}$ online to correctly retrieve the necessary entries in  $\tilde{\theta}_{\mathbf{A}}(\boldsymbol\mu)\in \mathbb{R}^{m_{\mathbf{A}}}$  for each new choice of parameter $\boldsymbol\mu \in \mathcal{P}$.
\end{itemize}
 
Regarding the latter task, it is necessary to interrelate the selected indices $p_{\ell} \in \left\{1,2,\dots, \mathcal{N}^2\right\}$ ($\ell=1, \dots , m_{\mathbf{A}}$) to standard degrees of freedom in $ \left\{1,2,\dots, \mathcal{N}\right\}$ and corresponding stiffness matrix entries $\mathbf{A}_{\mu}(i_{\ell},j_{\ell})$ as follows:
\begin{equation}\label{pell}
 p_{\ell}=\mathcal{N}(i_{\ell}-1)+j_{\ell}.
 \end{equation}
 Hence, denoting $\text{DOF}_{\tau}$ the set that includes the degrees of freedom corresponding to element $\mathcal{E}_{\tau}$ ($\tau=1, \dots ,\text{\#e}$), where  $\text{\#e}$ denotes the number of elements in the mesh, an element $\mathcal{E}_{\tau}$ should participate in the reduced mesh if and only if the set $\text{DOF}_{\tau}$ includes both indices $i_{\ell},j_{\ell}$ in (\ref{pell})  for some $\ell=1, \dots ,m_{\mathbf{A}}$.

Returning to the online stage,  contributions to the online approximation time include operations required to
\begin{itemize}
\item  Assemble the reduced bilinear form: $\tilde{\theta}_{\mathbf{A}}(\boldsymbol\mu)=P^T_{\mathbf{A}}\text{vec}(\mathbf{A}_{\mu}) \in \mathbb{R}^{m_{\mathbf{A}}}${\magenta{;}}
\item Perform the multiplication in (\ref{DEIM}): $\underbrace{U_{\mathbf{A}}(P_{\mathbf{A}}^TU_{\mathbf{A}})^{-1}}_{\mathcal{N}^2\times m_{\mathbf{A}}: \ \text{offline}}\underbrace{\tilde{\theta}_{\mathbf{A}}(\boldsymbol\mu)
}_{m_{\mathbf{A}} \times 1 : \ \text{online}}{\magenta{;}}$
\item Reshape the result ($\mathcal{N}^2$-dimensional vector) in $\mathcal{N} \times \mathcal{N}$ matrix form.
\end{itemize}

At this point, note an additional complication implied by the specific CutFEM bilinear form (\ref{bilinear}):
the diffusion and the Nitsche terms certainly pose no additional challenges. However, the final ghost penalty term in (\ref{ghost}) acts on the jumps of the gradients over the specific boundary zone element facets $\mathcal{F}_G(\boldsymbol\mu)$ in (\ref{FG}). Such jump--terms require integrators calculating element matrices for both elements at a given facet at once.  Hence, for each of the marked elements in the boundary zone, it is necessary to detect the corresponding facets upon which ghost penalty should be applied. Then, care should be taken to include neighboring elements to these facets (in case these neighboring elements have not already been marked)  in the reduced mesh as well. 
  
Now, rewriting (\ref{DEIM}) as
\begin{equation}\label{over}
\text{vec}(\mathbf{A}_{\boldsymbol\mu})\simeq \text{vec}(\widehat{\mathbf{A}}_{\boldsymbol\mu}):=U_{\mathbf{A}}\theta_{\mathbf{A}}(\boldsymbol\mu),
\end{equation}
with $U_{\mathbf{A}}:=\begin{bmatrix}u_{\mathbf{A}}^1& \dots & u_{\mathbf{A}}^{m_{\mathbf{A}}}\end{bmatrix} \in \mathbb{R}^{\mathcal{N}^2\times m_{\mathbf{A}}}$ and $\theta_A(\boldsymbol\mu):=\left(\theta_{{\mathbf{A}}}^j(\boldsymbol\mu)\right)_{j=1}^{m_{\mathbf{A}}} \in \mathbb{R}^{m_{\mathbf{A}}}$ the corresponding coefficient vector, partitioning
{\magenta{
and letting $U_{\mathbf{A}}^j
\in \mathbb{R}^{\mathcal{N}^2}$ for $j=1, \dots , m_{\mathbf{A}}$}}, the expansion $\text{vec}(\widehat{\mathbf{A}_{\boldsymbol\mu}})=\sum_{j=1}^{m_\mathbf{A}}\theta_{{\mathbf{A}}}^j(\boldsymbol\mu)\text{vec}\left(U_{\mathbf{A}}^j\right)$ is immediate. Hence, going back to the original matrix $\mathbf{A}_{\boldsymbol\mu}$, we conclude
$$
\mathbf{A}_{\boldsymbol\mu}\simeq \sum_{j=1}^{m_\mathbf{A}}\theta_{{\mathbf{A}}}^j(\boldsymbol\mu)U_{\mathbf{A}}^j.
$$

Repeating the procedure for the remaining components, we obtain analogue approximate decompositions:
$$
\mathbf{M}_{\boldsymbol\mu}\simeq\sum_{j=1}^{m_\mathbf{M}}\theta_{{\mathbf{M}}}^j(\boldsymbol\mu)U_{\mathbf{M}}^j, \ \  \mathbf{b}_{\boldsymbol\mu}\simeq\sum_{j=1}^{m_\mathbf{b}}\theta_{{\mathbf{b}}}^j(\boldsymbol\mu)u_{\mathbf{b}}^j, \ \ \mathbf{c}_{\boldsymbol\mu}\simeq\sum_{j=1}^{m_\mathbf{c}}\theta_{{\mathbf{c}}}^j(\boldsymbol\mu)u_{\mathbf{c}}^j, 
$$
with $\left\{U_{\mathbf{M}}^j\right\}_{j=1}^{m_\mathbf{A}} \subset \mathbb{R}^{\mathcal{N}  \times \mathcal{N}}$, $\left\{u_{\mathbf{b}}^j\right\}_{j=1}^{m_\mathbf{b}}\subset \mathbb{R}^{\mathcal{N}}, \left\{u_{\mathbf{c}}^j\right\}_{j=1}^{m_\mathbf{c}} \subset \mathbb{R}^{\mathcal{N}}$.

\subsection{Offline -- Online decoupling.}\label{online} The previous analysis shows that the formation and resolution of the reduced order model (\ref{ROM}) can be decoupled in two stages. Indeed, combination of all previous approximations permits us to approximate the system matrix $\mathbf{\mathcal{A}}(\boldsymbol\mu)$ and the right--hand side $\beta(\boldsymbol\mu)$ of (\ref{system}) as sums of products between given $\boldsymbol\mu$-dependent functions and $\boldsymbol\mu$-independent forms, i.e., $$\mathbf{\mathcal{A}}(\boldsymbol\mu)\simeq\sum_{j=1}^{Q_{\mathcal{A}}}\Theta^{\mathcal{A}}_{j}(\boldsymbol\mu)\mathbf{\mathcal{A}}_j \ \ \ \text{and} \  \ \ \beta(\boldsymbol\mu)\simeq\sum_{j=1}^{Q_{\mathcal{\beta}}}\Theta^{\beta}_{j}(\boldsymbol\mu)\beta_j
$$
with $\mathbf{\mathcal{A}}_j \in \mathbb{R}^{\mathcal{N}\times \mathcal{N}}$, $\beta_j \in \mathbb{R}^{\mathcal{N}}$ and $\left\{\Theta^{q}_{j}\right\}_{j=1}^{Q_q}$ ($q=\mathcal{A}, \beta$)  real functions with values determined during the online DEIM phase. The numbers of terms in the expansions are clearly given by  $$Q_{\mathcal{A}}:=4m_\mathbf{M}+2m_\mathbf{A}, \ \ \  Q_{\mathcal{\beta}}:=m_\mathbf{b}+m_\mathbf{c}. 
$$
For instance, for the expansion of the right hand--side vector, let 
$$
\left\{\beta_j\right\}_{j=1}^{Q_{\beta}}:=\left\{\begin{pmatrix} u_{\mathbf{b}}^j\\ 0_{2 \mathcal{N}}\end{pmatrix}\right\}_{j=1}^{m_\mathbf{b}} \cup \left\{\begin{pmatrix} 0_{2 \mathcal{N}} \\ u_{\mathbf{c}}^j\end{pmatrix}\right\}_{j=1}^{m_\mathbf{c}}
$$ 
and 
\begin{equation}\label{coef}\small
\Theta_j^{\beta}(\boldsymbol\mu)=\theta_{{\mathbf{b}}}^j(\boldsymbol\mu) \  (j=1, \dots ,m_{{\mathbf{b}}}),  \  \ \Theta_j^{\beta}(\boldsymbol\mu)=\theta_{{\mathbf{c}}}^{j-m_{{\mathbf{b}}}}(\boldsymbol\mu) \ \ (j=m_{{\mathbf{b}}}+1, \dots ,m_{{\mathbf{b}}}+m_{{\mathbf{c}}}).
\end{equation}\normalsize
The expansion of the system matrix $\mathbf{\mathcal{A}}(\mu)$ is completely analogous. This variable separation property is clearly inherited by the related quantities 
\begin{equation}\label{affine}\small
\mathbf{\mathcal{A}}_N(\boldsymbol\mu):=V^T\mathbf{\mathcal{A}}(\boldsymbol\mu)V\simeq\sum_{j=1}^{Q_{\mathcal{A}}}\Theta^{\mathcal{A}}_{j}(\boldsymbol\mu)\underbrace{V^T\mathbf{\mathcal{A}}_jV}_{\mathbf{\mathcal{A}}_j^N} \ \ \ \text{and} \ \ \ \beta_N(\boldsymbol\mu):=V^T\beta(\boldsymbol\mu)\simeq\sum_{j=1}^{Q_{\beta}}\Theta^{\beta}_{j}(\boldsymbol\mu)\underbrace{V^T\beta_j}_{\beta_j^N}
\end{equation}\normalsize
 in the reduced model (\ref{ROM}), where, for instance, the $[2(N_y+N_p)+N_u]$--dimensional vectors 
 $$
 \left\{\beta_j^N\right\}_{j=1}^{Q_{\beta}}:=\left\{\begin{pmatrix} V_{yp}^Tu_{\mathbf{b}}^j\\ 0_{N_y+N_p+N_u}\end{pmatrix}\right\}_{j=1}^{m_\mathbf{b}} \cup \left\{\begin{pmatrix} 0_{N_y+N_p+N_u} \\ V_{yp}^Tu_{\mathbf{c}}^j\end{pmatrix}\right\}_{j=1}^{m_\mathbf{c}}
 $$ 
 and the $[2(N_y+N_p)+N_u] \times [2(N_y+N_p)+N_u]$ matrices $\left\{\mathbf{\mathcal{A}}_j^N\right\}_{j=1}^{Q_{\mathcal{A}}}$ are $\mu$--independent. The latter are  similarly  defined in terms of $\left\{U_{\mathbf{A}}^j\right\}_{j=1}^{m_\mathbf{A}}$, $\left\{U_{\mathbf{M}}^j\right\}_{j=1}^{m_\mathbf{M}}$ and the diagonal blocks of $ V=V_{yp}\oplus V_{u}\oplus V_{yp}$.
 
Summarizing all previous considerations, in an offline stage we compute and store:
\begin{enumerate}
\item  the RB matrix $V$ in (\ref{RB}),
\item the matrices $\left\{\mathbf{\mathcal{A}}_j^N\right\}_{j=1}^{Q_{\mathcal{A}}}$ and vectors  $\left\{\beta_j^N\right\}_{j=1}^{Q_{\beta}}$,
\item the row selection matrices $P_{q} \in \mathbb{R}^{\mathcal{N}^2 \times m_q}$ ($q=\mathbf{A}, \mathbf{M}$), $P_{q} \in \mathbb{R}^{\mathcal{N} \times m_q}$ ($q= \mathbf{b}, \mathbf{c}$), and
\item the   matrices $U_q\left(P^T_{q}U_q\right)^{-1} \in \mathbb{R}^{\mathcal{N}^2 \times m_q}$ ($\mathcal{N}^2 \times m_q$--dimensional for $q=\mathbf{A}, \mathbf{M}$ and $\mathcal{N} \times m_q$--dimensional for $q=\mathbf{b}, \mathbf{c}$)
\end{enumerate}
 
Afterwards, in the online stage and for each new $\boldsymbol\mu \in \mathcal{P}$, these precomputed quantities are utilized to assemble the full matrix $\mathbf{\mathcal{A}}_N(\boldsymbol\mu)$ and vector $\beta_N(\boldsymbol\mu)$ from (\ref{affine}) with operations fully independent of $\mathcal{N}$. Indeed, the only ingredients missing to complete the computations in (\ref{affine}) are $\left\{\Theta_{j}^{\mathcal{A}}(\boldsymbol\mu)\right\}_{j=1}^{Q_{\mathcal{A}}}$ and $\left\{\Theta_{j}^{\beta}(\boldsymbol\mu)\right\}_{j=1}^{Q_{\beta}}$. These quantities are directly related to the respective coefficient vectors $\theta_q(\boldsymbol\mu) \in \mathbb{R}^{m_q}$ ($q=\mathbf{A}, \mathbf{M}, \mathbf{b}, \mathbf{c}$); recall, for instance, the interplay for the coefficients for the computation of $\beta_N(\boldsymbol\mu)$ in (\ref{coef}). As explained in the previous subsection,  $\theta_q(\boldsymbol\mu)$ may be obtained for different values of $\boldsymbol\mu$ simply by interpolation of the full--dimensional vectors $\text{vec}(\mathbf{A}_{\boldsymbol\mu}), \text{vec}(\mathbf{M}_{\boldsymbol\mu}) \in \mathbb{R}^{\mathcal{N}^2}$, $\mathbf{b}_{\boldsymbol\mu}, \mathbf{c}_{\boldsymbol\mu} \in \mathbb{R}^{\mathcal{N}}$ on the $m_q$ ($q=\mathbf{A}, \mathbf{M}, \mathbf{b}, \mathbf{c}$) indices selected by the DEIM procedure as \small
$$
\theta_q(\boldsymbol\mu)= P^T_q\text{vec}(q_{\boldsymbol\mu}), \ \text{for} \ q=\mathbf{A}, \mathbf{M} \ \ \text{or} \ \  \theta_q(\boldsymbol\mu)=P^T_qq_{\boldsymbol\mu}, \ \text{for} \ q=\mathbf{b}, \mathbf{c}.
$$
\normalsize
  
 \section{Numerical Experiments}\label{numer}

In this section, we  illustrate the effectiveness of the proposed framework. The following numerical simulations have been implemented in a python environment, using the open--source Netgen/NGSolve finite element software and were performed on an 1.7 GHz Intel Core i7 processor with 8 GB of RAM.
 
Henceforth, for the sake of simplicity, we focus on \itshape Dirichlet \normalfont boundary conditions; i.e., we consider the formulation (\ref{1}) -- (\ref{The Poisson Problem}) with $g_d \equiv g$, $\Gamma_N=\emptyset$ and Tikhonov regularization parameter $\alpha=0.0001$. The Nitsche stabilization and ghost penalty parameters in (\ref{Discr}) are taken to be equal to $\gamma_D=10$ and $\gamma_1=0.1$, respectively.


In particular, we consider   the parametrized square domain $\Omega(\boldsymbol\mu)$ centered at $(1,1)$ with side length $2\times \boldsymbol\mu$ ranging in the interval $[0.8,1]$. 
{\green{Specifically, the domain is $\Omega(\mu) = {(x, y) : \phi(x, y) < 0}$, where
\begin{equation}\label{level1}
\phi(x, y) = \left|x-1\right|+\left|y-1\right|+\left|\left|x-1\right|-\left|y-1\right|\right|-2\boldsymbol\mu,
\end{equation}
and $\mu \in [0.4,0.5]$.}}
 Following the CutFEM discretization paradigm for the high fidelity solver, we embed $\Omega(\boldsymbol\mu)$ in the  background mesh $\mathcal{B}_h=[-0.3,2.3]\times [-0.3,2.3]$ for all $\boldsymbol\mu\in [0.4, 0.5]$. Additionally,  we consider force $f(x,y)=xy$, homogeneous Dirichlet conditions, desired state $y_{d}=\frac{1}{2\pi}\sin(\pi x)\cos(\pi x)$, and distributed control  throughout  the parametrized by  (\ref{level1}) square domain  $\Omega(\boldsymbol\mu)$. 

The discretization parameter for the high---fidelity solver is set $h=0.09$, resulting in a background mesh with 1944 elements and a system with $\mathcal{N}=1031$ degrees of freedom for linear finite elements; hence, a full order system of order $3\times \mathcal{N}=3093$.   
Henceforth, the \itshape dimension of the POD \normalfont refers to the dimension of the POD basis used for projecting the full--order model; the \itshape dimension of DEIM \normalfont refers to the number of interpolation indices used in the  approximation of the nonaffine terms $\mathbf{A}_{\boldsymbol\mu}$,  $\mathbf{M}_{\boldsymbol\mu}$,  $\mathbf{b}_{\boldsymbol\mu}$, and  $\mathbf{c}_{\boldsymbol\mu}$.  Since, to the authors' best knowledge, this is the first time that the DEIM procedure has been applied to the matrices arising from a CutFEM discretization, before tackling the  full-scale optimal control problem, we first discuss the application of the DEIM procedure on the nonaffine terms $\mathbf{A}_{\boldsymbol\mu}$,  $\mathbf{M}_{\boldsymbol\mu}$,  $\mathbf{b}_{\boldsymbol\mu}$, and  $\mathbf{c}_{\boldsymbol\mu}$ in more detail.

\subsection{DEIM approximation of the stiffness matrix $\mathbf{A}_{\boldsymbol\mu}$} \label{DEIMAnew}

\begin{figure}
\centering
\subfloat[$m_{\mathbf{A}}=5$ modes.]{%
\resizebox*{4cm}{!}{\includegraphics{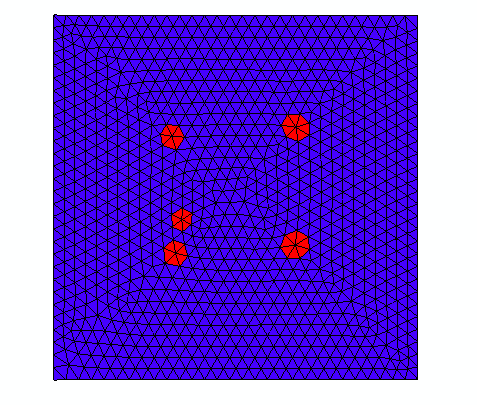}}}\hspace{5pt}
\subfloat[$m_{\mathbf{A}}=10$ modes]{%
\resizebox*{4.15cm}{!}{\includegraphics{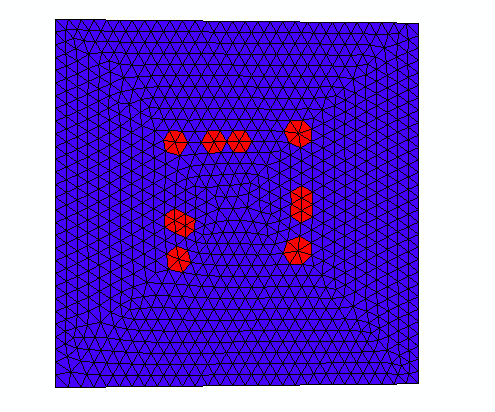}}}
\subfloat[$m_{\mathbf{A}}=15$ modes.]{%
\resizebox*{4cm}{!}{\includegraphics{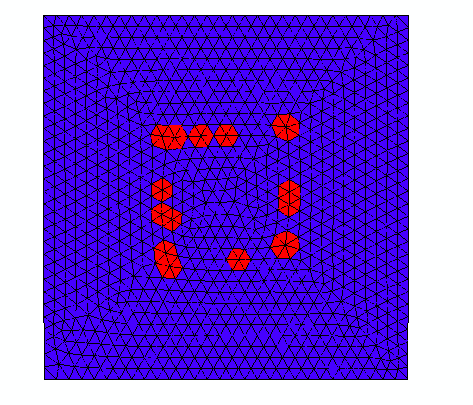}}}
\\
\centering
\subfloat[$m_{\mathbf{A}}=20$ modes.]{%
\resizebox*{4cm}{!}{\includegraphics{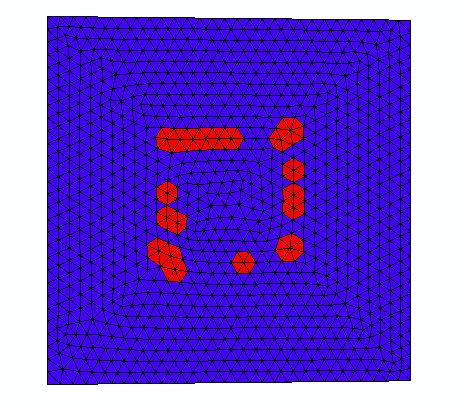}}}\hspace{5pt}
\subfloat[$m_{\mathbf{A}}=25$ modes]{%
\resizebox*{4.15cm}{!}{\includegraphics{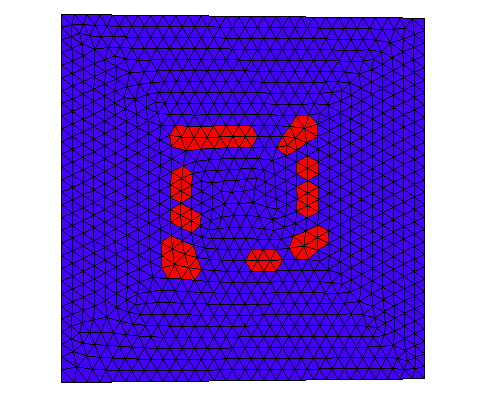}}}
\subfloat[$m_{\mathbf{A}}=30$ modes.]{%
\resizebox*{4cm}{!}{\includegraphics{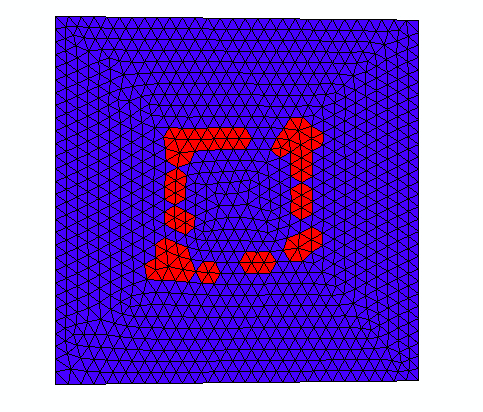}}}
\\
\centering
\subfloat[$m_{\mathbf{A}}=40$ modes.]{%
\resizebox*{4cm}{!}{\includegraphics{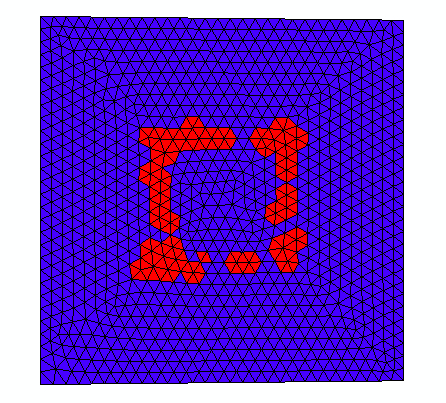}}}\hspace{5pt}
\subfloat[$m_{\mathbf{A}}=60$ modes]{%
\resizebox*{4.1cm}{!}{\includegraphics{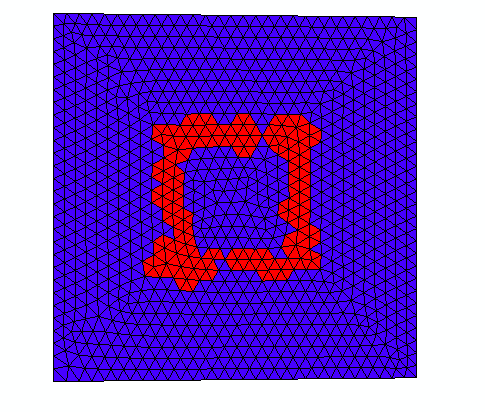}}}
\subfloat[$m_{\mathbf{A}}=90$ modes.]{%
\resizebox*{4cm}{!}{\includegraphics{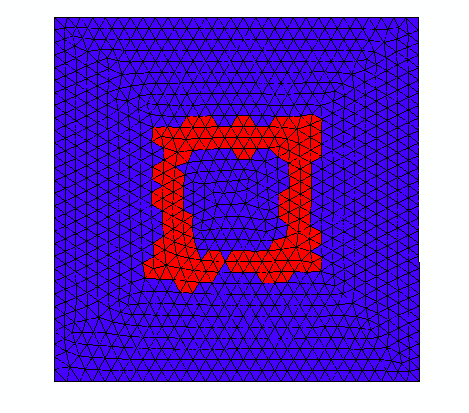}}}
\caption{Illustrations of reduced meshes for the DEIM approximation of the stiffness matrix $\mathbf{A}_{\boldsymbol\mu}$ with progressively increasing DEIM dimension $m_{\mathbf{A}}$ obtained though an offline stage with $M=370$ {\green{training snapshots  onto  parameter samples chosen randomly inside the parameter space $
 [0.4,0.5]$}}.}\label{A_reduced_mesh}
\end{figure}

As described in Section \ref{DEIMsection}, the offline stage of the DEIM procedure for the  stiffness matrix $\mathbf{A}_{\boldsymbol\mu}$ results in the determination of  $m_{\mathbf{A}}$  interpolation indices, which indicate its most important entries. These indices are \itshape global\normalfont, in the sense that these positions are \itshape independent \normalfont of any particular parameter value $\boldsymbol\mu \in \mathcal{P}$; the respective entries will be collected online in the vector $\tilde{\theta}_{\mathbf{A}}(\boldsymbol\mu)= P^T_{\mathbf{A}}\text{vec}(\mathbf{A}_{\boldsymbol\mu})\in\mathbb{R}^{  m_{\mathbf{A}}}$ to recover the  DEIM approximation $\widehat{\text{vec}(\mathbf{A}_{\boldsymbol\mu})}$ to the vectorized form of  $\mathbf{A}_{\boldsymbol\mu}$ in (\ref{DEIM}). Hence, there is no need to assemble the bilinear form $\alpha(y_h(\boldsymbol\mu), \upsilon_h(\boldsymbol\mu))$ on the whole background mesh, but only on the \itshape reduced mesh\normalfont, a judiciously chosen subset, so that only those $m_{\mathbf{A}}$ entries of $\mathbf{A}_{\boldsymbol\mu}$ selected by the row selection matrix $P^T_{\mathbf{A}}\in \mathbb{R}^{m_{\mathbf{A}} \times \mathcal{N}^2}$ in $\tilde{\theta}_{\mathbf{A}}(\boldsymbol\mu)$  are retrieved. Reduced meshes for progressively increasing DEIM dimension  $m_{\mathbf{A}}$ are depicted in Figure \ref{A_reduced_mesh}. Observe that in all cases, most of the information encoded in $\mathbf{A}_{\boldsymbol\mu}$ is captured by   entries which correspond to degrees of freedom in the band $\left\{\Gamma(\boldsymbol\mu):\boldsymbol\mu \in \mathcal{P}\right\}$ of boundary variation for different parameter values $\boldsymbol\mu \in \mathcal{P}=[0.4,0.5]$. This is to be expected, reflecting the fact that the bilinear form $\alpha(y_h(\boldsymbol\mu), \upsilon_h(\boldsymbol\mu))$ is augmented with ghost penalty terms, which are applied on facets intersecting the boundary zone.  On the other hand, degrees of freedom in the interior of $\Omega(\boldsymbol\mu)$ do not seem as important for the DEIM approximation of $\mathbf{A}_{\boldsymbol\mu}$.

 From another viewpoint, Figure \ref{A_reduced_mesh}  offers the opportunity to visualize the impact of increasing DEIM dimension $m_{\mathbf{A}}$ on the respective computational cost to partially assemble $\alpha(y_h(\boldsymbol\mu), \upsilon_h(\boldsymbol\mu))$  on gradually less compact reduced meshes to recover $\tilde{\theta}_{\mathbf{A}}(\boldsymbol\mu)$. A closer inspection reveals another opportunity that should be taken into consideration. For low DEIM dimension $m_{\mathbf{A}}=5$, the reduced mesh simply involves  an equal number of  disjoint patches on the mesh; the
 degree of freedom related to the central node in each patch clearly corresponds to each of the $m_{\mathbf{A}}=5$ entries of $\tilde{\theta}_{\mathbf{A}}(\boldsymbol\mu)$ that need be computed.  For larger DEIM dimension, for instance, for $m_{\mathbf{A}}=10$, these umbrella--like patches begin to intersect. Hence, to fully take advantage of the potential computational savings of the partial online assemble on the reduced mesh, care must be taken to avoid superfluous  calculations for element matrices that are repeated for different entries in $\tilde{\theta}_{\mathbf{A}}(\boldsymbol\mu)$. Hence, in our implementation, in the online stage, all necessary element matrices are computed once and then summed up accordingly to yield  the components of  $\tilde{\theta}_{\mathbf{A}}(\boldsymbol\mu)$.

\begin{figure}
\centering
\subfloat[Reduced mesh characteristics.]{%
\resizebox*{5.5cm}{!}{\includegraphics{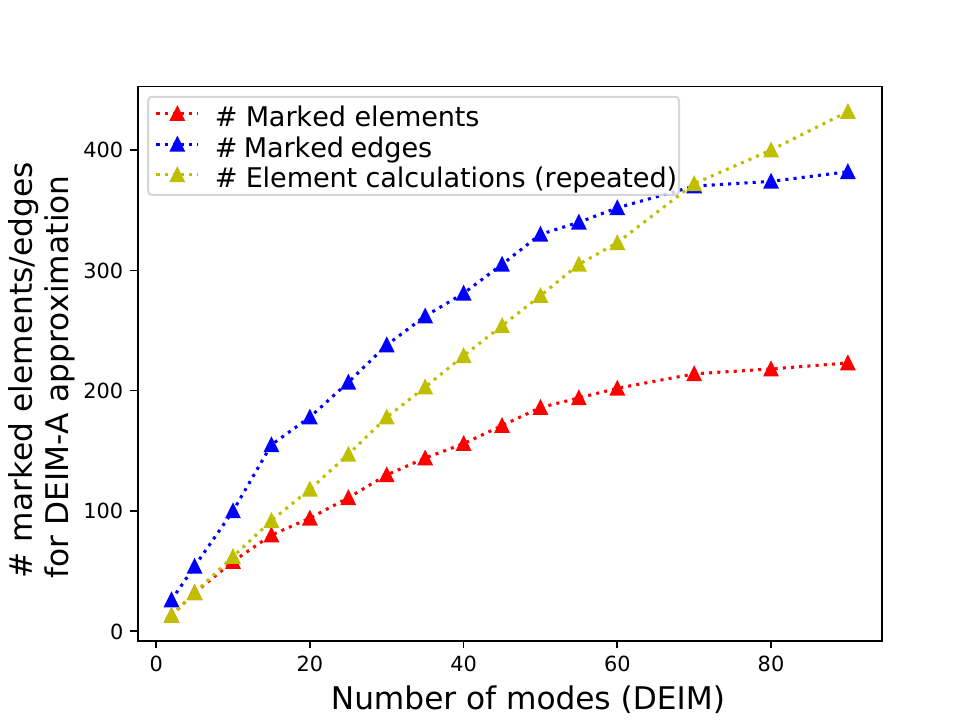}}}\hspace{5pt}  
\subfloat[Online reduced assemble timing.]{%
\resizebox*{5.5cm}{!}{\includegraphics{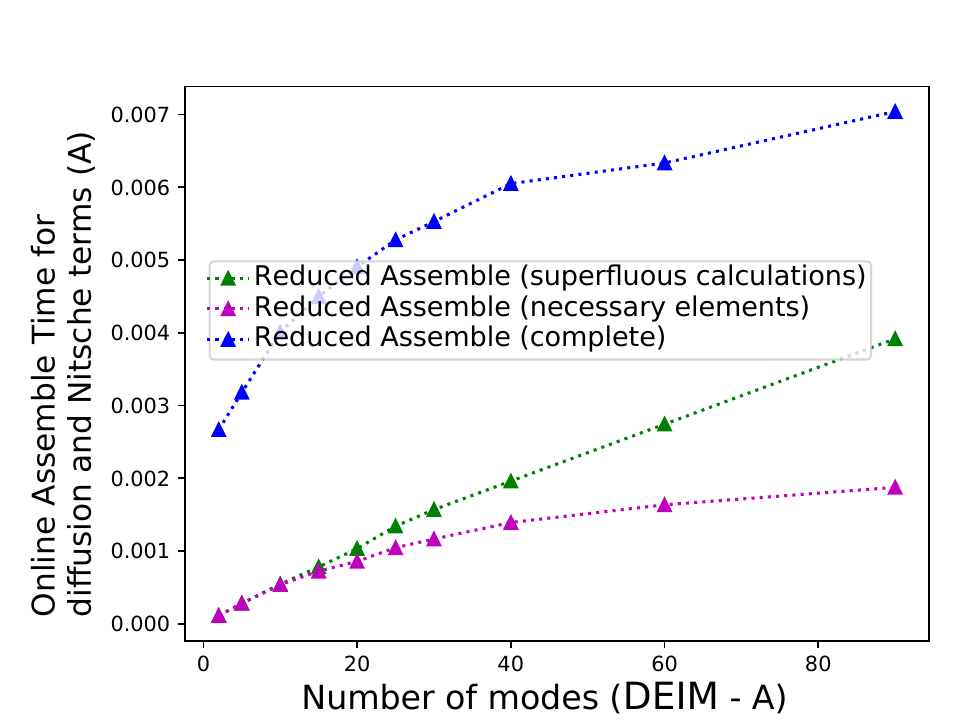}}} \hspace{5pt}
\subfloat[Reduced assemble speedup.]{%
\resizebox*{5.5cm}{!}{\includegraphics{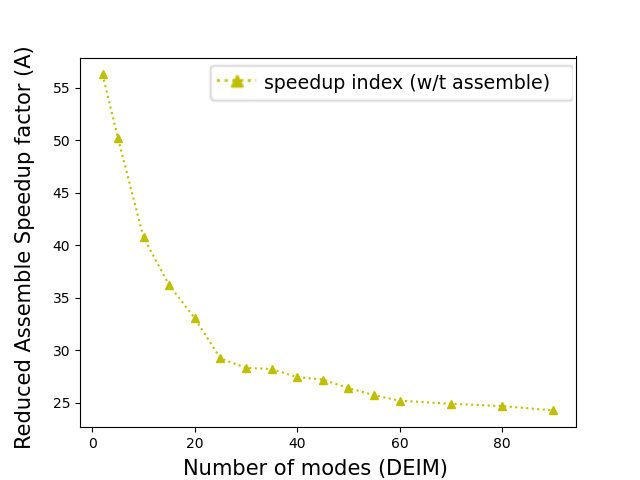}}}\hspace{5pt}
\caption{Reduced mesh characteristics  for progressively increasing DEIM dimension $m_{\mathbf{A}}$ and  their effect on computational performance for the online reduced assemble of $\tilde{\theta}_{\mathbf{A}}(\boldsymbol\mu_0)\in\mathbb{R}^{  m_{\mathbf{A}}}$ for the parameter value $\boldsymbol\mu_0=0.45$.} 
\label{redmeshA2}\label{redmeshA2:1}\label{redmeshA2:2}\label{redmeshA2:3}
\end{figure}

The gradual expansion of the reduced mesh is made more precise by the indicators of  reduced mesh compactness in Table \ref{redmeshA}. The second column describes the number of elements the reduced mesh consists of. Hence, this column also indicates the minimum number of necessary element matrix calculations to correctly assemble the diffusion and Nitsche term components of the bilinear form  $\alpha(y_h(\boldsymbol\mu), \upsilon_h(\boldsymbol\mu))$ on the reduced mesh. The additional computational burden of a more naive implementation involving superfluous element matrix calculations in patch intersections is reflected in the third column of Table \ref{redmeshA}. 
Since most elements in the reduced mesh are clustered around the band of boundary variation $\left\{\Gamma(\boldsymbol\mu):\boldsymbol\mu \in \mathcal{P}\right\}$, patch intersections appear rather early and become more prominent for larger DEIM dimensions; indeed, as verified in Figure \ref{redmeshA2:2} (compare the green and magenta curves), this can be a cause of significant delays and should be avoided.

Regarding the final component of $\alpha(y_h(\boldsymbol\mu), \upsilon_h(\boldsymbol\mu))$, namely the ghost penalty terms, the respective online integrations  should be restricted exclusively on those edges of the reduced mesh that are actually intersected by the boundary $\Gamma(\boldsymbol\mu)$ for each new parameter value. The \itshape total \normalfont number of edges in the reduced mesh is described in the fourth column of  Table \ref{redmeshA}; the \itshape exact \normalfont subset with nontrivial $\Gamma(\boldsymbol\mu)$--intersection to be used as integration domain will have to be chosen online. For a more visual depiction of the same data, refer to Figure \ref{redmeshA2:1}. As can be deduced from Figure \ref{redmeshA2:2}, the majority  of the computational time for the online calculation of $\tilde{\theta}_{\mathbf{A}}(\boldsymbol\mu)\in\mathbb{R}^{  m_{\mathbf{A}}}$ (comparing the blue and magenta curves) is taken up by ghost penalty integrations on reduced mesh edges that are intersected by $\Gamma(\boldsymbol\mu)$. In fact, noting for higher DEIM dimensions in Figure \ref{A_reduced_mesh}  that the reduced mesh gradually fills the boundary variation zone $\left\{\Gamma(\boldsymbol\mu):\boldsymbol\mu \in \mathcal{P}\right\}$, few savings with respect to the full assemble of $\alpha(y_h(\boldsymbol\mu), \upsilon_h(\boldsymbol\mu))$ are expected for large $m_{\mathbf{A}}$. However, due to fewer diffusion--term and Nitsche--term integrations, the comparative advantage of the reduced assemble in still retained for larger $m_{\mathbf{A}}$, as can be seen in Figure \ref{redmeshA2:3}, albeit somewhat weakened. In Figure \ref{redmeshA2:3}, the speedup factor is simply the ratio of the time required to assemble the bilinear form on the full mesh to the reduced mesh assemble time; i.e., speedup indicates the number of  reduced mesh assembles that can be performed in the respective time for a \itshape single \normalfont full order assemble. Here, the full assemble time for the parameter value $\boldsymbol\mu_0=0.45$ has been computed $0.15684$ (as the arithmetic mean of 11 observed full assemble times) and compared to the respective reduced assemble times for the same parameter value $\boldsymbol\mu_0$ for the various DEIM dimensions. As can be verified in Figure \ref{redmeshA2:3}, speedup declines by almost 50\% when proceeding from $m_{\mathbf{A}}=2$ to $m_{\mathbf{A}}=90$ DEIM modes. However, the greatest part of this decrease is observed for the first 20 retained DEIM modes and is not so pronounced for reduced mesh expansion thereafter.   
    
\begin{table}[h!]
\begin{center}
\caption{Reduced mesh characteristics for the partial online assemble of $\tilde{\theta}_{\mathbf{A}}(\boldsymbol\mu)\in\mathbb{R}^{  m_{\mathbf{A}}}$ for progressively increasing DEIM dimension $m_{\mathbf{A}}$.}\label{redmeshA}
{\scriptsize
\begin{tabular}[scale=0.80]{cccc}
\hline
\# DEIM modes   & \# Elements  & \# Element calculations   & \# Edges  \\
  $m_{\mathbf{A}}$ &   &  (superfluous) &    \\
\hline
\hline
$1$ & $6 $ & $6$ & $12$  \\ 
$2$ & $13 $ & $13$ & $26$  \\ 
$5$ & $32$ & $32$   & $54$ \\ 
$10$ & $58$ & $62$ & $100$ \\ 
$15$ & $80$ & $92$  & $155$  \\ 
$20$ & $94$ & $118$   & $178$ \\ 
$25$ & $111$ & $147$  & $207$  \\ 
$30$ & $130$ & $178$  & $238$  \\ 
$35$ & $144$ & $203$  & $262$  \\ 
$40$ & $156$ & $229$  & $281$  \\ 
$60$ & $202$ & $323$  & $352$  \\ 
$90$ & $223$ & $432$  & $382$  \\  
\hline
\end{tabular}
}
\end{center}
\end{table}

%
%

Having considered the online assemble of $\tilde{\theta}_{\mathbf{A}}(\boldsymbol\mu)\in\mathbb{R}^{  m_{\mathbf{A}}}$, to complete our investigation of the computational performance of the DEIM procedure,  as described in Section \ref{DEIMsection}, we still have to take two more factors into account. Namely,  the time required to:
 \begin{itemize}
 \item perform the multiplication $\tilde{U}_{\mathbf{A}}\tilde{\theta}_{\mathbf{A}}(\boldsymbol\mu)\equiv \widehat{\text{vec}(\mathbf{A}_{\boldsymbol\mu})}$ in (\ref{DEIM}), where the factor  $\tilde{U}_{\mathbf{A}}=U_{\mathbf{A}}(P^T_{\mathbf{A}}\tilde{U}_{\mathbf{A}})^{-1}\in \mathbb{R}^{\mathcal{N}^n \times m_{\mathbf{A}}}$ is available from the offline stage.
 \item reshape the above approximation $\widehat{\text{vec}(\mathbf{A}_{\boldsymbol\mu})}$ to the vectorized form of the stiffness matrix $\mathbf{A}_{\boldsymbol\mu}$ in $\mathcal{N}\times \mathcal{N}$ matrix form. This operation is, of course, $m_{\mathbf{A}}$--independent.
 \end{itemize}
 
\begin{figure}
\centering
\subfloat[Online DEIM execution timing.]{%
\resizebox*{5.5cm}{!}{\includegraphics{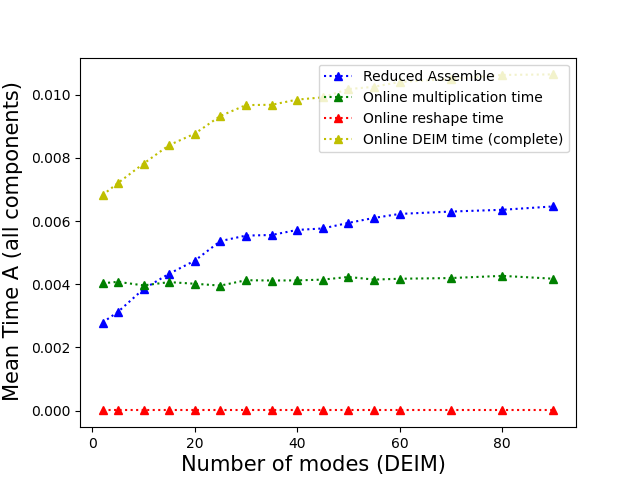}}}\hspace{5pt}
\subfloat[Reduced assemble speedup.]{%
\resizebox*{5.5cm}{!}{\includegraphics{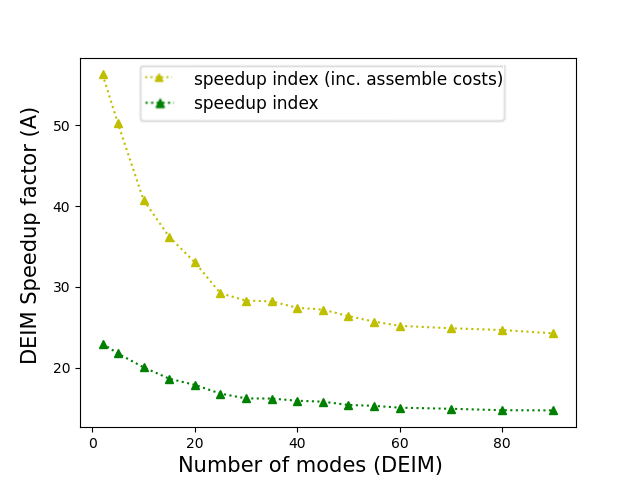}}}
\caption{ Efficiency of online DEIM procedure for the stiffness matrix $\mathbf{A}_{\boldsymbol\mu_0}$\normalfont. Combining all components: online times for $\tilde{\theta}_{\mathbf{A}}(\boldsymbol\mu) $ assemble, multiplication/reshape times and   their effect on the efficiency  of $\mathbf{A}_{\boldsymbol\mu_0}$--approximation for the parameter value $\boldsymbol\mu_0=0.45$ with respect to  increasing DEIM dimension $m_{\mathbf{A}}$.} 
  \label{redmeshA3:1}\label{redmeshA3:2} \label{redmeshA3}
\end{figure}

While the second of these contributions to the overall online computational burden  is small enough (mean reshape time has been computed  as $5.05149
e^{-06}$, indeed of different order when compared to reduced assemble times in Figure \ref{redmeshA2:2}), the multiplication operation remains a significant bottleneck which significantly compromises the efficiency of the online DEIM procedure. As revealed in Figure \ref{redmeshA3:1}, multiplication times (green curve) show a rather weak positive relationship with DEIM dimension and are comparable with online reduced assemble costs (blue curve, carried over from Figure \ref{redmeshA2:2}). Total online DEIM timings (yellow curve), computed as the sum of these three components, naturally remain competitive when compared with the full assemble procedure; recall that the mean full assemble time is $0.15684$. The effect of the additional computational costs, especially that of the multiplication operation on the overall efficiency of the online DEIM procedure is summarized in Figure \ref{redmeshA3:2}. The green curve in Figure \ref{redmeshA3:2} encodes the  speedup of the complete online DEIM procedure, which, although greatly decreased (compare with yellow curve, carried over from Figure \ref{redmeshA2:3}), nevertheless shows that the DEIM procedure remains almost 15 times faster than the full assemble procedure, even for as many DEIM modes as $m_{\mathbf{A}}=90$.
%

\begin{figure}[h!] 
  \centering
  \includegraphics[width=.4\linewidth]{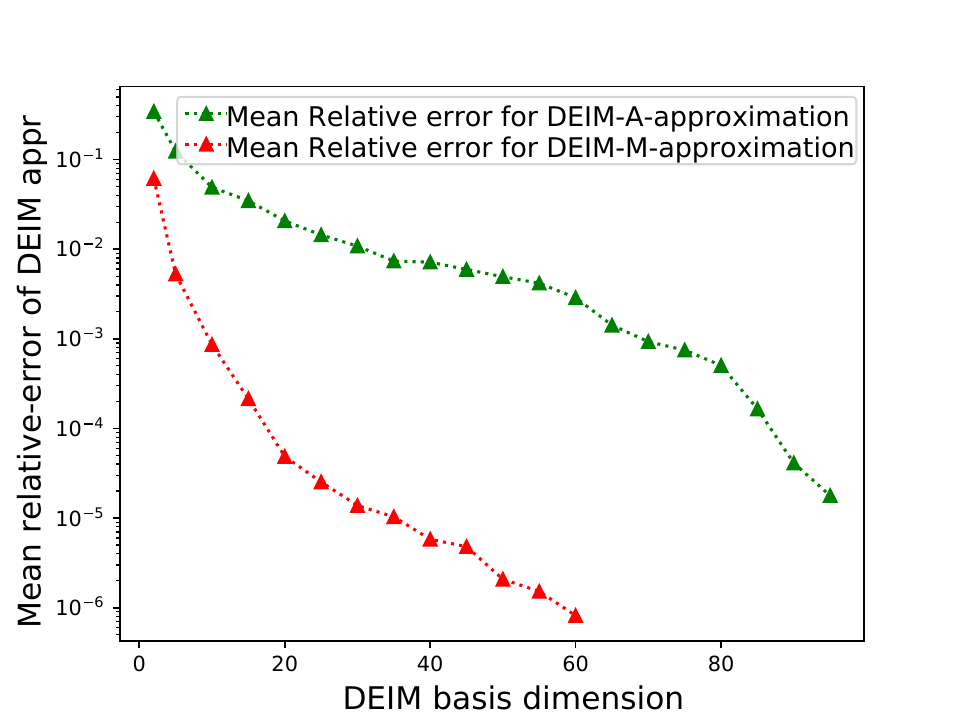}
  \caption{Reliability of   DEIM procedure for the nonaffine components of the system matrix $\mathcal{A}(\boldsymbol\mu)$\normalfont. Average relative $\left\|\cdot\right\|_2$--errors of the DEIM approximations of $\mathbf{A}_{\boldsymbol\mu}$ and $\mathbf{M}_{\boldsymbol\mu}$ for 30 randomly chosen parameter values $\boldsymbol\mu \in \mathcal{P}=[0.4,0.5]$  with respect to  increasing DEIM dimensions.}
  \label{redmeshA4} 
\end{figure}

Of course, the preceding efficiency analysis would have been in vain, were the DEIM procedure not capable to yield an accurate approximation to the stiffness matrix $\mathbf{A}_{\boldsymbol\mu}$. Its reliability is readily verified in Figure \ref{redmeshA4}, where the average relative $\left\|\cdot\right\|_2$--errors of the DEIM approximation of $\mathbf{A}_{\boldsymbol\mu}$ are shown to converge (green curve) to zero for increasing DEIM dimension $m_{\mathbf{A}}$. The average relative errors have been computed over a test sample of 30 randomly chosen parameters $\boldsymbol\mu \in \mathcal{P}=[0.4,0.5]$, selected during the online stage and different from the ones used offline to compute the stiffness matrix snapshots.  
{\green{
We recall that we start by the generation of a set of full order solutions of the parametrized problem under a parameter
values random choice. The final objective of reduced basis methods is to emulate any member of this
solution set with a low number of basis functions and this is based on a two-stage procedure,
the offline and the online stage. The reduced basis methods are predictive in the sense that the mean relative errors for a specific number
of random samples which are not the same as any of the samples used in the training stage,
allow these mean relative errors to remain the same after repeating the procedure and for
different samples. We clarify that we could use fewer snapshots, although, it is more favorable to employ more snapshots during the offline phase (which is computed only once) achieving in this way good accuracy for very few modes in the online phase rather than employing fewer snapshots which will cause less accuracy on the online stage.
}}

More precise details about the trade--off between accuracy and computational cost are summarized in Table \ref{relA}.

\begin{table}[h!]
\begin{center}
\caption{Trade--off between accuracy and computational cost of DEIM approximation to $\mathbf{A}_{\boldsymbol\mu}$ with respect to increasing DEIM dimension $m_{\mathbf{A}}$.}\label{relA}
{\scriptsize
\begin{tabular}[scale=0.80]{cc|cc|cc}
\hline
\# DEIM & $\left\|\cdot\right\|_2$-rel.  &  Reduced  & Reduced  &    DEIM time & DEIM  \\
modes $m_{\mathbf{A}}$  & error  & assemble time & assemble speedup& (online) & speedup \\
\hline
\hline
$1$ & $  6.11e-1 $ & - &-  &  - & - \\  
$2$ & $  3.43e-1 $ & $0.00279$& $56.29 $  &  $0.00682$ & $22.99 $  \\  
$5$ & $  1.24e-1$ &  $0.00312$ & $ 50.23$  & $0.00720$& $ 21.76$  \\  
$10$ & $  4.91e-2$ & $ 0.00385 $ & $ 40.75$  &  $0.00782 $& $20.07 $  \\
$15$ & $  4.91e-2$ &   $0.00384$& $ 36.21$  & $ 0.00840$& $ 18.67$  \\
$20$ & $  3.48e-2$ & $0.00385$ & $ 33.05$  &  $ 0.00877$& $ 17.89$  \\
$25$ & $  2.07e-2 $ &  $ 0.00547$  & $ 29.23$  &$0.00933 $& $ 16.81$  \\
$30$ & $  1.45e-2$ & $0.00554 $& $ 28.31$  &  $ 0.00967$ & $ 16.22$  \\
$35$ & $  1.09e-2$ & $ 0.00556$ & $ 28,22$  &  $ 0.00968$& $ 16.20$  \\
$40$ & $  7.26e-3 $ & $ 0.00572$& $ 27.43$  &  $ 0.00894$ & $ 15.93$  \\
$60$ & $  4.18e-3$ &  $ 0.00622$ & $25.19$ & $ 0.01040$  & $ 15.08$  \\
$90$ &   $1.72e-4$ 
 & $ 0.00646$ & $24.27 $  &  $ 0.01064$ &14.74 $ $  \\
\hline
\end{tabular}
}
\end{center}
\end{table}

\subsection{DEIM approximations of the mass matrix $\mathbf{M}_{\boldsymbol\mu}$ and right--hand side vectors  $\mathbf{b}_{\boldsymbol\mu}$ and $\mathbf{c}_{\boldsymbol\mu}$} \label{DEIMMnew}

In this subsection, we proceed to summarize the computational details from the application of the considerations in the previous section to the remaining non--affine components of $\mathcal{A}_N(\boldsymbol\mu)$ and $\beta_N(\boldsymbol\mu)$ for the system (\ref{ROM}). 

\begin{figure}
\centering
\subfloat[$m_{\mathbf{M}}=2$ modes.]{%
\resizebox*{4.2cm}{!}{\includegraphics{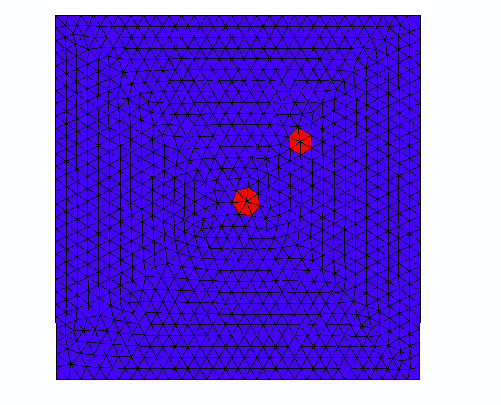}}}\hspace{5pt}
\subfloat[$m_{\mathbf{M}}=10$ modes.]{%
\resizebox*{4cm}{!}{\includegraphics{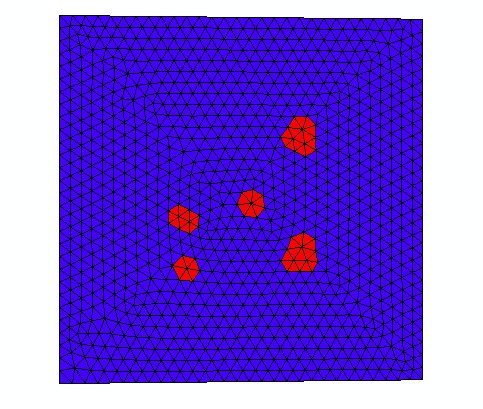}}}
\subfloat[$m_{\mathbf{M}}=35$ modes.]{%
\resizebox*{4cm}{!}{\includegraphics{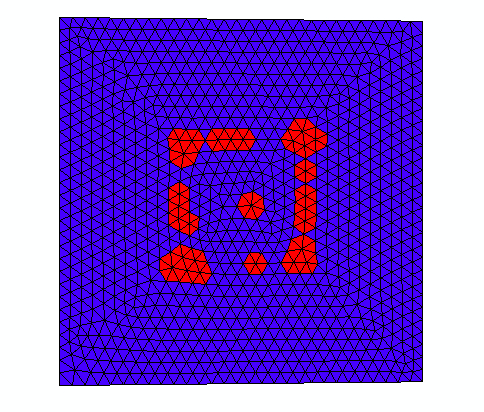}}}
\caption{Illustrations of reduced meshes for the DEIM approximation of the mass matrix $\mathbf{M}_{\boldsymbol\mu}$ with progressively increasing DEIM dimension $m_{\mathbf{M}}$ obtained though an offline stage with $M=370$ snapshots.}\label{M_reduced_mesh}
\end{figure}

A distinctive feature of the reduced mesh for the mass matrix $\mathbf{M}_{\boldsymbol\mu}$ is that the first DEIM mode to be selected in fact corresponds to a central degree of freedom  of $\Omega(\boldsymbol\mu)$, for all $\boldsymbol\mu \in [0.4,0.5]$. All remaining elements in the reduced mesh are still clustered around the band of boundary variation, as in the other cases for $\mathbf{A}_{\boldsymbol\mu}$,  $\mathbf{b}_{\boldsymbol\mu}$ and $\mathbf{c}_{\boldsymbol\mu}$; see Figure \ref{M_reduced_mesh}. This noteworthy remark conforms with the intuition that the single most important information for the mass matrix has to arise from a central point for all domain configurations; the remaining DEIM modes gradually sketch the domain boundary.

%

\begin{figure}
\centering
\subfloat[Online DEIM execution timing.]{%
\resizebox*{5cm}{!}{\includegraphics{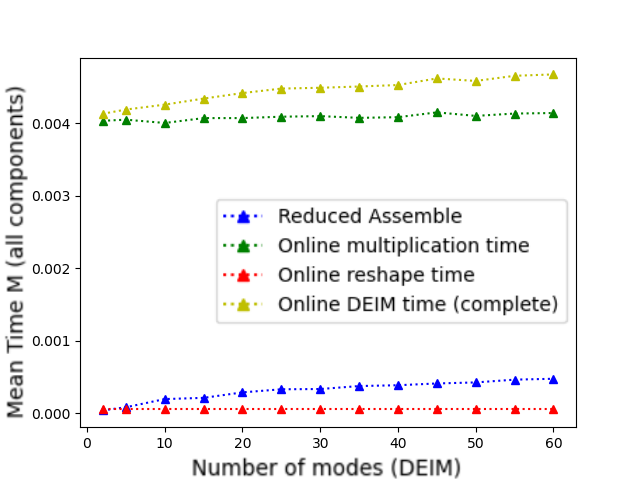}}}\hspace{5pt}
\subfloat[Reduced assemble speedup.]{%
\resizebox*{5cm}{!}{\includegraphics{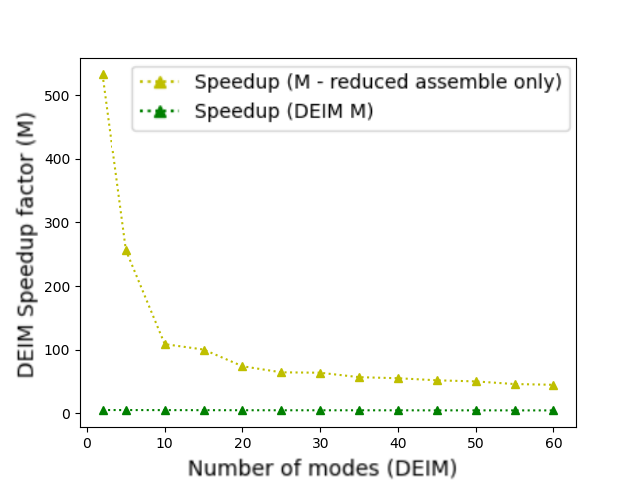}}}
\caption{Efficiency of online DEIM procedure for the mass matrix $\mathbf{M}_{\boldsymbol\mu_0}$\normalfont. Online times  and   their effect on the efficiency of $\mathbf{M}_{\boldsymbol\mu_0}$--approximation for the parameter value $\boldsymbol\mu_0=0.45$ with respect to  increasing DEIM dimension $m_{\mathbf{M}}$.} 
\label{redmeshM:1}\label{redmeshM:2}\label{redmeshM}
\end{figure}

Another interesting, yet rather obvious, feature of the mass matrix DEIM approximation is that, due to lack of ghost penalty integrations, the relative importance of the contribution of the operation $\tilde{U}_{\mathbf{M}}\tilde{\theta}_{\mathbf{M}}(\boldsymbol\mu)$ in overall execution times is much higher. Indeed, the  online DEIM--$\mathbf{M}_{\boldsymbol\mu}$ approximation times curve in Figure \ref{redmeshM:1} closely follows that of the multiplication operation, the reduced assemble and reshape times being much shorter and  insignificant. This is also reflected to the much compromised speedup index in Figure \ref{redmeshM:2}  (green curve), when the multiplication operation is taken into consideration.

%

More precise details about the trade--off between accuracy and computational cost for the DEIM approximation of the mass matrix $\mathbf{M}_{\boldsymbol\mu}$  are summarized in Table \ref{relM}. It is clear that, in view of lack of ghost penalty integration and the domination of multiplication operation cost, the online DEIM procedure for the approximation mass matrix is both \itshape more accurate \normalfont and \itshape less efficient \normalfont  than that for the CutFEM stiffness matrix; refer also to Figure \ref{redmeshA4}.

\begin{table}[h!]
\begin{center}
\caption{Trade--off between accuracy and computational cost of DEIM approximation to $\mathbf{M}_{\boldsymbol\mu}$ with respect to increasing DEIM dimension $m_{\mathbf{M}}$.}\label{relM}
{\scriptsize
\begin{tabular}[scale=0.80]{cc|cc|cc}
\hline
\# DEIM & $\left\|\cdot\right\|_2$-rel.  &  Reduced  & Reduced  &    DEIM time & DEIM  \\
modes $m_{\mathbf{M}}$  & error  & assemble time & assemble speedup& (online) & speedup \\
\hline
\hline
$1$ & $   3.18e-1 $ & - &-  &  - & - \\  
$2$  & $  6.16e-2$ &  $ 0.00004$ & $ 532.45$  & $ 0.00413$ & $5.17$  \\  
$5$  & $  5.28e-3 $ &  $ 0.00008$ & $ 256.31$  & $ 0.00419$ & $5.09$  \\  
$10$ & $ 8.61e-4  $ &  $ 0.00020$ & $ 108.69$  & $ 0.00425$ & $5.01$  \\
$15$ & $ 2.10e-4  $ &  $ 0.00021$ & $ 100.17$  & $ 0.00434$ & $4.92$  \\
$20$ & $ 4.88e-5  $ &  $ 0.00029$ & $ 073.93$  & $ 0.00441$ & $4.83$  \\
$25$ & $ 2.53e-5  $ &  $ 0.00033$ & $ 064.40$  & $ 0.00448$ & $4.77$  \\
$30$ & $ 1.38e-5   $ &  $ 0.00033$ & $ 063.76$  & $ 0.00449$ & $4.75$  \\
$35$ & $ 1.03e-5  $ &  $ 0.00038$ & $ 056.76$  & $ 0.00450$ & $4.73$  \\
$40$ & $ 5.81e-6   $ &  $ 0.00039$ & $ 055.12$  & $ 0.00452$ & $4.71$  \\
$50$ & $ 2.07e-6  $ &  $ 0.00043$ & $ 050.11$  & $ 0.00458$ & $4.65$  \\
$60$ & $ 8.17e-7 $ &  $ 0.00048$ & $ 044.79$  & $ 0.00467$ & $4.57$  \\
\hline
\end{tabular}
}
\end{center}
\end{table}

\begin{figure}
\centering
\subfloat[Online DEIM execution timing ($\mathbf{b}_{\boldsymbol\mu_0}$).]{%
\resizebox*{5.5cm}{!}{\includegraphics{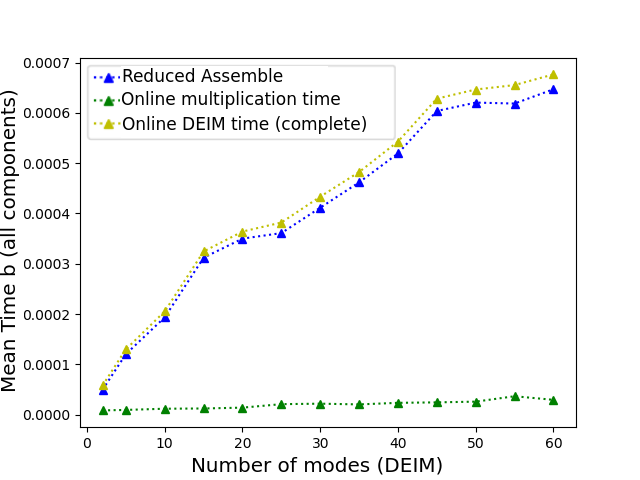}}}\hspace{5pt}
\subfloat[Reduced assemble speedup ($\mathbf{b}_{\boldsymbol\mu_0}$).]{%
\resizebox*{5.5cm}{!}{\includegraphics{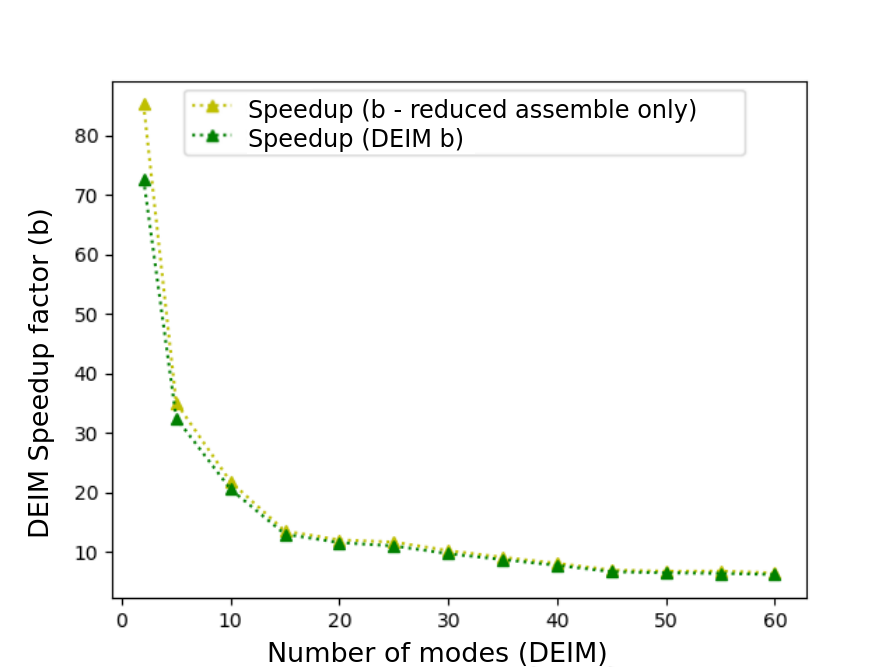}}}\hspace{5pt}
\subfloat[Online DEIM execution timing ($\mathbf{c}_{\boldsymbol\mu_0}$).]{%
\resizebox*{5.5cm}{!}{\includegraphics{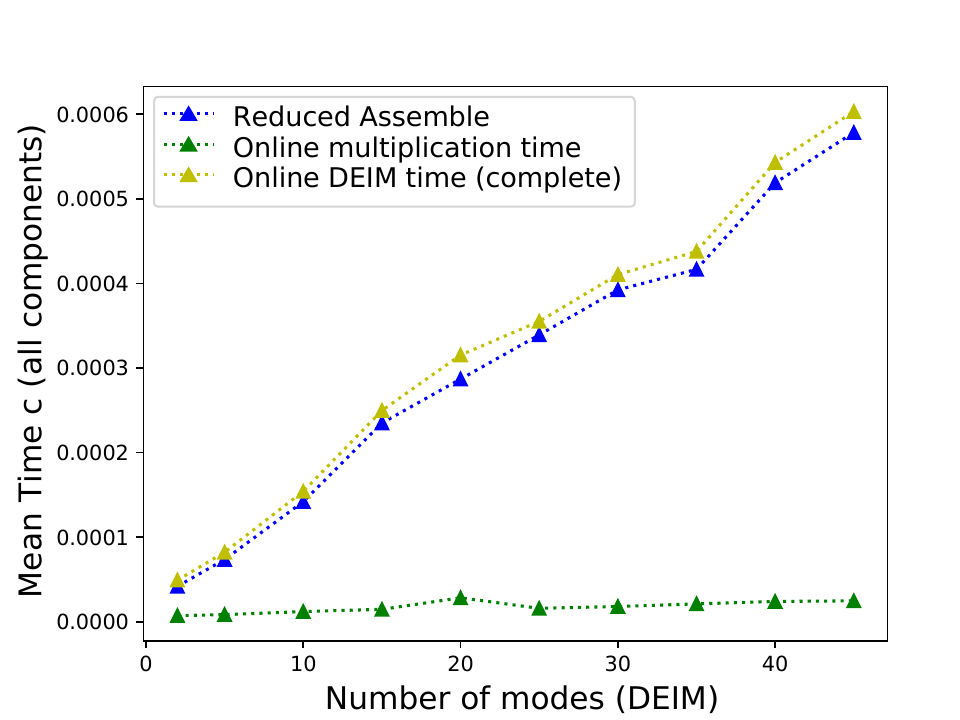}}}\hspace{5pt}
\subfloat[Reduced assemble speedup ($\mathbf{c}_{\boldsymbol\mu_0}$).]{%
\resizebox*{5.5cm}{!}{\includegraphics{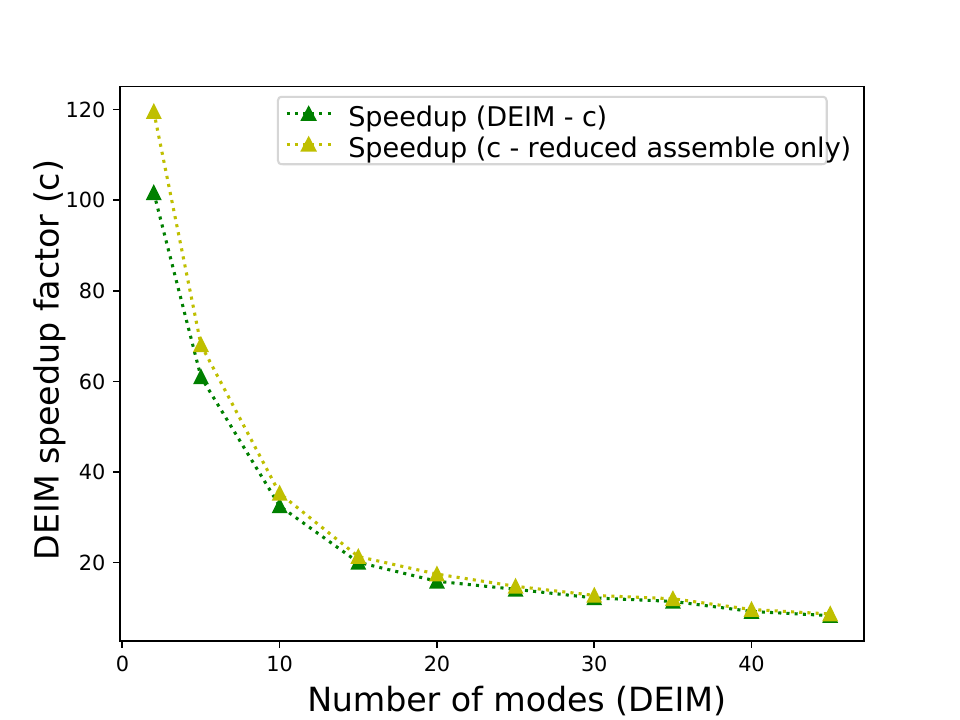}}}
\caption{Efficiency of online DEIM procedure  for the nonaffine components of the RHS vector $\beta(\boldsymbol\mu)$\normalfont. Online times  and   their effect on the efficiency of the DEIM approximations of $\mathbf{b}_{\boldsymbol\mu_0}$ and $\mathbf{c}_{\boldsymbol\mu}$ for the parameter value $\boldsymbol\mu_0=0.45$ with respect to  increasing DEIM dimensions $m_{\mathbf{b}}$ and $m_{\mathbf{c}}$, respectively.} 
\label{redmeshc:1}\label{redmeshc:2}\label{redmeshb:1}\label{redmeshb:2}\label{redmeshb}
\end{figure}

\begin{figure}[h!] 
  \centering
  \includegraphics[width=.4\linewidth]{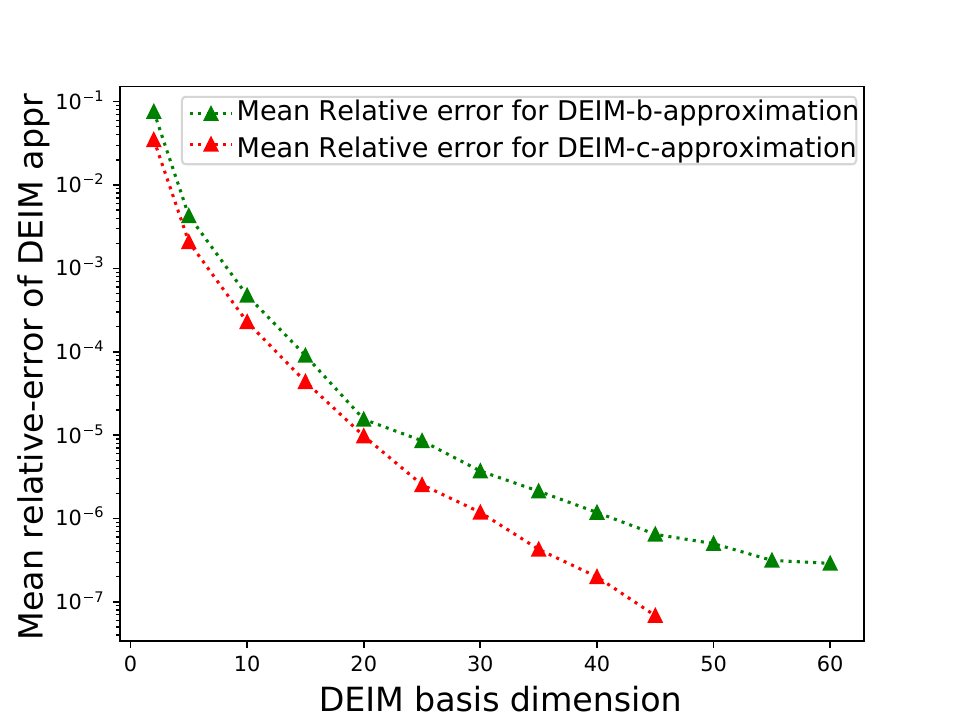}
  \caption{\itshape Reliability of   DEIM procedure for the nonaffine components of the RHS vector $\beta(\boldsymbol\mu)$\normalfont. Average relative $\left\|\cdot\right\|_2$--errors of the DEIM approximations of $\mathbf{b}_{\boldsymbol\mu}$ and $\mathbf{c}_{\boldsymbol\mu}$ for 30 randomly chosen parameter values $\boldsymbol\mu \in \mathcal{P}=[0.4,0.5]$  with respect to  increasing DEIM dimensions $m_{\mathbf{b}}$ and $m_{\mathbf{c}}$, respectively.}
  \label{redmeshb4} 
\end{figure}

Proceeding to the nonaffine components  of $\beta_N(\boldsymbol\mu)$, similar computational evidence is presented in Figure \ref{redmeshb}. As can be seen in Figures \ref{redmeshb:1} (a) and \ref{redmeshc:1} (c), in both cases of vectors $\mathbf{b}_{\boldsymbol\mu}$ and $\mathbf{c}_{\boldsymbol\mu}$, the online DEIM approximation times follow closely that of the reduced assemble. This indicates a noteworthy contrast in the application of DEIM procedure to the components of $\mathcal{A}_N(\boldsymbol\mu)$ and $\beta_N(\boldsymbol\mu)$; indeed, the multiplication operation for $\mathcal{N}$--dimensional vectors is apparently much cheaper than that for $\mathcal{N}^2$--dimensional vectors arising from matricial vectorization. Hence, the corresponding multiplication time (green curves) does not constitute a significant portion of the total DEIM execution time for either $\mathbf{b}_{\boldsymbol\mu}$ or $\mathbf{c}_{\boldsymbol\mu}$ and reduced assemble speedup factors remain relatively unaffected by its inclusion in the computations; see Figures \ref{redmeshb:2} (b) and \ref{redmeshc:2} (d). This is the reason we have opted to forfeit the inclusion of reduced assemble times and related speedup factors for $\mathbf{b}_{\boldsymbol\mu}$ and $\mathbf{c}_{\boldsymbol\mu}$   in Table \ref{relb}. The relevant DEIM reliability diagram can be found in Figure \ref{redmeshb4}.

%
%

%

\begin{table}[h!]
\begin{center}
\caption{Trade--off between accuracy and computational cost of DEIM approximations to $\mathbf{b}_{\boldsymbol\mu}$  and $\mathbf{c}_{\boldsymbol\mu}$ with respect to increasing DEIM dimension $m_{\mathbf{b}}$ and $m_{\mathbf{c}}$, respectively.}\label{relb}
{\scriptsize
\begin{tabular}[scale=0.80]{c|ccc|ccc}
\hline
\# DEIM & $\left\|\cdot\right\|_2$-rel.  &  DEIM time  & DEIM   &  $\left\|\cdot\right\|_2$-rel. &  DEIM time & DEIM  \\
modes   & error  & (online) &   speedup& error&  (online) & speedup \\ 
$m_{\mathbf{b}}$/$m_{\mathbf{c}}$ & ($\mathbf{b}_{\boldsymbol\mu}$) & ($\mathbf{b}_{\boldsymbol\mu}$) & ($\mathbf{b}_{\boldsymbol\mu}$)  & ($\mathbf{c}_{\boldsymbol\mu}$)  & ($\mathbf{c}_{\boldsymbol\mu}$)  &($\mathbf{c}_{\boldsymbol\mu}$)  \\
\hline
\hline
$1$ & $  3.18e-1  $  & - & -  & $  2.08e-1  $   & - & -  \\  
$2$ &  $ 7.60e-2  $ & $0.00006$ & $ 72.65$  & $ 3.49e-2  $  & $0.00005$ & $101.65$  \\  
$5$ &  $ 4.29e-3  $ & $0.00013$ & $ 32.53$  & $ 2.08e-3  $  & $0.00008$ & $060.99$   \\  
$10$ & $ 4.76e-4  $ & $0.00021$ & $ 20.57$  & $ 2.30e-4  $  & $0.00015$ & $032.52 $   \\
$15$ & $ 9.05e-5  $ & $0.00032$ & $ 13.03$  & $ 4.38e-5  $  & $0.00025$ & $020.08$   \\
$20$ & $ 1.55e-5  $ & $0.00036$ & $ 11.61$  & $ 9.79e-6  $  & $0.00032$ & $015.90$   \\
$25$ & $ 8.57e-6  $ & $0.00038$ & $ 11.07$  & $ 2.54e-6  $  & $0.00036$ & $014.12$    \\
$30$ & $ 3.72e-6  $ & $0.00043$ & $ 09.75$  & $ 1.19e-6  $  & $0.00041$ & $012.22$  \\
$35$ & $ 2.13e-6  $ & $0.00048$ & $ 08.76$  & $ 4.28e-7  $  & $0.00043$ & $011.45$  \\
$40$ & $ 1.18e-6  $ & $0.00054$ & $ 07.79$  & $ 2.00e-7  $  & $0.00054$ & $009.24$  \\
$45$ & $ 6.48e-7 $ & $0.00063$ & $ 06.73$  &  $ 6.86e-8 $  & $0.00060$ & $008.32$  \\
$60$ & $ 2.90e-7  $ & $0.00066$ & $ 06.24$  &  $ - $ & $ -  $ &$ - $   \\
\hline
\end{tabular}
}
\end{center}
\end{table} 

\subsection{POD/DEIM approximation: The optimal control problem} 

We are now ready to combine the previous components and proceed to full-scale  optimal control problem (\ref{ROM}) with force $f(x,y)=xy$, homogeneous Dirichlet conditions, desired state $y_{d}=\frac{1}{2\pi}\sin(\pi x)\cos(\pi x)$ and control distributed throughout   the square domain  $\Omega(\boldsymbol\mu)$ parametrized by the {\magenta{level set}} function (\ref{level1}). Recall that the common background domain is taken $\mathcal{B}=[-0.3,2.3]\times [-0.3,2.3]$ with mesh parameter $h=0.09$, resulting in a corresponding mesh $\mathcal{B}_h$ with 1944 elements and a high--fidelity space $V_h$ of dimension $\mathcal{N}=1031$. All calculations have been realized  in a python3 framework, enriched by the ngsolve/ngsxfem software packages, \cite{Scho14,LeHePreWa21}. 
 
\begin{figure}
\centering
\subfloat[Normalized eigenvalues of the correlation matrices related to $(y,u,p)$.]{%
\resizebox*{5.5cm}{!}{\includegraphics{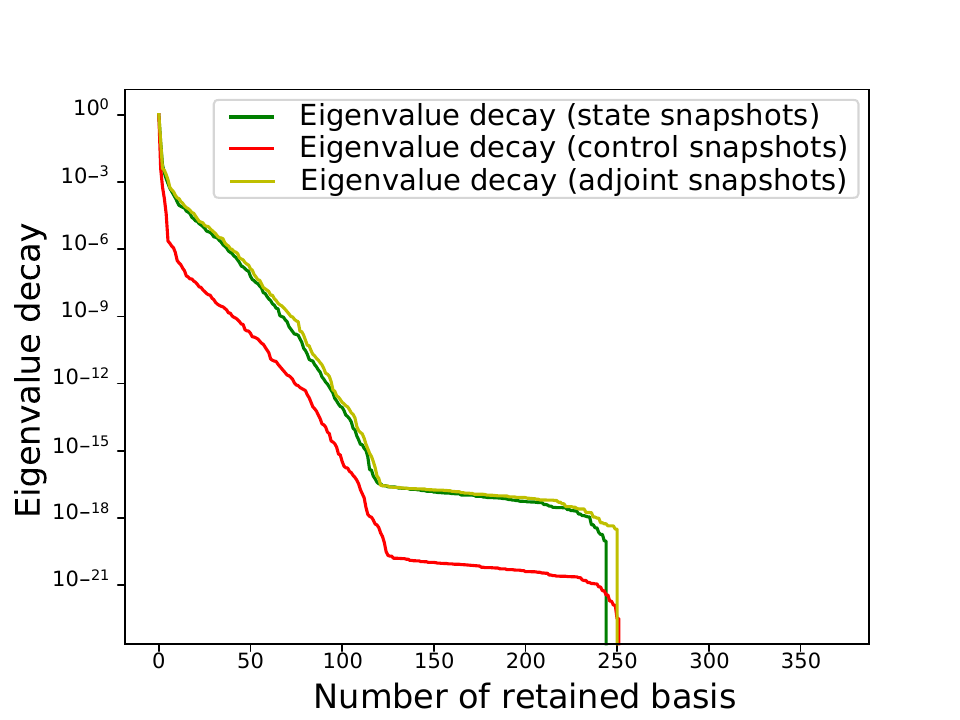}}}\hspace{5pt}
\subfloat[Normalized eigenvalues of the correlation matrices related to the nonaffine components of the system matrix $\mathcal{A}(\mu)$ and the RHS--vector $\beta(\mu)$.]{%
\resizebox*{5.5cm}{!}{\includegraphics{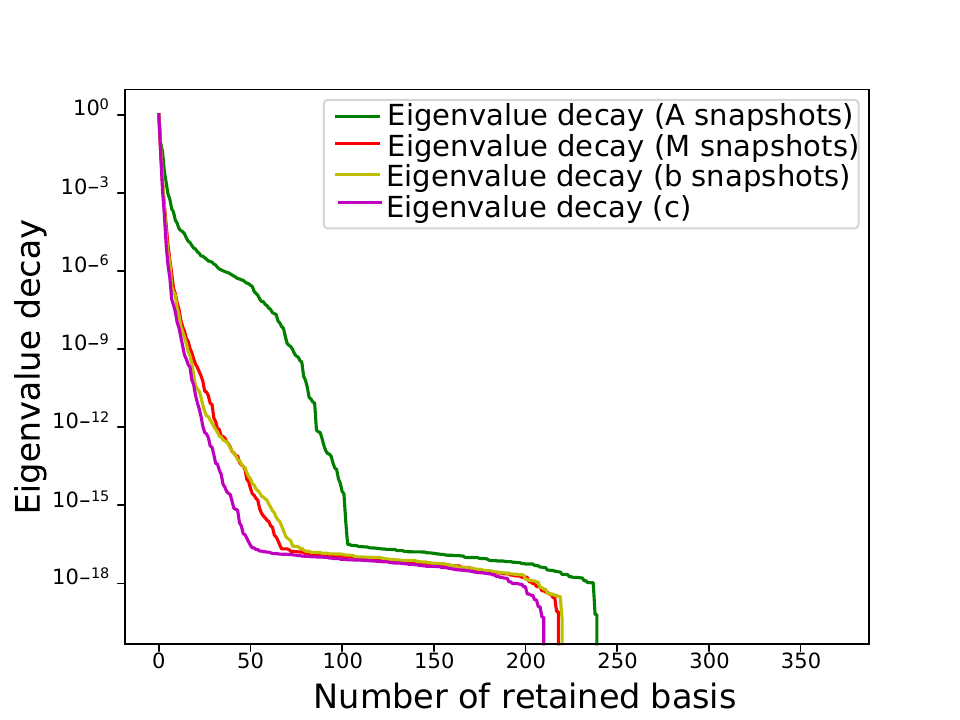}}}
\caption{Eigenvalue decay for the POD/DEIM reduction  of the OCP with non--homogeneous Dirichlet conditions on the boundary defined by the {\magenta{level set}} curve with $\phi(x,y;\boldsymbol\mu)=\left|x-1\right|+\left|y-1\right|+\left|\left|x-1\right|-\left|y-1\right|\right|-2\boldsymbol\mu$. 370 snapshots have been used for different parameter values $\boldsymbol\mu \in   [0.4, 0.5]$. }
  \label{fig:eigs_POD1a}\label{fig:eigs_POD2a}\label{homo_errorsa} 
\end{figure}

To {\magenta{setup the POD and the DEIM model in the offline phase and to}} explore the solution manifold of the {\magenta{optimal control problem (OCP)}}, a training set of $M=370$ solution snapshots was {\magenta{based on $370$ parameters}} generated through the full--order $3093 \times 3093$ ($=3\mathcal{N} \times 3 \mathcal{N}$)--dimensional linear system (\ref{system}) for random values of  $\boldsymbol\mu \in  \mathcal{P}= [0.4, 0.5]$. 
{\magenta{The offline phase required time $508.25$ seconds.}}
The decay of the eigenvalues $\left\{\lambda_{ji}\right\}_{i=1}^{370}$ of the snapshot correlation matrices $\mathbf{C}_j$ ($j=y,u,p$) related to the state, control, and adjoint state respectively is depicted in Figure  \ref{fig:eigs_POD1a}, plotting in each case the points $\left(i, \frac{\lambda_{ji}}{\lambda_{j1}}\right)$.  Here, the eigenvalues $\left\{\lambda_{ji}\right\}_{i=1}^{370}$ are indexed in non--increasing order.  The decay is rather steep and shows that a lot of the energy of the full--order model may be retained, using only a few basis functions.  

\begin{figure}
\centering
\subfloat[1st POD basis (state).]{%
\resizebox*{4.5cm}{!}{\includegraphics{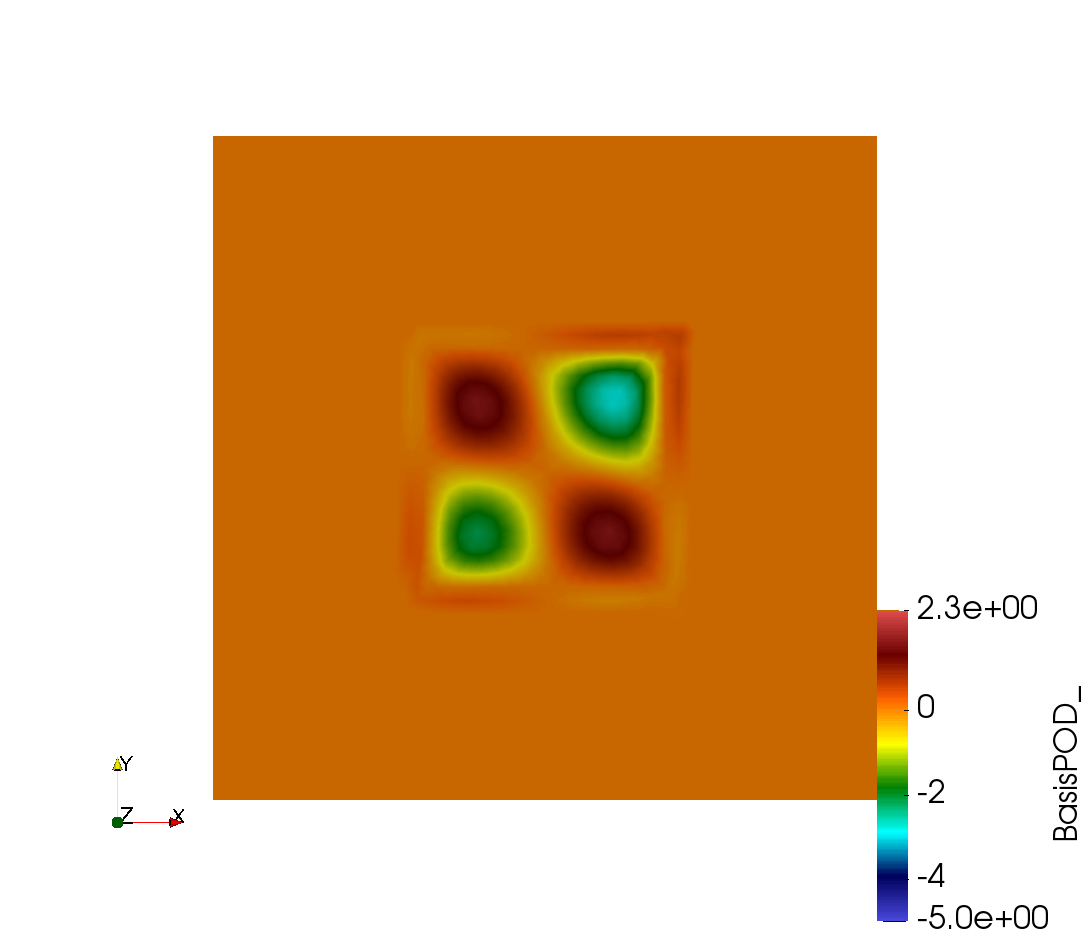}}}\hspace{5pt}
\subfloat[2nd POD basis (state).]{%
\resizebox*{4.5cm}{!}{\includegraphics{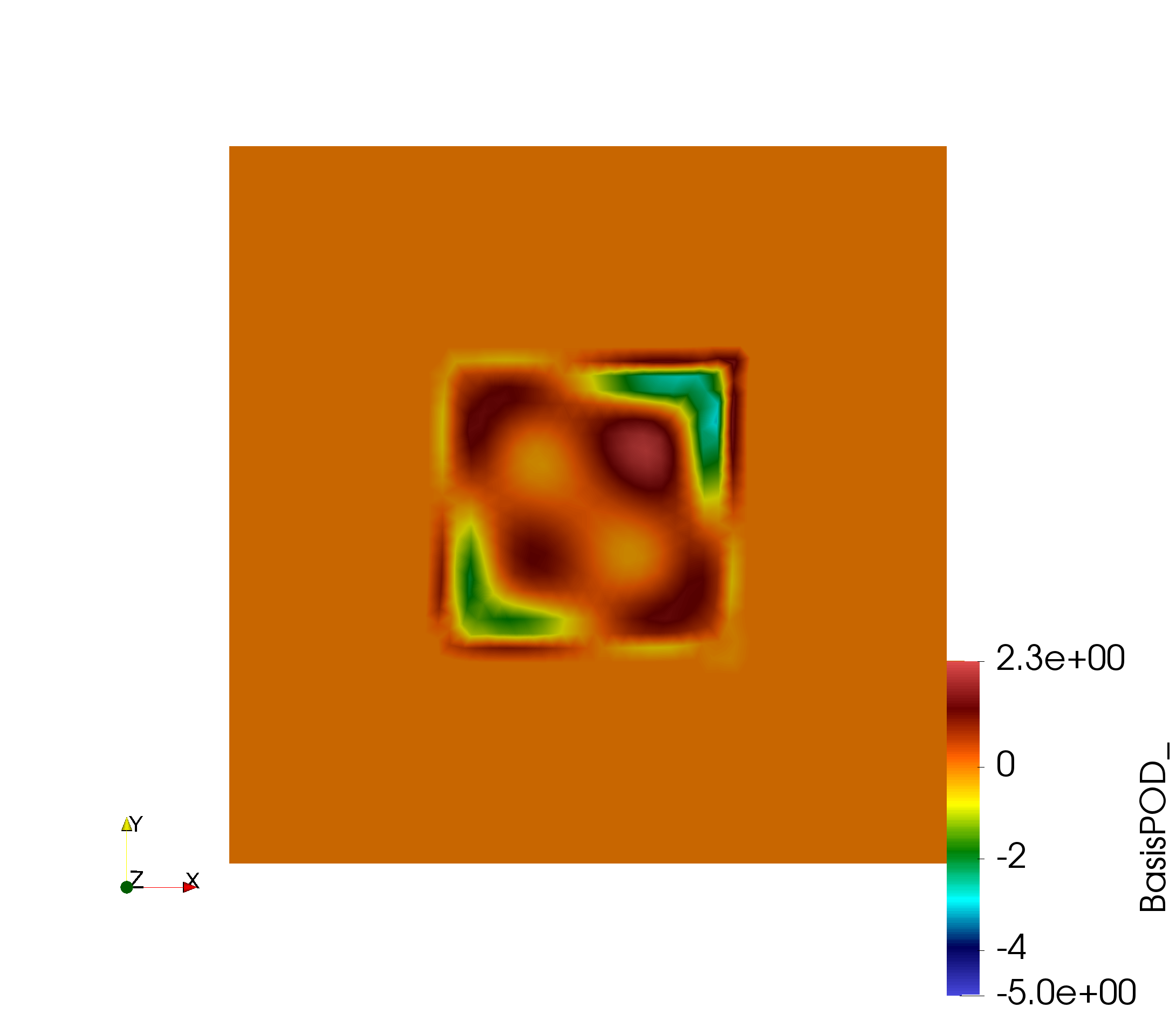}}}\hspace{5pt}
\subfloat[3rd POD basis (state).]{%
\resizebox*{4.5cm}{!}{\includegraphics{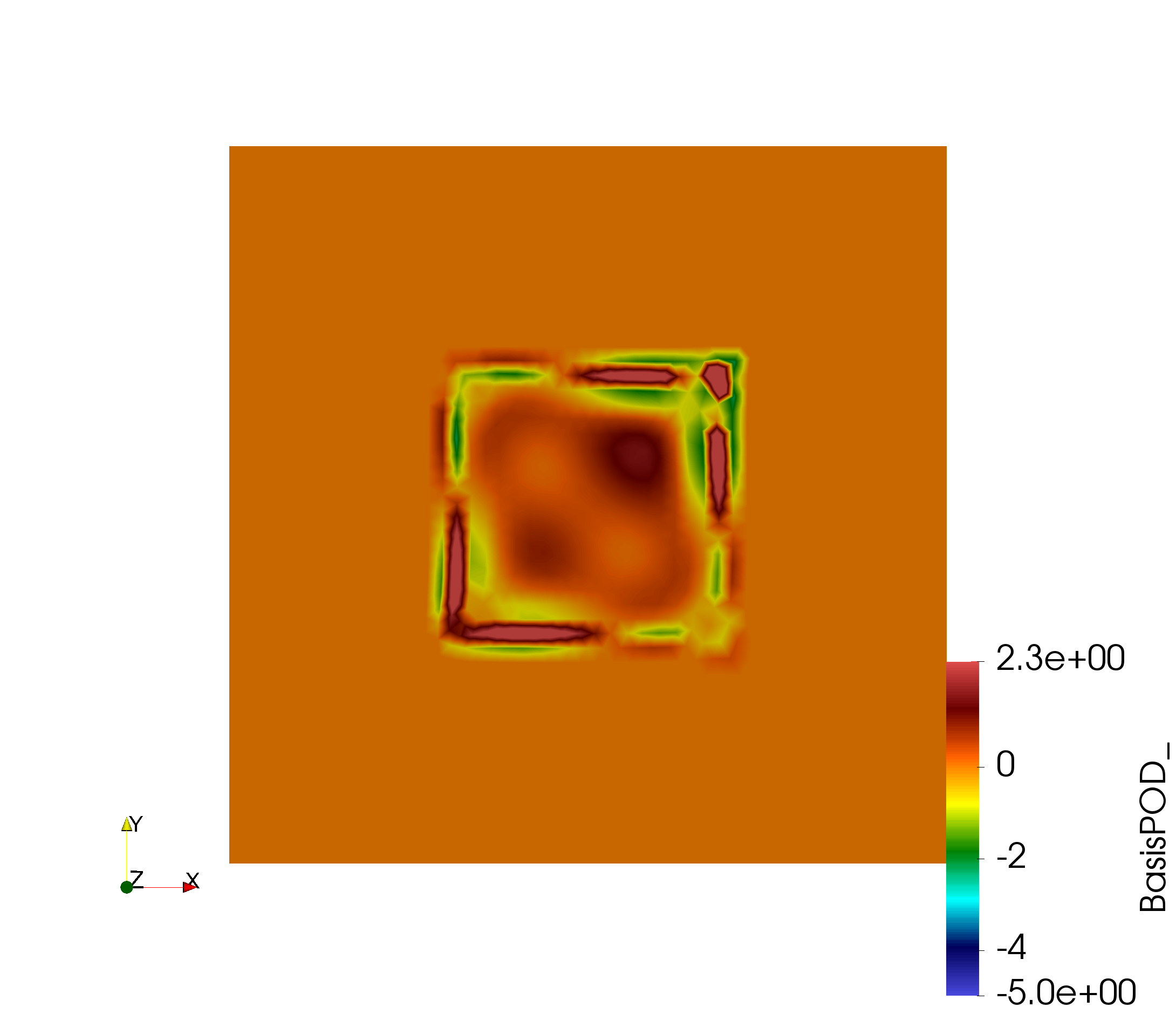}}}
\\
\subfloat[1st POD basis (control).]{%
\resizebox*{4.5cm}{!}{\includegraphics{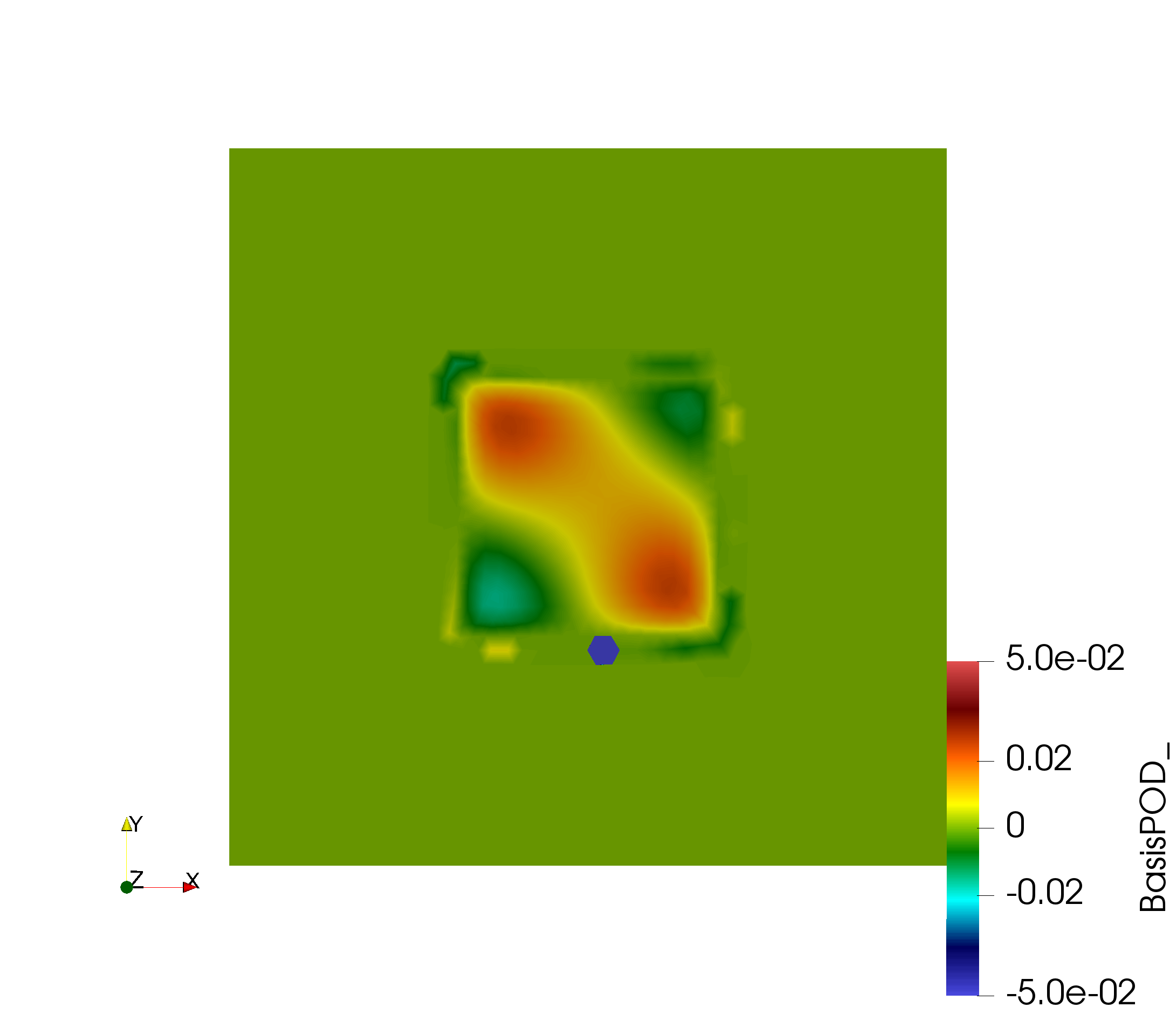}}}\hspace{5pt}
\subfloat[2nd POD basis (control).]{%
\resizebox*{4.5cm}{!}{\includegraphics{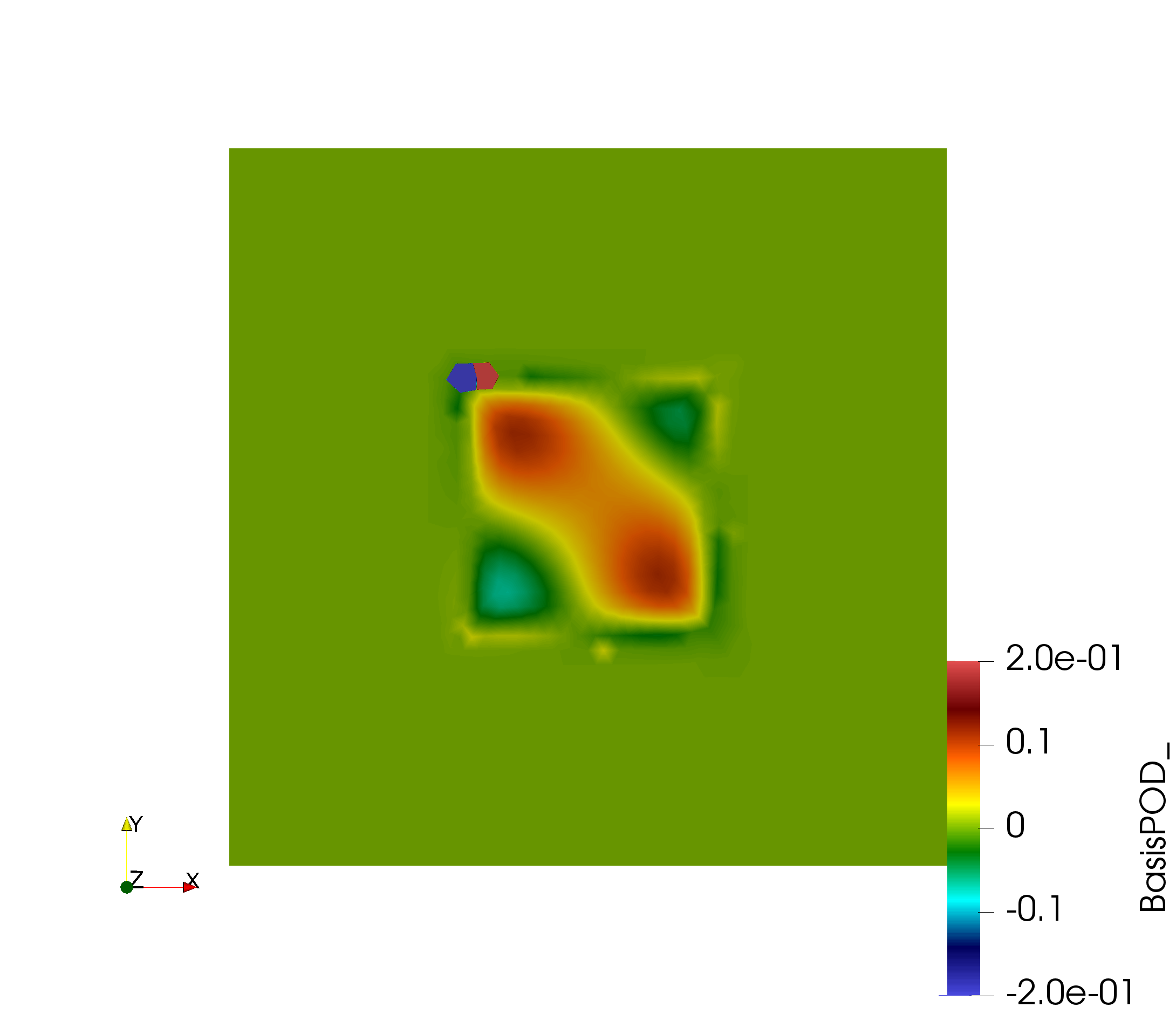}}}\hspace{5pt}
\subfloat[3rd POD basis (control).]{%
\resizebox*{4.5cm}{!}{\includegraphics{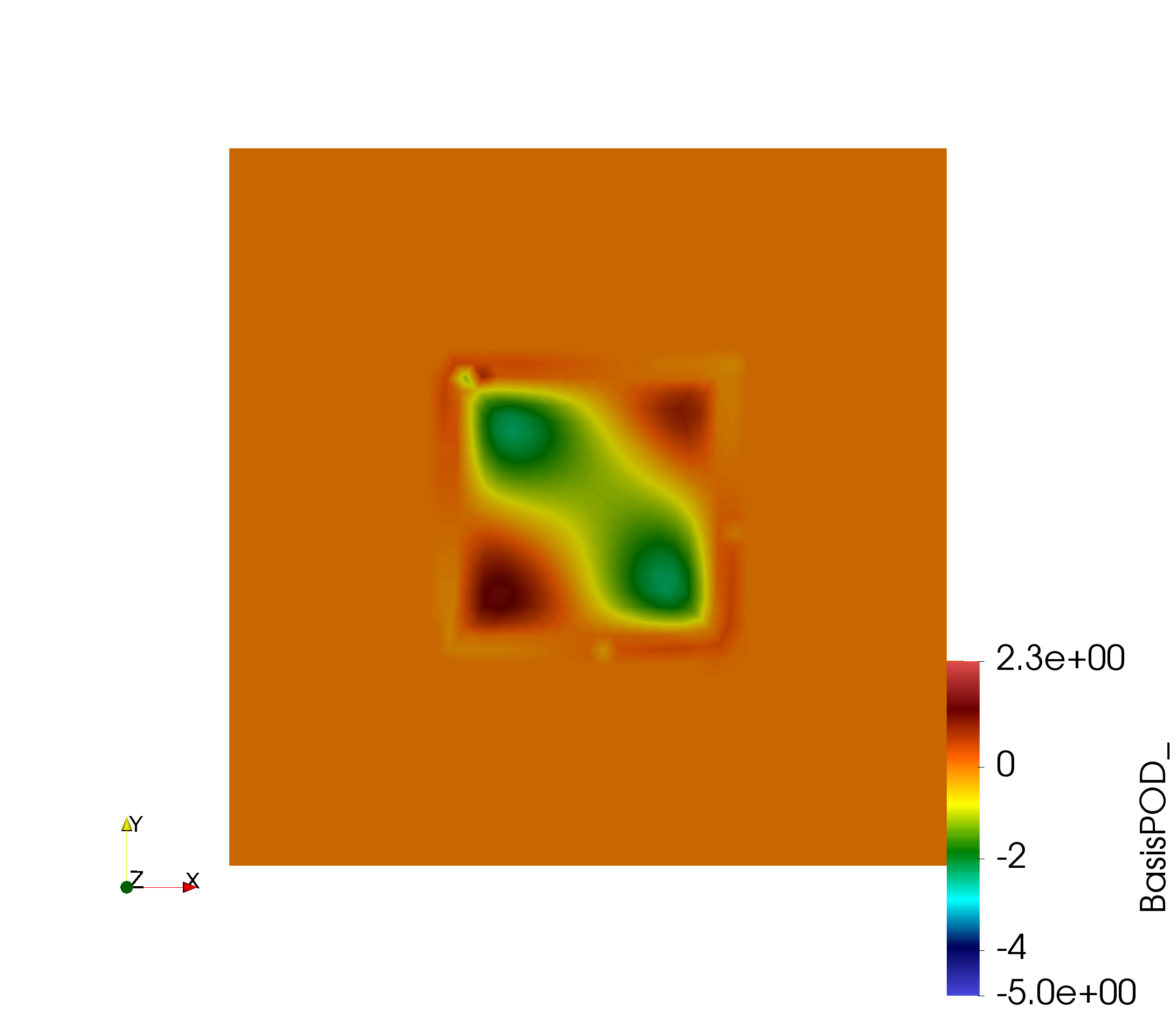}}}
\\
\subfloat[1st POD basis (adjoint).]{%
\resizebox*{4.5cm}{!}{\includegraphics{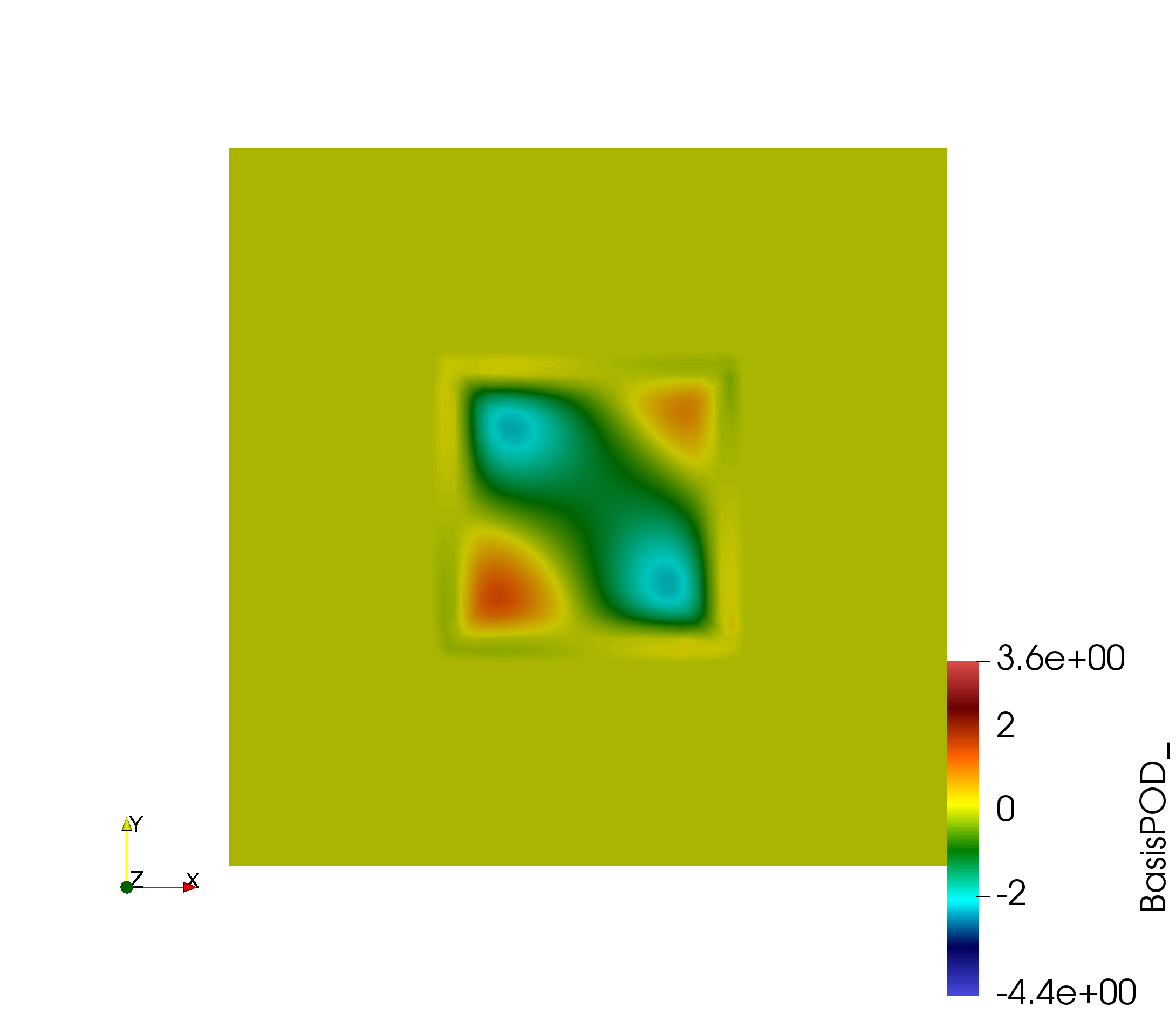}}}\hspace{5pt}
\subfloat[2nd POD basis (adjoint).]{%
\resizebox*{4.5cm}{!}{\includegraphics{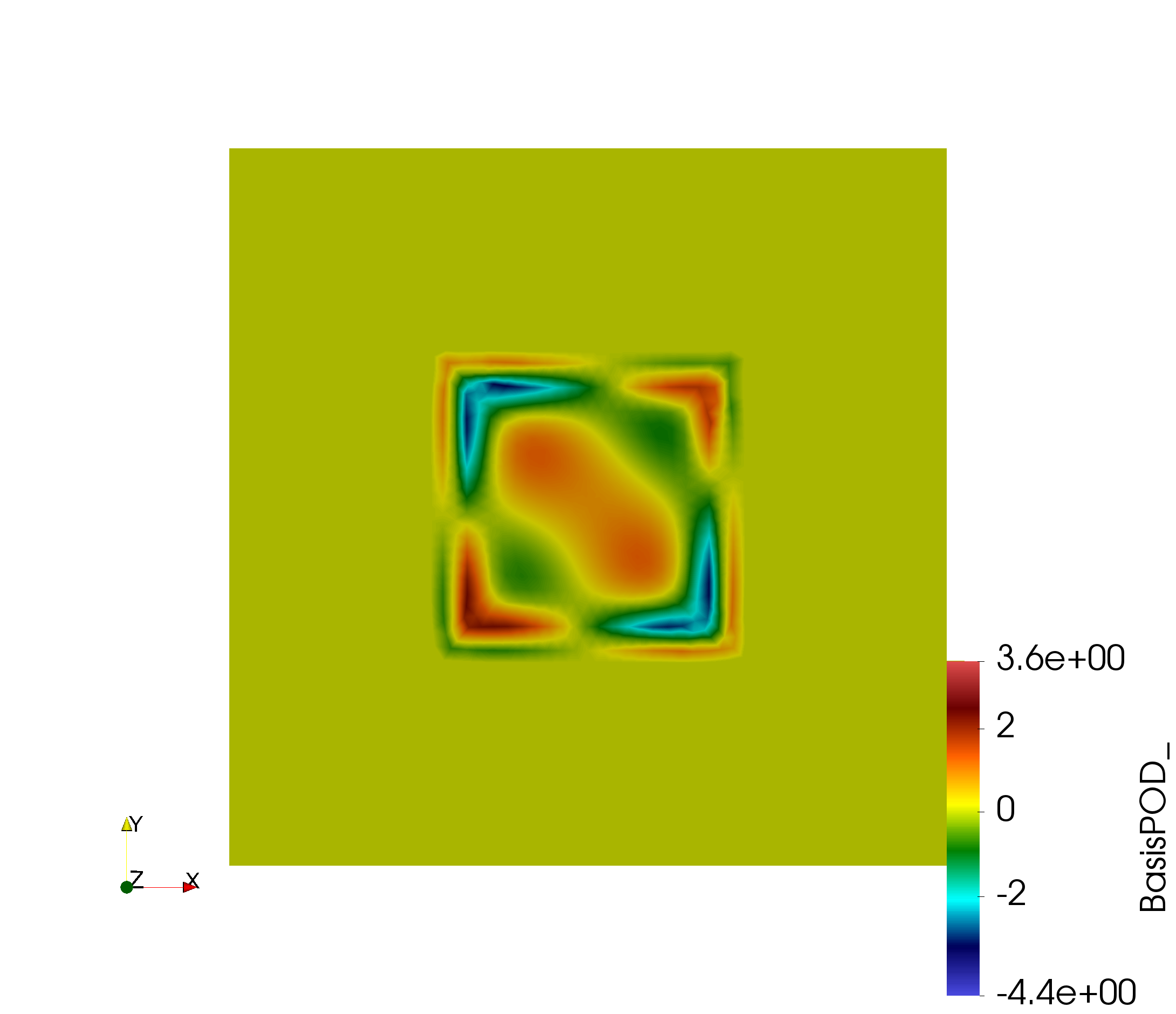}}}\hspace{5pt}
\subfloat[3rd POD basis (adjoint).]{%
\resizebox*{4.5cm}{!}{\includegraphics{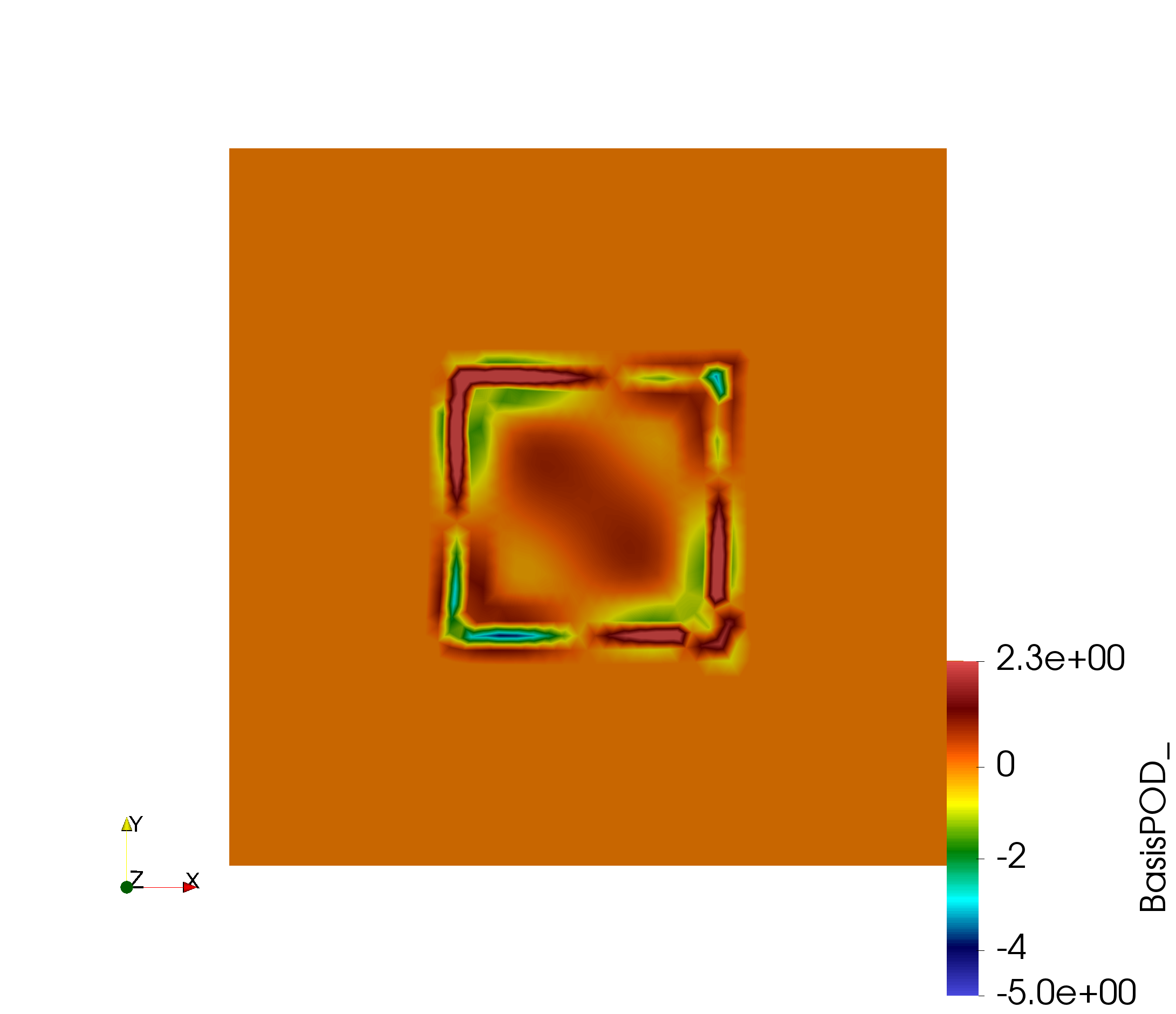}}}
    \caption{Leading POD basis functions for state/control/adjoint state reduction for the optimal control problem based on $M=370$ snapshots.}\label{fig:base1}\label{fig:base2}\label{fig:base2a}
    \label{fig:base3}\label{fig:base4}\label{fig:base5}\label{fig:base6}\label{fig:base6a}\label{PODbasis}
\end{figure}

To perform the Galerkin--POD reduction, we consider a tolerance $\epsilon_{POD}=10^{-6}$ and choose in each case ($j=y, u , p$) the basis dimension $N_j$ as the minimum integer such that
\begin{equation*}
\frac{\sum_{i=1}^{N_j}\lambda_{ji}}{\sum_{i=1}^{M}\lambda_{ji}}\geq 1-\epsilon_{POD}.
\end{equation*}
According to this criterion with tolerance $\epsilon_{POD}=10^{-5}$, the numbers of retained POD--modes for state, control and adjoint were determined $\left(N_y, N_u, N_p\right)=(31,9,31)$. Thus, a ROM of order  $2N_y+N_u+2N_p=133$ is obtained in  (\ref{ROM}), as opposed to the full 3093(=$3\mathcal{N}$)--dimensional system. The  leading POD basis functions are illustrated in Figure \ref{PODbasis}.

%

To ensure its efficient resolution for each new parameter value $\boldsymbol\mu \in [0.4,0.5]$ through rapid online assembly of the corresponding system matrix $\mathcal{A}_N(\boldsymbol\mu)$ and right--hand size $\beta_N(\boldsymbol\mu)$, DEIM approximations were performed on all their nonaffine blocks $\mathbf{A}_{\boldsymbol\mu}$,  $\mathbf{M}_{\boldsymbol\mu}$,  $\mathbf{b}_{\boldsymbol\mu}$, and  $\mathbf{c}_{\boldsymbol\mu}$. As explained in Sections \ref{DEIMAnew}--\ref{DEIMMnew}, these components are approximated in lower--dimensional subspaces of respective dimensions 
$$(m_{\mathbf{A}},m_{\mathbf{M}},m_{\mathbf{b}},m_{\mathbf{c}})=(83, 25, 21, 19)$$ determined offline by POD on corresponding sets of $M=370$ snapshots with tolerance $\epsilon_{DEIM}=10^{-10}$. The decay of the eigenvalues for the respective correlation matrices is depicted in Figure \ref{fig:eigs_POD2a}. Figures   \ref{redmeshA4} and   \ref{redmeshb4}  indicate the average $\left\|\cdot\right\|_2$ errors of the DEIM approximations for the blocks of the system matrix and its RHS--vector, respectively. The error has been computed over a test sample of 30 parameters, chosen during the online stage, different from the ones used to compute the snapshots. 
The approximation error is most pronounced for the stiffness matrix $\mathbf{A}_{\boldsymbol\mu}$, while the mass matrix $\mathbf{M}_{\boldsymbol\mu}$ and the vectors $\mathbf{b}_{\boldsymbol\mu}$, $\mathbf{c}_{\boldsymbol\mu}$ are markedly easier to approximate. Nevertheless, Figures   \ref{redmeshA4} and   \ref{redmeshb4}   highlight the fact  that the proposed procedure is able to {\magenta{accurately}} capture the profile of the system matrix and its RHS--vector, relying on a rather restricted number of basis functions. 

This procedure results in affine expansions for $\mathcal{A}_N(\boldsymbol\mu)$ and $\beta_N(\boldsymbol\mu)$ with $Q_{\mathcal{A}}=m_{\mathbf{A}}+m_{\mathbf{M}}=108$
 and $Q_{\beta}=m_{\mathbf{b}}+m_{\mathbf{c}}=40$ terms respectively, as in Section \ref{online}. Hence, for each online computation, we need only interpolate the system matrix $\mathcal{A}_N(\boldsymbol\mu)$ and right--hand size $\beta_N(\boldsymbol\mu)$ on $Q_{\mathcal{A}}=108$  and $Q_{\beta}=40$  selected components by the DEIM algorithm, instead of assembling the full 3093(=$3\mathcal{N}$)--dimensional system.
 
\begin{figure}
\centering
\subfloat[High--fidelity state.]{%
\resizebox*{4.5cm}{!}{\includegraphics{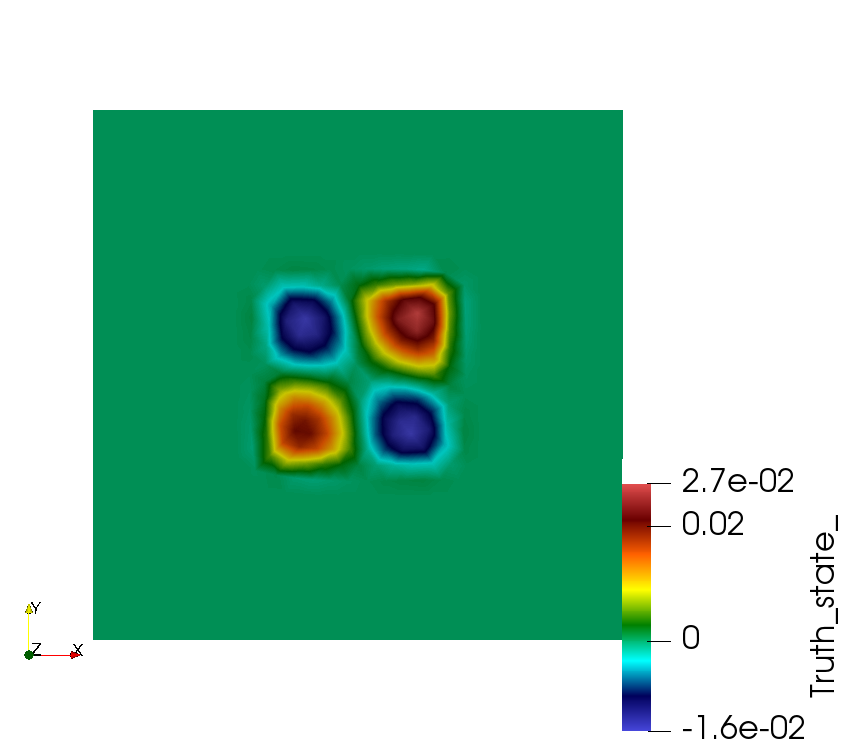}}}\hspace{5pt}
\subfloat[Reduced order state.]{%
\resizebox*{4.5cm}{!}{\includegraphics{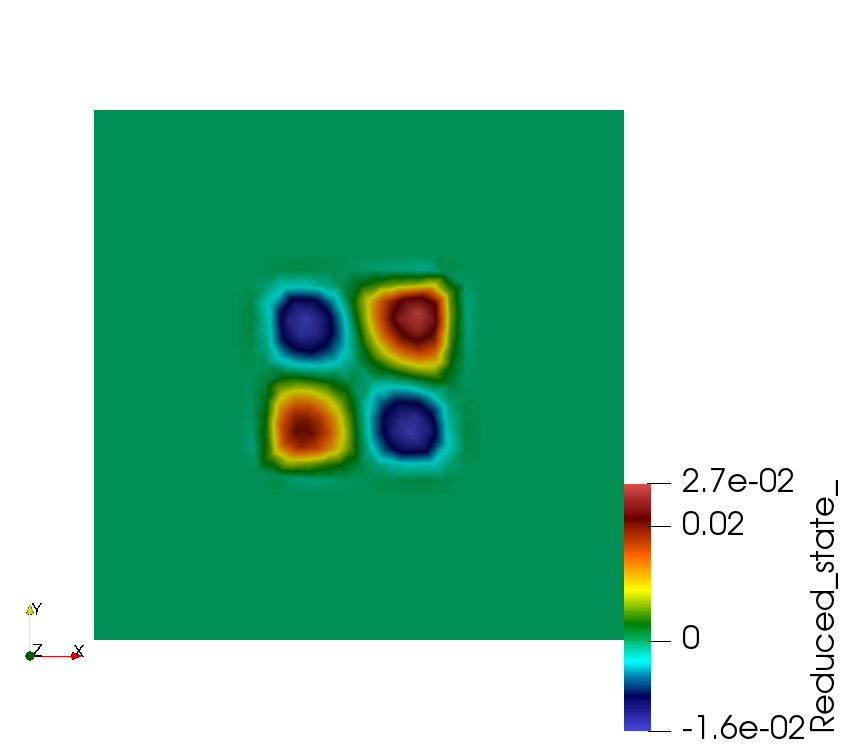}}}\hspace{5pt}
\subfloat[Error (state).]{%
\resizebox*{4.5cm}{!}{\includegraphics{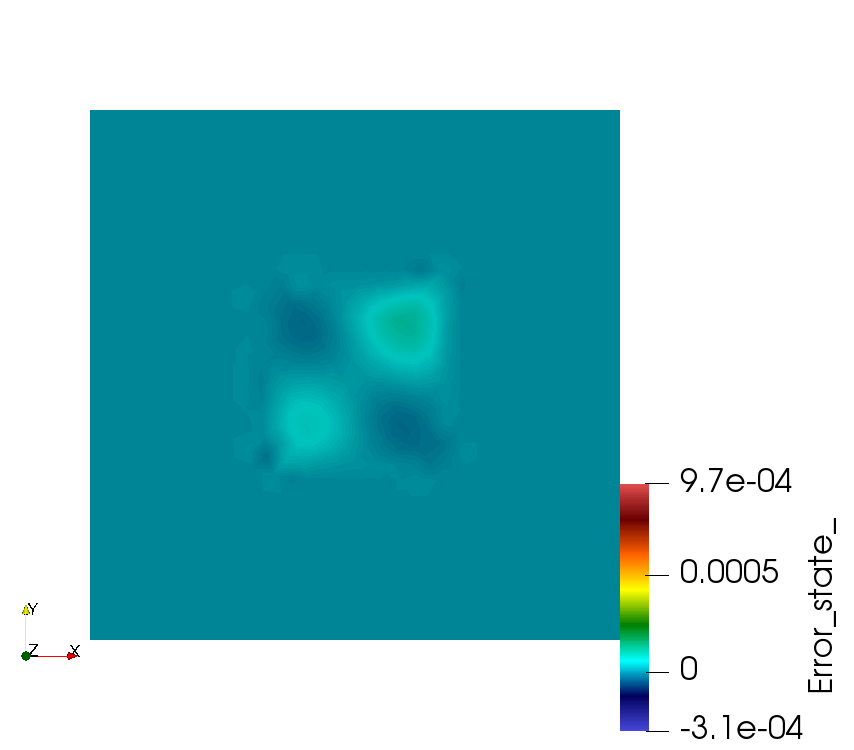}}}
\\
\subfloat[High--fidelity control.]{%
\resizebox*{4.5cm}{!}{\includegraphics{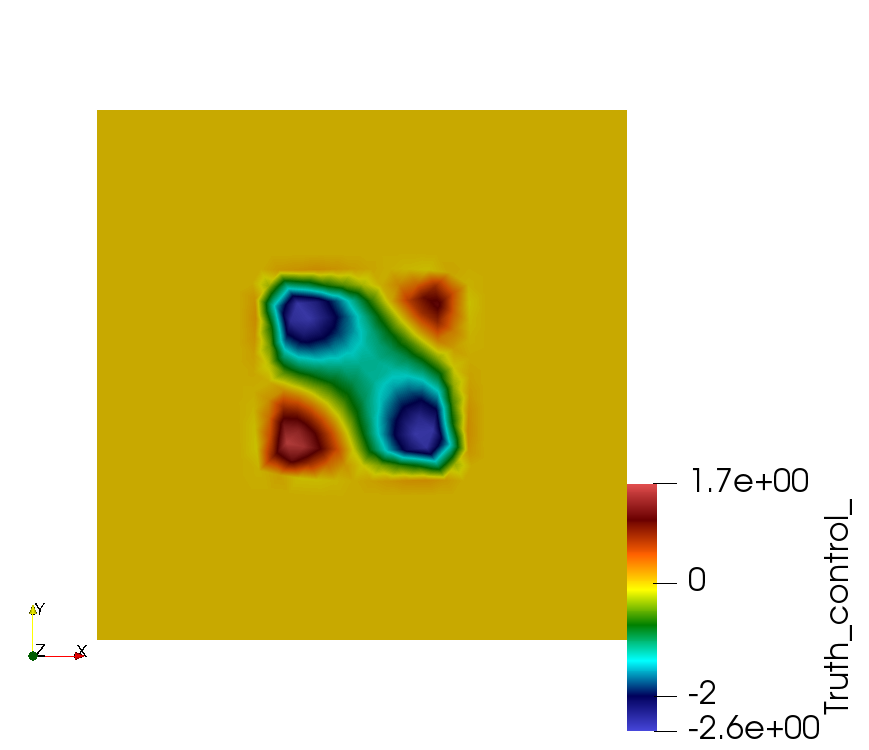}}}\hspace{5pt}
\subfloat[Reduced order control.]{%
\resizebox*{4.5cm}{!}{\includegraphics{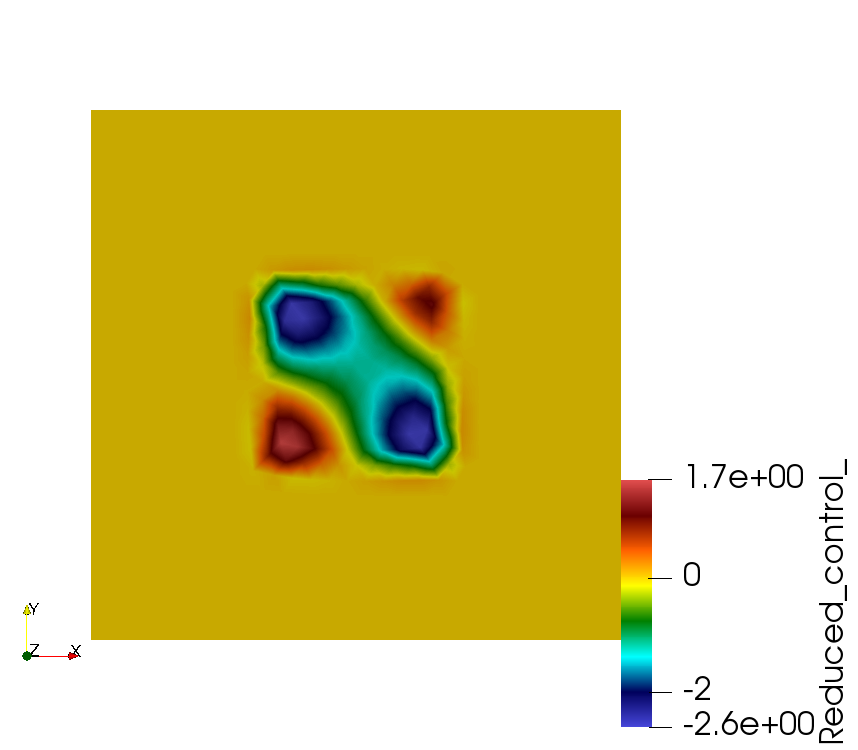}}}\hspace{5pt}
\subfloat[Error (control).]{%
\resizebox*{4.5cm}{!}{\includegraphics{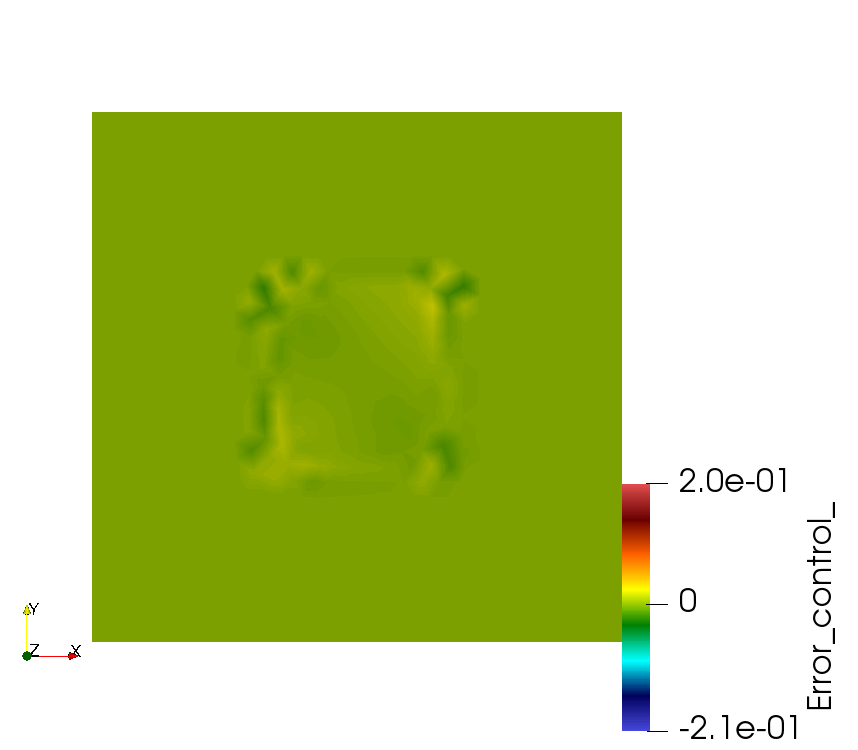}}}
\\
\subfloat[High--fidelity adjoint.]{%
\resizebox*{4.5cm}{!}{\includegraphics{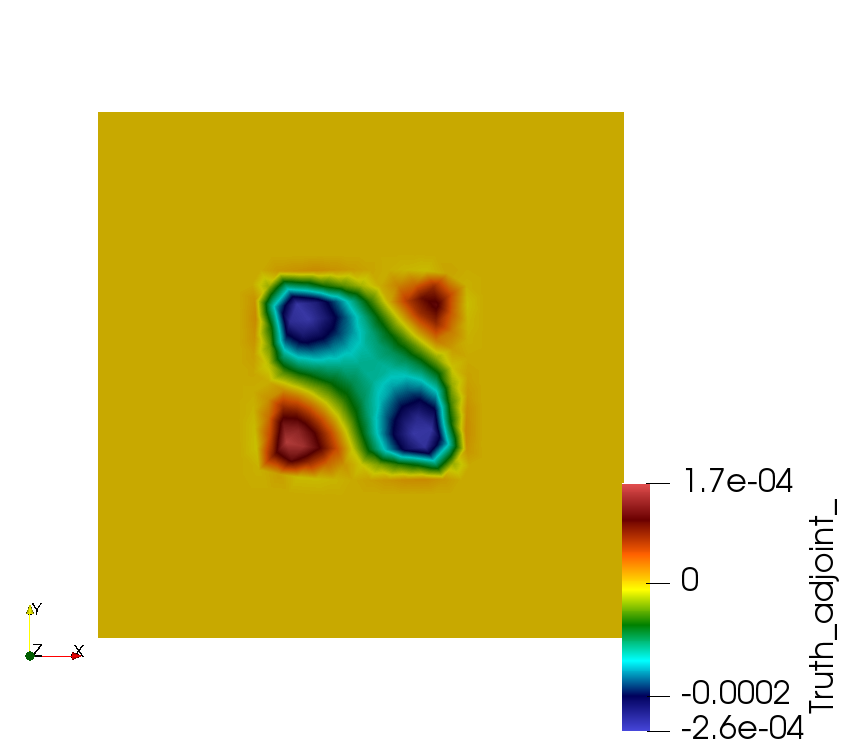}}}\hspace{5pt}
\subfloat[Reduced order adjoint.]{%
\resizebox*{4.5cm}{!}{\includegraphics{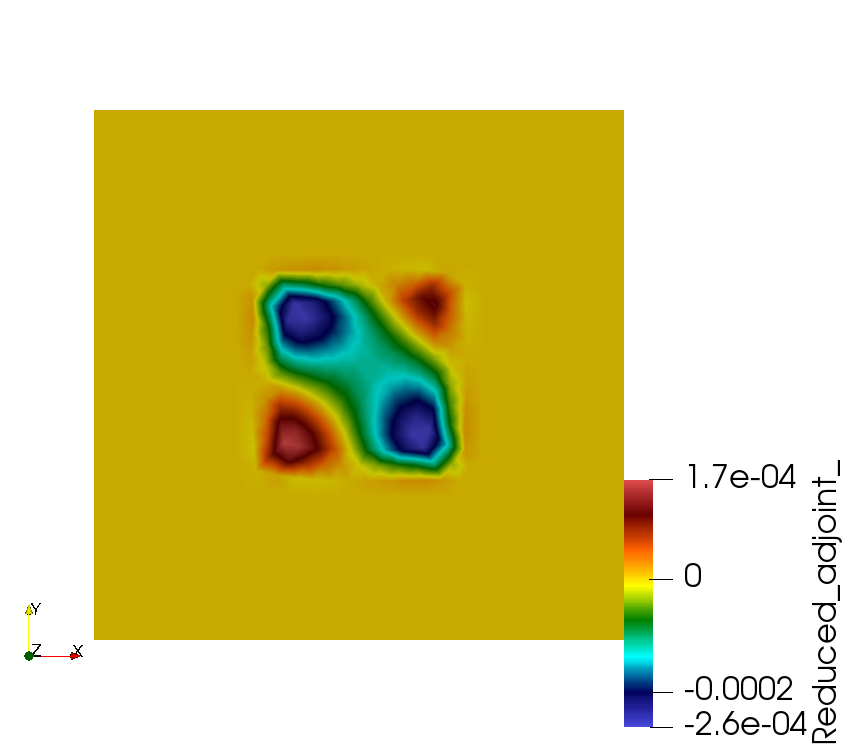}}}\hspace{5pt}
\subfloat[Error (adjoint).]{%
\resizebox*{4.5cm}{!}{\includegraphics{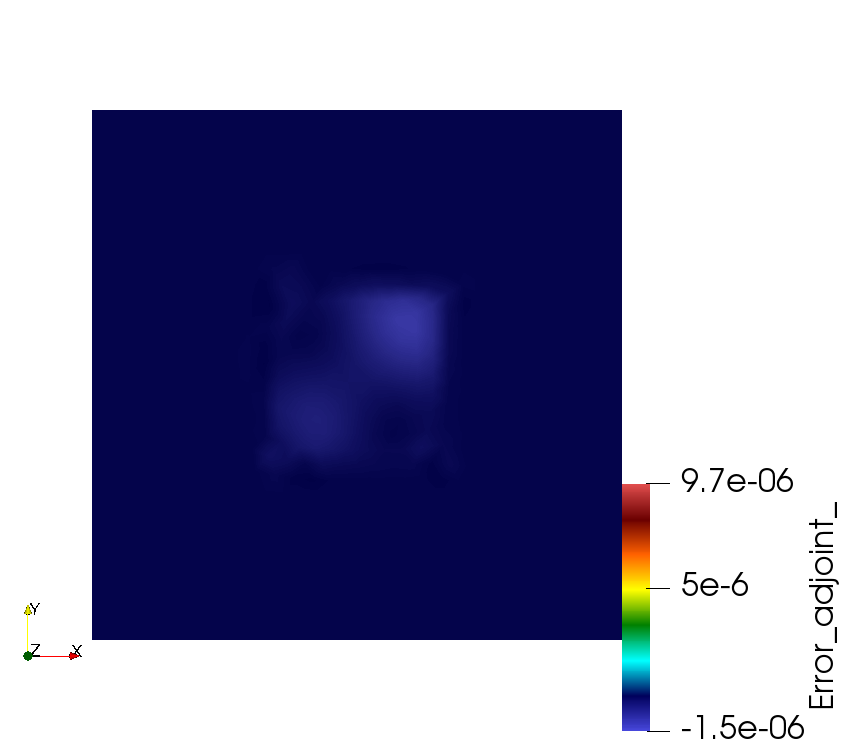}}}
\caption{High fidelity solutions (left), POD/DEIM approximations (middle), and their differences (right) for the OCP on the domain $\Omega(\boldsymbol\mu_0)$ boundary defined by the square centered at  $(1,1)$ and side length $2\times \boldsymbol\mu_0$, where $\boldsymbol\mu_0= 0.4757$.}\label{DEIM_vtk1a}
\end{figure}

In the following, we provide details regarding the computational performance of the method for the test case of the randomly chosen  parameter value $\boldsymbol\mu_0=0.4757$. Regrading DEIM approximations to the nonaffine components of the OCP system, the relevant $\left\|\cdot\right\|_2$--errors have been computed as follows:
\begin{itemize}
\item   DEIM approximation $\left\|\cdot\right\|_2$--error for $\mathbf{A}_{\boldsymbol\mu_0}$: $5.55e^{-07}$,
\item  DEIM approximation $\left\|\cdot\right\|_2$--error for $\mathbf{M}_{\boldsymbol\mu_0}$: $1.23e^{-05}$,
\item DEIM approximation $\left\|\cdot\right\|_2$--error for $\mathbf{b}_{\boldsymbol\mu_0}$: $1.95e^{-06}$,
\item  DEIM approximation $\left\|\cdot\right\|_2$--error for $\mathbf{c}_{\boldsymbol\mu_0}$: $1.11e^{-05}$.
\end{itemize}
Figure \ref{DEIM_vtk1a} illustrates the truth and ROM approximations for $\boldsymbol\mu_0=0.4757$, along with their differences. The relative errors between the truth and ROM solutions are as follows: \small
\begin{align*}
\frac{\left\|\mathbf{y}_{\boldsymbol\mu}^h-V_{yp}\mathbf{y}^N_{\boldsymbol\mu}\right\|_{\mathbf{M}_{\boldsymbol\mu}}}{\left\|\mathbf{y}_{\boldsymbol\mu}^h\right\|_{\mathbf{M}_{\boldsymbol\mu}}}= 0.00320, \ 
\frac{\left\|\mathbf{u}_{\boldsymbol\mu}^h-V_u\mathbf{u}^N_{\boldsymbol\mu}\right\|_{\mathbf{M}_{\boldsymbol\mu}}}{\left\|\mathbf{u}_{\boldsymbol\mu}^h\right\|_{\mathbf{M}_{\boldsymbol\mu}}}= 0.00414, \ 
\frac{\left\|\mathbf{p}_{\boldsymbol\mu}^h-V_{yp}\mathbf{p}^N_{\boldsymbol\mu}\right\|_{\mathbf{M}_{\boldsymbol\mu}}}{\left\|\mathbf{p}_{\boldsymbol\mu}^h\right\|_{\mathbf{M}_{\boldsymbol\mu}}}= 0.00424.
\end{align*} 
\normalsize
Regarding execution times, we remark that the full order model resolution for $\boldsymbol\mu_0=0.4757$ took 0.02881 seconds to complete, while the corresponding time for the ROM was 0.00211 seconds, leading to significant computational savings and a speed up factor of the order of $13.68$.

\begin{figure}
\centering
\subfloat[Error decay.]{%
\resizebox*{5.5cm}{!}{\includegraphics{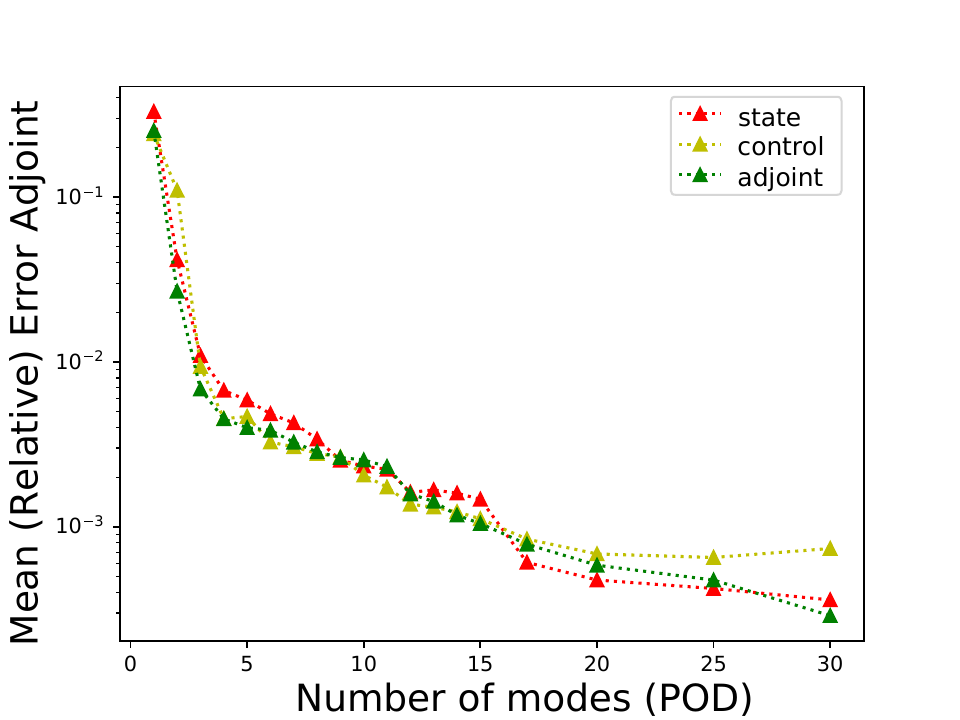}}}\hspace{5pt}
\subfloat[Resolution speedup.]{%
\resizebox*{5.5cm}{!}{\includegraphics{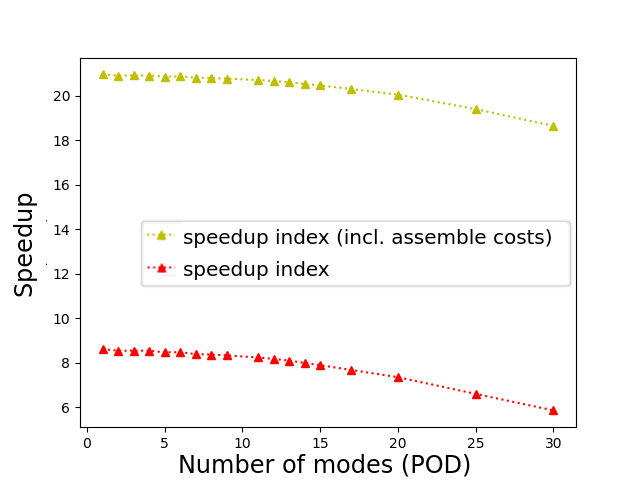}}}
   \caption{Reliability and efficiency of the POD/DEIM reduction for the optimal (distributed) control problem. \normalfont Mean relative errors between the truth and POD/DEIM--reduced approximations of  the  state/control/adjoint state, computed over a testing set of 30 randomly chosen parameters $\boldsymbol\mu \in (0.4,0.5)$, and resolution speedup (taking into account and disregarding assemble times), varying the number of retained POD modes.}
 \label{DEIMa2}\label{DEIMerrorsa}
\end{figure}

During an online phase, we assess the reliability of the POD/DEIM reduction by solving the full and reduced  order models for progressively increasing POD dimension. Average relative errors  between the high fidelity and POD/DEIM reduced approximations are depicted in   Figure \ref{DEIMerrorsa}. The error analysis has been carried out over a testing set of 30 randomly chosen parameters $\boldsymbol\mu \in [0.4,0.5]$ which has been selected during the on--line phase and is independent of the training set used for the construction of the solution snapshots. 
The convergence of the reduced solution with respect to the number of POD basis functions is verified. Figure \ref{DEIMerrorsa} suggests that we need less than $9$ basis functions for state, control and adjoint state, in order to correctly reproduce the high fidelity solution. Hence, instead of considering a ROM with order $2(N_y+N_p)+N_u=133$, we could have used a ROM with dimension $ 5\min\left\{N_y, N_u, N_p\right\}=95$ with even more competitive results. 

%

 More details can be found in Table \ref{poddeim_table}. Its first four columns are dedicated to reliability and error decay versus POD dimension; see also Figure \ref{DEIMa2}. Regarding efficiency, we mention that in compiling   column 5, we have considered three main contributions to overall solution time:
 \begin{enumerate}
 \item Time to construct the block system matrix $\mathcal{A}_N(\boldsymbol\mu)=V^T\mathcal{A}(\boldsymbol\mu)V$ and RHS block vector $\beta_N(\boldsymbol\mu)=V^T\beta(\boldsymbol\mu)$ in (\ref{ROM}) from their individual sub--blocks.
 \item Time to actually solve the reduced order system (\ref{ROM}).
 \item Time to retrieve the solution triple $(\mathbf{y}_{\boldsymbol\mu}, \mathbf{u}_{\boldsymbol\mu}, \mathbf{p}_{\boldsymbol\mu})$ from the solution vector $\begin{bmatrix}\mathbf{y}_{\boldsymbol\mu}^N \\ \mathbf{u}_{\boldsymbol\mu}^N \\ \mathbf{p}_{\boldsymbol\mu}^N \end{bmatrix}$.
 \end{enumerate}
 
 Then, in column 6, the speedup index has been computed, considering the relevant operations (1--3) for the full--order system (\ref{system}); for a graphical representation, refer to the red curve in Figure \ref{DEIMa2}. Hence, the results in columns 5/6 of Table  \ref{poddeim_table} disregard the time required to actually assemble the individual blocks of $\mathcal{A}_N(\boldsymbol\mu)$, $\beta_N(\boldsymbol\mu)$ and $\mathcal{A}(\boldsymbol\mu)$, $\beta(\boldsymbol\mu)$ of the reduced and full--order models, respectively. These are taken into account in columns 7/8, using the tools developed in subsections \ref{DEIMAnew}--\ref{DEIMMnew}. It turns out that the relative efficiency of DEIM in retrieving the components $\mathbf{A}_{\boldsymbol\mu}$, $\mathbf{M}_{\boldsymbol\mu}$,  $\mathbf{b}_{\boldsymbol\mu}$ and $\mathbf{c}_{\boldsymbol\mu}$ necessary to assemble $\mathcal{A}_N(\boldsymbol\mu)$ and $\beta_N(\boldsymbol\mu)$, as opposed to full-scale assemble, results in a significant improvement of speedup indices from column 6 for all POD dimensions. Resolution time improvement is visualized in Figure \ref{DEIMa2}.

\begin{table}[h!]
\begin{center}
\caption{Mean relative errors  between the truth and POD/DEIM-reduced approximations of $(y,u,p)$ and speed-up index over a testing set of 30 randomly chosen parameters $\mu \in   [0.4, 0.5]$, varying the number of retained POD modes.}\label{poddeim_table}
{\scriptsize
\begin{tabular}[scale=0.80]{cccc|cc|cc}
\hline
\# POD & Mean Rel.    & Mean Rel.  & Mean Rel.  & Solver  & Speedup & Solver  & Speedup\\
Modes &  Error $y$   &  Error  $u$ &  Error  $p$ &  timing &  index &  timing &  index\\
\hline
\hline
$1$ & $0.32975$ & $0.24152$    & 0.25266 & 0.00261& 8.62 & 0.00993& 20.96\\ 
$2$ & $0.04160$ & $0.10980$  & 0.02684&0.00263 & 8.54 & 0.00996& 20.91\\ 
$3$ & $0.01092$ & $0.00939$   & 0.00685& 0.00263 & 8.54 & 0.00996& 20.91\\ 
$4$ & $0.00672$ & $0.00453$   &0.00452 & 0.00263 & 8.54& 0.00996& 20.91\\ 
$5$ & $0.00587 $ & $0.00467$   &0.00399& 0.00266 & 8.47 & 0.00998& 20.86\\ 
$6$ & $0.00484$ & $0.00326$   & 0.00385& 0.00265 & 8.48 & 0.00998& 20.87\\ 
$7$ & $0.00427$ & $0.00306$   &0.00326 & 0.00268 & 8.39 & 0.01001& 20.81\\ 
$8$ & $0.00340$ & $0.00279$   & 0.00285&0.00269 & 8.37 & 0.01001& 20.79\\ 
$9$ & $0.00253$ & $0.00265$   & 0.00263& 0.00270 & 8.33 & 0.01002& 20.77\\ 
  $11$ &    $0.00224$ & $0.00174$ & 0.00232&0.00273 & 8.23 & 0.01005& 20.71\\
$12$ &    $0.00162$ & $0.00137$ & 0.00158& 0.00275 & 8.17 & 0.01008&  20.67\\
$13$ &    $0.00167$ & $0.00132$ &0.00143 & 0.00278 & 8.09 & 0.01010& 20.61\\
$14$ &    $0.00160$ & $0.00123$ & 0.00118&0.00282 & 7.98 & 0.01014&  20.53\\
$15$ &    $0.00147$ & $0.00111$ &0.00105& 0.00285 & 7.89 & 0.01018& 20.46\\
$17$ &    $0.00061$ & $0.00084$ &0.00078&  0.00293 & 7.67 & 0.01025& 20.31\\
$20$ &    $0.00048$ & $0.00068$ & 0.00059&0.00306 & 7.35 & 0.01038& 20.05\\
$25$ &   $0.00042$ & $0.00065$ &0.00048 & 0.00341& 6.58 & 0.01073& 19.40\\
$30$ &   $0.00036$ & $0.00074$ & 0.00029& 0.00384 & 5.85 & 0.01117& 18.64\\
\hline
\end{tabular}
}
\end{center}
\end{table}

\section{Conclusions and perspectives}

In this work, we have developed a reduced basis framework for the efficient solution of parametrized linear/quadratic optimal control problem governed by elliptic PDEs in  parametrically dependent domains. More precisely, we proposed a POD--Galerkin procedure to explore the solution manifold and to define a suitable low--dimensional reduced basis for the projection, using an aggregated space strategy. The basis functions are built via POD on optimal snapshots for all state, adjoint state and control  variables. Indeed, application of an embedded FEM as a truth solver allows these snapshots  to  be defined in a unified background mesh for different domain configurations. The well--posedeness of the model follows by exploiting a suitable saddle--point formulation. Application of DEIM for the system matrix and right--hand side vector allows the recovery of a full offline/online decomposition strategy, ensuring the online efficiency of the method. The main contribution of this work is the combination of POD Galerkin/DEIM procedure with an unfitted method to alleviate the problem of offline/online decoupling in the case of  parametrized domains. To the authors' best knowledge, this is the first time this has been attempted in the literature.

Our numerical results indicate the possibility of obtaining large computational savings in the online stage, in comparison to classical truth solvers.
The presented examples demonstrated how reduced order methods can be a useful tool in cases where parametrized simulations are typically very computationally demanding and costly.  A possible drawback resides in the offline stage, which demands large computational resources. Further developments on the present work would include time--dependent optimal control problems, considering also the non--linear case. 
{\blue{An important extension and future work is to investigate the efficiency of an unassembled form, namely UDEIM based on embedded FEMs and fixed background meshes, which may
lead to additional gains in the online cost of the reduced
order models although with additional costs in the offline stage similarly to  \cite{TiRi13, AHS14}.
%
%
}}

\section*{Acknowledgments}
 This project has received funding from the Hellenic Foundation for Research and Innovation (HFRI) and  the  General  Secretariat  for  Research  and  Technology (GSRT), under  grant agreement No[1115] (PI: E. Karatzas),  the ”First Call for H.F.R.I. Research Projects to support Faculty members and Researchers and the procurement of high-cost research equipment” grant 3270 and the support  of the National Infrastructures for Research and Technology S.A. (GRNET S.A.) in the National HPC facility - ARIS - under project ID pa190902.

\end{document}